\documentclass[12pt]{article}
\usepackage[margin=0.75in]{geometry}
\usepackage{amsfonts,amsmath,amssymb}
\usepackage{color}
\usepackage{lmodern}
\usepackage{textcomp} 
\usepackage{verbatim}
\usepackage{mathtools}
\usepackage{hyperref}
\usepackage{empheq}
\usepackage[normalem]{ulem}

\allowdisplaybreaks

\newcommand*\widefbox[1]{\fbox{\hspace{2em}#1\hspace{2em}}}

\newcommand{\bff}{\boldsymbol{f}}
\newcommand{\bj}{\boldsymbol{j}}
\newcommand{\bq}{\boldsymbol{q}}
\newcommand{\bu}{\boldsymbol{u}}
\newcommand{\bv}{\boldsymbol{v}}
\newcommand{\bw}{\boldsymbol{w}}
\newcommand{\bx}{\boldsymbol{x}}
\newcommand{\bz}{\boldsymbol{z}}
\newcommand{\bF}{\boldsymbol{F}}
\newcommand{\divx}{\mathrm{div}_{\boldsymbol{x}}}
\newcommand{\divq}{\mathrm{div}_{\boldsymbol{q}}}
\newcommand{\nabx}{\nabla_{\boldsymbol{x}}}
\newcommand{\nabq}{\nabla_{\boldsymbol{q}}}
\newcommand{\nabxq}{\nabla_{\bx, \bq}}

\newcommand{\them}{\theta_{\mathrm{min}}}
\newcommand{\bDD}{\mathbb{D}}
\newcommand{\bII}{\mathbb{I}}
\newcommand{\bSS}{\mathbb{S}}
\newcommand{\bTT}{\mathbb{T}}
\newcommand{\cU}{\mathcal{U}}
\newcommand{\dd}{\,\mathrm{d}}
\newcommand{\wtD}{\widetilde{D}}

\newtheorem{remark}{Remark}
\newtheorem{lemma}{Lemma}
\newtheorem{definition}{Definition}
\newtheorem{theorem}{Theorem}

\numberwithin{equation}{section}

\title{
Weak sequential stability of solutions to a nonisothermal kinetic model for incompressible dilute polymeric fluids}
\author{Miroslav Bul\'{\i}\v{c}ek\thanks{
Faculty of Mathematics and Physics, Charles University, Sokolovská 83, CZ 186 75 Praha 8–Karlín, Czech Republic. Email: \texttt{mbul8060@karlin.mff.cuni.cz}}, ~Josef M\'{a}lek\thanks{
Faculty of Mathematics and Physics, Charles University, Sokolovská 83, CZ 186 75 Praha 8–Karlín, Czech Republic. Email: \texttt{malek@karlin.mff.cuni.cz}}, ~and ~Endre S\"uli\thanks{
Mathematical Institute, University of Oxford, Radcliffe Observatory Quarter, Woodstock Road, Oxford OX2 6GG, UK. Email: \texttt{endre.suli@maths.ox.ac.uk}
\\
\textit{Acknowledgements:}  M. Bul\'{\i}\v{c}ek and J.
M\'{a}lek acknowledge the support of the project No. 25-16592S financed by the
Czech Science Foundation (GA\v{C}R). M. Bul\'{\i}\v{c}ek and J. M\'{a}lek
are members of the Ne\v{c}as Center for Mathematical Modeling. 
E. S\"uli is grateful to the Faculty of Mathematics and Physics of the Charles University and the Ne\v{c}as Center for Mathematical Modeling for hosting his visit during the period 7/9/2025--5/10/2025.}}
\date{\today}
\begin{document}
\maketitle
\begin{abstract}
The paper is concerned with the mathematical analysis of a class of thermodynamically consistent kinetic models for nonisothermal flows of dilute polymeric fluids, based on the identification of energy storage mechanisms and entropy production mechanisms in the fluid under consideration. The model involves a system of nonlinear partial differential equations coupling the unsteady incompressible temperature-dependent Navier--Stokes equations to a temperature-dependent generalization of the classical Fokker--Planck equation and an evolution equation for the absolute temperature. Sequences of smooth solutions to the initial-boundary-value problem, satisfying the available bounds that are uniform with respect to the given data of the model, are shown to converge to a global-in-time large-data weak solution that satisfies an energy inequality, where the absolute temperature satisfies a renormalized variational inequality, implying weak sequential stability of the mathematical model. 
\end{abstract}
\section{Introduction}\label{sec:1}

Since the pioneering contributions of Werner Kuhn \cite{Kuhn}, Hans Kramers \cite{Kramers} and other scientists working at the interface of polymer chemistry and statistical physics during the first half of the twentieth century, kinetic models have been widely and successfully used to describe the motion of polymeric fluids; see, in particular, the second volume of the two-volume monograph by Bird et al.~\cite{bird1987dynamics2}
and the book by \"Ottinger~\cite{MR1383323} for further details. During the past two decades significant progress has been made with the mathematical analysis of kinetic models of dilute polymers that involve the coupling of the incompressible or compressible Navier--Stokes equations to the Fokker--Planck equation, and proofs of the existence of large-data global weak solutions to these models are now available in various physical  settings; (see, for example, \cite{MR2660717,BS2011,MR2902849,MR3698145,MR3046297,MR4542472,MR4927782,MR4701225} as well as \cite{MR2221204,MR2039220,MR2456183,MR2843224,MR3010381,MR3916792,MR2454611} and the references cited therein for the case of coupled incompressible Navier--Stokes--Fokker--Planck systems, and 
\cite{MR3458248,MR3487270,MR3508492,MR4104467} for the mathematical analysis of coupled compressible Navier--Stokes--Fokker--Planck systems). While for fully macroscopic models of nonisothermal flows of  viscoelastic fluids there is now a growing body of mathematical theory focused on proofs of the existence of a global-in-time large-data weak solutions (cf., for example, \cite{MR4505743,MR4711550,MR1827894,MR3818666,MR4285769,MR2807424}; see also \cite{MR2729028,Malek}), the mathematical analysis of nonisothermal kinetic models of dilute polymers is in a much less satisfactory state. Indeed, all of the publications concerned with the mathematical analysis of coupled Navier--Stokes--Fokker--Planck systems cited above are restricted to isothermal flows: the temperature field is tacitly assumed to be homogeneous in space, and the evolution equation describing the temporal and spatial variations of the temperature field is excluded from the model.

The aim of this paper is to develop the mathematical analysis of a class of thermodynamically consistent kinetic models for nonisothermal flows of dilute polymeric fluids, formulated for compressible polymeric fluids by Dostal\'{\i}k et al.~in \cite{Dostalik}. While the model derived in~\cite{Dostalik} is essentially the same as the one discussed by Grmela and \"{O}ttinger in \cite{MR1492399,MR1492400},\footnote{See, in particular, Sec. III of \cite{MR1492400},  entitled \textit{Nonisothermal kinetic theory of polymeric fluids.}} its derivation is substantially different: instead of the so-called GENERIC formalism of Grmela and \"{O}ttinger, it is based on the identification of energy storage and entropy production mechanisms in the fluid under consideration, which then yields explicit formulae for the Cauchy stress tensor and for all the fluxes involved.  We shall confine our attention here to the physical situation where the solvent is an incompressible heat-conducting Newtonian fluid of constant specific density, resulting in the coupling of the unsteady incompressible temperature-dependent Navier--Stokes equations to a temperature-dependent generalization of the classical Fokker--Planck equation and an evolution equation for the absolute temperature. 
For clarity and completeness of the exposition, we shall revisit the derivation of the model presented in~\cite{Dostalik}, but in the case of an incompressible polymeric fluid, which is of interest to us here. Following the derivation of the model, we shall collect the equations into a system of nonlinear partial differential equations, supplement them with boundary and initial conditions, and state the precise conditions on the data, in preparation for the subsequent mathematical analysis of the model. The main result of the paper is the proof of weak sequential stability of solutions to the mathematical model in the sense that sequences of smooth solutions to the initial-boundary-value problem, satisfying the available bounds that are uniform with respect to the given data, converge to a global-in-time weak solution that satisfies an energy inequality, where the absolute temperature satisfies a renormalized variational inequality. While the result is not an explicit proof of the existence of weak solutions, it is a key step towards the development of a rigorous existence theory for the model. The next section focuses on the derivation of the model. The structure of the rest of the paper is outlined at the end of Section \ref{sec:1a}.

\section{Derivation of the model}\label{sec:1a}

\smallskip
The aim of this subsection is to revisit the derivation of the model presented in \cite{Dostalik}; however, instead of the compressible heat-conducting dilute polymeric fluid flow model considered there, we shall focus on the case of an incompressible heat-conducting dilute polymeric fluid. 
Since the precise values of the various physical constants
included in the model stated in \cite{Dostalik} do not influence the mathematical analysis pursued later on in this paper, we shall set all of them (viz. the specific density $\rho_s>0$, the Boltzmann constant $k_{\mathrm B}>0$, the constant reference temperature $\theta_{\mathrm{ref}}>0$, the specific heat at constant volume $c_{\rm V}>0$,
and the reference polymer molecule length $q_{\mathrm{ref}}>0$) equal to $1$, with the exception of the (constant) hydrodynamic drag coefficient $\zeta$, which we set equal to $1/2$. 

Suppose that $T>0$ and let $\Omega$ be a bounded open simply-connected Lipschitz domain in $\mathbb{R}^d$, where $d \in \{2,3\}$.
For our state variables, we choose the velocity $\bv$, the internal energy $e$, and the nonnegative probability density function $\varphi$. We have the following governing equations: 
\begin{alignat}{2}
\begin{aligned}\label{eq:nond}
\mbox{div}_{\bx}\,\bv &= 0  &&\qquad \mbox{in $(0,T) \times \Omega$},\\
\dot{\bv}&= \mbox{div}_{\bx}\, \mathbb{T} + \bff  &&\qquad \mbox{in $(0,T) \times \Omega$},\\
\dot{e} &= \mathbb{T} : \mathbb{D}(\bv)- \mbox{div}_{\bx}\, \bj_e  &&\qquad \mbox{in $(0,T) \times \Omega$},\\
\dot\varphi + \mbox{div}_{\bx} \bj_{\varphi,\bx} + \mbox{div}_{\bq}\big((\nabx \bv)\bq \varphi + \bj_{\varphi,\bq}\big)&=0  &&\qquad \mbox{in $(0,T) \times \Omega \times D$},
\end{aligned}
\end{alignat}
where $\dot{\bv}$, $\dot{e}$ and $\dot{\varphi}$ signify 
the material derivatives\footnote{The material-derivative $\dot{u}$ of a real-valued function $u=u(t,\bx)$ is defined by $\dot{u}= \partial_t u + \bv \cdot \nabx u$, and the material derivative $\dot \bu$ of an $\mathbb{R}^d$-valued vector function $\bu=\bu(t,\bx)$ is defined by $\dot{\bu} = \partial_t {\bu} + (\bv \cdot \nabx) \bu$. In particular, when $\divx\,\bv=0$ as is the case here, $\dot{u} = \partial_t u  + \divx(u\bv)$ and $\dot{\bv} = \partial_t {\bv} + \divx (\bv \otimes\bv)$, where $\bv \otimes \bv \in \mathbb{R}^{d \times d}$ with $(\bv \otimes \bv)_{i,j}:= v_i v_j$.} of $\bv$, $e$ and $\varphi$, respectively,  $\bff$ is the given density of body forces, $\bTT$ is the Cauchy stress tensor, $\bj_e$, $\bj_{\varphi, \bx}$ and $\bj_{\varphi,\bq}$ are flux-vectors, to be determined in the course of the derivation of the model, and $\mathbb{D}(\bv):=\frac{1}{2}(\nabx \bv + (\nabx \bv)^{\mathrm{T}})$ is the symmetric velocity gradient. The partial differential equation  \eqref{eq:nond}$_4$ 
is the temperature-dependent generalization of the Fokker--Planck equation for the evolution of the probability density function $\varphi=\varphi(t,\bx,\bq)$, modelling, at time $t \in [0,T)$ and spatial location $\bx \in \Omega$, the random configuration (i.e., orientation vector) $\bq$ of polymer molecules idealized as elastic dumbbells suspended in the solvent (which are assumed not to interact with each other and which move without self-interaction). In eq. \eqref{eq:nond}$_4$, $D:=B(\boldsymbol{0},\sqrt{b}) \subset \mathbb{R}^d$ is a ball of radius $\sqrt{b}$ centred at the origin, and $b>0$.\footnote{The extension of the mathematical analysis of the dumbbell model considered here to a bead-spring-chain model, involving $K+1$ massless beads instead of just a pair of massless beads, which are linearly coupled with $K$ elastic springs, does not involve any conceptual or technical difficulties. For the sake of clarity of the exposition we shall therefore confine our attention to the dumbbell model, corresponding to the case of $K=1$.}  The form of the partial differential equation \eqref{eq:nond} is motivated by the structure of the Fokker--Planck equation in the isothermal setting (see, for example,  Sec.~1.1.4 of the doctoral thesis of Lozinski \cite{lozinski} and the introductory section of \cite{MR2338493} for the derivation of the Fokker--Planck equation associated with the dumbbell model from Brownian dynamics). 

The spring-force associated with the dumbbell is assumed to be of the form
\begin{align}\label{eq:F} 
\boldsymbol{F}(\bq):= \nabq 
\left[\, \cU_e(\bq) +  \theta \,\cU_\eta(\bq)\right]\!,
\end{align}
where $\theta>0$ is the absolute temperature. The function $\mathcal{U}$ defined by 
\[
\mathcal{U}(\bq):= \cU_e (\bq) + \theta\,\cU_\eta(\bq)
\]
is the total spring potential of the dumbbell that represents the polymer molecule suspended in the solvent, which is additively decomposed into its elastic part $\mathcal{U}_e$
 and its entropic part $\theta \,\mathcal{U}_\eta$.
 
In what follows, we shall adopt the notational convention that the same symbol, $|\cdot|$, denotes the absolute value of a real number, the Euclidean norm of a $d$-component vector, or the Frobenius norm of a $d \times d$ matrix; it will be clear from the context which of these three interpretations is intended in a particular case. The functions $\cU_e$ and $\cU_\eta$ are defined by
\begin{align}\label{eq:4a} 
\cU_e(\bq):= U_e\bigg(\frac{1}{2}\left|\bq \right|^2\bigg)\quad \mbox{for $\bq \in  D$} \quad \mbox{and}\quad \cU_\eta(\bq):= U_\eta\bigg(\frac{1}{2}\left|\bq \right|^2\bigg)\quad \mbox{for $\bq \in \overline{D}$},
\end{align}
where $\overline{D}$ is the closure of $D$, and $U_e \in C^2([0,b/2); \mathbb{R}_{\geq 0}) \cap L^1((0,b/2))$ and $U_\eta \in C^1([0,b/2];\mathbb{R}_{\geq 0})$ are monotonically increasing functions. We shall make further, physically relevant assumptions on $U_e$ and $U_\eta$ later on, but for now the above ones will suffice. We note in passing that because 
\[\boldsymbol{F} \otimes \bq  = \left[U_e'\left(\frac{1}{2}|\bq|^2\right)+\theta \,U_\eta'\left(\frac{1}{2}|\bq|^2\right)\right](\bq \otimes \bq),
\]
we find that $\boldsymbol{F} \otimes \bq$ is a nonnegative scalar multiple of the symmetric positive semidefinite rank-one matrix $\bq \otimes \bq$. Therefore, $ \boldsymbol{F} \otimes \bq = (\boldsymbol{F} \otimes \bq)^{\mathrm{T}} = \bq \otimes \boldsymbol{F}$ is a symmetric positive semidefinite rank-one matrix.

Our objective is to determine the Cauchy stress tensor $\mathbb{T}$ and the fluxes $\bj_e$, $\bj_{\varphi,\bx}$ and $\bj_{\varphi,\bq}$
in the model \eqref{eq:nond}
from the specified constitutive equation for the Helmholtz free energy, $\psi= \psi(\theta,\varphi)$, defined by 
\begin{align}\label{eq:helm}
\psi(\theta,\varphi):= - \theta(\log \theta - 1) + \int_D \cU_e(\bq)\varphi \dd \bq  + \theta \int_D [\,\cU_\eta(\bq) \varphi + \varphi \log \varphi] \dd \bq.
\end{align}
In order to fulfil this objective, we shall derive an evolution equation for the specific entropy $\eta$ of the form
\begin{align} \label{eq:ent0}
\dot{\eta} + \mbox{div}_{\bx}\, \bj_\eta = \xi,
\end{align}
where $\bj_\eta$ is the entropy-flux, and require that $\xi \geq 0$, in agreement with the second law of thermodynamics. In the process of doing so, we shall also specify both $\bj_\eta$ and $\xi$. 

We begin by noting that while both sides of the equality \eqref{eq:helm} are functions of $t$ and $\bx$, the Helmholtz free energy $\psi$ can also be understood as a functional that maps the pair of functions $\theta = \theta(t,\bx)$ and $\varphi = \varphi(t,\bx,\bq)$ to $\psi(\theta,\varphi)$, parametrized by $t$ and $\bx$.

The thermodynamic quantities $\psi$, $\theta$, $e$, and $\eta$ are related by means of the relationships
\begin{align} \label{eq:ther}
\psi:= e -  \eta \theta \qquad \mbox{and}\qquad \frac{\partial \psi}{\partial \theta} = - \eta,
\end{align}
where $\eta= \eta(t,\bx)$ is the specific entropy and $\partial \psi/\partial \theta$ denotes the partial Gateaux derivative of $\psi$, understood as a functional, with respect to $\theta$. It then follows from \eqref{eq:ther}$_2$ that
\begin{align}\label{eq:eta}
\eta := - \frac{\partial \psi}{\partial \theta} = \log \theta - \int_D [\,\cU_\eta(\bq) \varphi  + \varphi \log \varphi ]\dd \bq.
\end{align}
Therefore, by \eqref{eq:ther}$_1$,  \eqref{eq:helm} and \eqref{eq:eta}, the internal energy $e$ can be expressed as
\begin{align}
\label{eq:e} 
e= \psi + \eta \theta =\psi - \frac{\partial \psi}{\partial \theta} \theta = \theta + \int_D \cU_e(\bq) \varphi \dd \bq.
\end{align}

The total energy is
\begin{align}\label{eq:E+} 
E:=  \frac{1}{2}|\bv|^2 + e = \frac{1}{2}|\bv|^2 + \theta + \int_D \cU_e(\bq) \varphi \dd \bq.
\end{align}
Taking the scalar product of \eqref{eq:nond} with $\bv$, adding the resulting equality to \eqref{eq:nond}$_3$ and noting \eqref{eq:E+}, it follows that the evolution of the total energy $E$ is governed by the partial differential equation
\begin{align}\label{eq:E} \dot{E} + \divx(\bj_e - \mathbb{T}\bv) = \bff \cdot \bv.
\end{align}
If, as will be the case here, the equation \eqref{eq:E} is supplemented with the boundary condition
\begin{align}\label{eq:Ebc}
(E \bv + \bj_e - \bTT \bv)\cdot \boldsymbol{n}_{\bx} = 0 \quad \mbox{on $(0,T) \times \Omega$,}    
\end{align}
where $\boldsymbol{n}_{\bx}$ is the unit outward normal vector to the boundary $\partial \Omega$ of $\Omega$, then \eqref{eq:nond}$_1$, \eqref{eq:E} and \eqref{eq:Ebc} imply that 
\begin{align}\label{eq:en}
\frac{\dd}{\dd t}\int_\Omega E(t,\bx) \dd \bx  = \int_\Omega \bff(t,\bx) \cdot \bv(t,\bx) \dd \bx,
\end{align}
expressing the balance of the total energy. 

Next, we derive an evolution equation for the specific entropy $\eta$. Taking the material derivative of both sides of the equality \eqref{eq:ther}$_1$, using the chain rule on the left-hand side and the product rule on the right-hand side, we have
\begin{align}\label{eq:ther1}\frac{\partial \psi}{\partial \theta}\dot{\theta} + \frac{\partial \psi}{\partial \varphi}\dot\varphi = \dot{e} - \theta \dot \eta - \eta \dot \theta,
\end{align}
where $(\partial \psi/\partial \theta)\dot\theta$ is the partial Gateaux derivative of $\psi$ with respect to $\theta$ acting on $\dot\theta$, while $(\partial \psi/\partial \varphi)\dot\varphi$ is the partial Gateaux derivative of $\psi$ with respect to $\varphi$ acting on $\dot\varphi$. However, since $(\partial\psi/\partial \theta)\dot\theta = - \eta\dot \theta$ thanks to \eqref{eq:ther}$_2$, the equality \eqref{eq:ther1} is simplified to 
\[ \theta \dot\eta = \dot e - \frac{\partial \psi}{\partial \varphi} \dot \varphi.\]
Next, because 
\[ \frac{\partial \psi}{\partial \varphi} \dot\varphi = \int_D \cU_e(\bq)\dot \varphi  \dd \bq + \theta \int_D [\cU_\eta(\bq) + \log \varphi + 1] \dot \varphi \dd \bq, \]
by inserting this into the above equality and using the equations for $\dot{e}$ and $\dot\varphi$, we find that
\begin{align*}
\theta \dot\eta &= \mathbb{T} : \mathbb{D}(\bv)- \mbox{div}_{\bx}\, \bj_e  + \int_D \big[\cU_e(\bq) + \theta \cU_\eta(\bq) + \theta(\log \varphi + 1)\big]\,\mbox{div}_{\bx} \bj_{\varphi,\bx}\dd \bq    \\
&\quad +  \int_D \big[\cU_e(\bq) + \theta \cU_\eta(\bq) + \theta(\log \varphi + 1)\big] \,\mbox{div}_{\bq}\big((\nabx \bv)\bq \varphi + \bj_{\varphi,\bq}\big) \dd \bq\\
&=  \mathbb{T} : \mathbb{D}(\bv)- \mbox{div}_{\bx} \left[ \bj_e  - \int_D \big[\cU_e(\bq) + \theta \cU_\eta(\bq) + \theta(\log \varphi + 1)\big] \bj_{\varphi,\bx}\dd \bq \right] \\
&\quad -\int_D  \nabx \big[\cU_e(\bq) + \theta \cU_\eta(\bq) + \theta(\log \varphi + 1)\big]\cdot  \bj_{\varphi,\bx}\dd \bq\\
&\quad - \int_D \nabq \big[\cU_e(\bq) + \theta \cU_\eta(\bq) + \theta(\log \varphi + 1)\big] \cdot\big((\nabx \bv)\bq \varphi + \bj_{\varphi,\bq}\big) \dd \bq\\
&\quad +  \int_{\partial D} \big[\cU_e(\bq) + \theta \cU_\eta(\bq) + \theta(\log \varphi + 1)\big] \,\big((\nabx \bv)\bq \varphi + \bj_{\varphi,\bq}\big)\cdot \boldsymbol{n}_{\bq} \dd S(\bq),
\end{align*}
following partial integration on $D$, where $\boldsymbol{n}_{\bq}:=\bq/|\bq|$ with $\bq \in D$ is the unit outward normal vector to the boundary $\partial D$ of $D$.
We proceed by imposing the boundary condition
\begin{align}\label{eq:bcfi1} \big((\nabx \bv)\bq \varphi + \bj_{\varphi,\bq}\big)\cdot \boldsymbol{n}_{\bq} = 0 \quad \mbox{on $(0,T) \times \Omega \times \partial D$},
\end{align}
which then annihilates the boundary integral over $\partial D$. Hence, 
using \eqref{eq:F} and the symmetry of the matrix-valued function $\boldsymbol{F} \otimes \bq$,
\begin{align*}
\theta \dot\eta &=  \mathbb{T} : \mathbb{D}(\bv)- \mbox{div}_{\bx} \left[ \bj_e  - \int_D \big[\cU_e(\bq) + \theta \cU_\eta(\bq) + \theta(\log \varphi + 1)\big] \bj_{\varphi,\bx}\dd \bq \right] \\
&\quad -\int_D \left\{ \big[ \cU_\eta(\bq) + (\log \varphi + 1)\big]\,\nabx \theta  + \theta \frac{\nabx \varphi}{\varphi}\right\} \cdot  \bj_{\varphi,\bx}\dd \bq\\
&\quad - \int_D \nabq \big[\cU_e(\bq) + \theta \cU_\eta(\bq) + \theta(\log \varphi + 1)\big] \cdot\big((\nabx \bv)\bq \varphi + \bj_{\varphi,\bq}\big) \dd \bq\\
&=  \mathbb{T} : \mathbb{D}(\bv)- \mbox{div}_{\bx} \left[ \bj_e  - \int_D \big[\cU_e(\bq) + \theta \cU_\eta(\bq) + \theta(\log \varphi + 1)\big] \bj_{\varphi,\bx}\dd \bq \right] \\
&\quad -\int_D \left\{ \big[ \cU_\eta(\bq) + (\log \varphi + 1)\big]\,\nabx \theta  + \theta \frac{\nabx \varphi}{\varphi}\right\} \cdot  \bj_{\varphi,\bx}\dd \bq\\
&\quad - \int_D (\boldsymbol{F} \otimes \bq)\varphi \dd \bq  : \mathbb{D}(\bv) -  \theta \int_D \nabq \varphi \cdot (\nabx \bv)\bq \dd \bq\\
&\quad - \int_D \big[\boldsymbol{F}+ \theta\frac{\nabq \varphi}{\varphi}\big] \cdot \bj_{\varphi,\bq} \dd \bq.
\end{align*}
We focus our attention on the second integral in the penultimate line. Because $\divq ((\nabx \bv)\bq) = 0$, 
partial integration over $D$ implies that 
\[ \int_D \nabq \varphi \cdot (\nabx \bv)\bq \dd \bq = \int_{\partial D} \varphi \big((\nabx \bv)\bq \cdot \boldsymbol{n}_{\bq}\big) \dd S(\bq).\]
For reasons that will become apparent after we have stated our precise assumptions on the spring potential $\cU_e$, the probability density function $\varphi$ must vanish on $(0,T) \times \Omega \times \partial D$. Since we have already imposed a boundary condition for $\varphi$ on this set (cf.~\eqref{eq:bcfi1} above), we cannot of course explicitly impose $\varphi = 0$ on $(0,T) \times \Omega \times \partial D$ as an additional boundary condition. Interestingly, the reason why $\varphi$ should vanish on $(0,T) \times \Omega \times \partial D$ is quite different and more indirect. Suffice to say at this point that it is a consequence of the fact that $\cU_e$ will be assumed to be a finitely extensible nonlinear elastic (FENE) type spring potential, which blows up on $\partial D$, and then the requirement that a term of the form $\int_{[0,T] \times \Omega \times D}|\nabq\,\cU_e(\bq)|^2 \varphi(t,\bx,\bq) \dd \bq \dd \bx \dd t$ that appears in the energy inequality associated with the model should remain finite is, in fact, the mechanism that forces $\varphi$ to vanish on $(0,T) \times \Omega \times \partial D$; the precise details of this will be given in Section \ref{sec:2}. In any case, as $\varphi$ vanishes on this set, the boundary integral stated above is equal to $0$. Consequently, 
\begin{align*}
\theta \dot \eta &=  \left[\mathbb{T} - \int_D (\boldsymbol{F} \otimes \bq)\varphi \dd \bq \right] : \mathbb{D}(\bv) - \mbox{div}_{\bx} \left[ \bj_e  - \int_D \big[\cU_e(\bq) + \theta \cU_\eta(\bq) + \theta(\log \varphi + 1)\big] \bj_{\varphi,\bx}\dd \bq \right] \\
&\quad -\int_D \left\{ \big[ \cU_\eta(\bq) + (\log \varphi + 1)\big]\,\nabx \theta  + \theta \frac{\nabx \varphi}{\varphi}\right\} \cdot  \bj_{\varphi,\bx}\dd \bq\\
&\quad - \int_D \big[\boldsymbol{F}+ \theta\frac{\nabq \varphi}{\varphi}\big] \cdot \bj_{\varphi,\bq} \dd \bq.
\end{align*}
Next, we divide this equality by $\theta$ with the aim of transforming it into an equation of the desired form \eqref{eq:ent0}, while defining $\mathbb{T}$, $\bj_e$, $\bj_{\varphi,\bx}$ and $\bj_{\varphi,\bq}$ so that the term $\xi$ appearing on the right-hand side of \eqref{eq:ent0} is nonnegative, as required. Thus, we have
\begin{align*}
\dot\eta &=   
\frac{1}{\theta}\left[\mathbb{T} - \int_D (\boldsymbol{F} \otimes \bq)\varphi \dd \bq \right] : \mathbb{D}(\bv)
- \frac{1}{\theta}\mbox{div}_{\bx} \left[ \bj_e  - \int_D \big[\cU_e(\bq) + \theta \cU_\eta(\bq) + \theta(\log \varphi + 1)\big] \bj_{\varphi,\bx}\dd \bq \right] \\
&\quad -\int_D \left\{ \big[ \cU_\eta(\bq) + (\log \varphi + 1)\big]\,\frac{\nabx \theta}{\theta}  + \frac{\nabx \varphi}{\varphi}\right\} \cdot  \bj_{\varphi,\bx}\dd \bq\\
&\quad - \int_D \big[\frac{1}{\theta}\boldsymbol{F}+ \frac{\nabq \varphi}{\varphi}\big] \cdot \bj_{\varphi,\bq} \dd \bq\\
&= \frac{1}{\theta}\left[\mathbb{T} - \int_D (\boldsymbol{F} \otimes \bq)\varphi \dd \bq \right] : \mathbb{D}(\bv)
- \mbox{div}_{\bx} \left[ \frac{\bj_e}{\theta}  - \frac{1}{\theta}\int_D \big[\cU_e(\bq) + \theta \cU_\eta(\bq) + \theta(\log \varphi + 1)\big] \bj_{\varphi,\bx}\dd \bq \right] \\
& \quad - \frac{\nabx \theta}{\theta^2} \cdot  \left[ \bj_e  - \int_D \big[\cU_e(\bq) + \theta \cU_\eta(\bq) + \theta(\log \varphi + 1)\big] \bj_{\varphi,\bx}\dd \bq \right]\\
&\quad -\int_D \left\{ \big[ \cU_\eta(\bq) + (\log \varphi + 1)\big]\,\frac{\nabx \theta}{\theta}  + \frac{\nabx \varphi}{\varphi}\right\} \cdot  \bj_{\varphi,\bx}\dd \bq\\
&\quad - \int_D \big[\frac{1}{\theta}\boldsymbol{F}+ \frac{\nabq \varphi}{\varphi}\big] \cdot \bj_{\varphi,\bq} \dd \bq.
\end{align*}
Hence, by transferring the second term appearing on the right-hand side to the left-hand side, 
\begin{align}\label{eq:et}
&\dot\eta + \mbox{div}_{\bx} \left[ \frac{\bj_e}{\theta}  - \frac{1}{\theta}\int_D \big[\cU_e(\bq) + \theta \cU_\eta(\bq) + \theta(\log \varphi + 1)\big] \bj_{\varphi,\bx}\dd \bq \right]\nonumber\\
&= 
\frac{1}{\theta}\left[\mathbb{T} - \int_D (\boldsymbol{F} \otimes \bq)\varphi \dd \bq \right] : \mathbb{D}(\bv)
 - \frac{\nabx \theta}{\theta^2} \cdot  \left[ \bj_e  - \int_D \big[\cU_e(\bq) + \theta \cU_\eta(\bq) + \theta(\log \varphi + 1)\big] \bj_{\varphi,\bx}\dd \bq \right] \nonumber\\
&\quad -\int_D \left\{ \big[ \cU_\eta(\bq) + (\log \varphi + 1)\big]\,\frac{\nabx \theta}{\theta}  + \frac{\nabx \varphi}{\varphi}\right\} \cdot  \bj_{\varphi,\bx}\dd \bq \nonumber\\
&\quad - \int_D \big[\frac{1}{\theta}\boldsymbol{F}+ \frac{\nabq \varphi}{\varphi}\big] \cdot \bj_{\varphi,\bq} \dd \bq \nonumber\\
& = 
\frac{1}{\theta}\left[\mathbb{T} - \int_D (\boldsymbol{F} \otimes \bq)\varphi \dd \bq \right] : \mathbb{D}(\bv) - \frac{\nabx \theta}{\theta^2} \cdot  \left[ \bj_e  - \int_D\cU_e(\bq) 
 \bj_{\varphi,\bx}\dd \bq \right] \nonumber\\
&\quad -\frac{1}{\theta}\int_D \frac{\nabx \varphi}{\varphi} \cdot  \bj_{\varphi,\bx}\dd \bq  - \int_D \big[\frac{1}{\theta}\boldsymbol{F}+ \frac{\nabq \varphi}{\varphi}\big] \cdot \bj_{\varphi,\bq} \dd \bq.
\end{align}
The left-hand side of the equality \eqref{eq:et} is now in the desired form. Therefore, it remains to choose
$\bTT$,  $\bj_e$, $\bj_{\varphi,\bx}$ and $\bj_{\varphi,\bq}$ so that each of the expressions in the last two lines of \eqref{eq:et} is nonnegative. Starting from the last term in the last line and proceeding backwards, we first define
\begin{align}\label{eq:7}
\bj_{\varphi,\bq}:= - 4\boldsymbol{F} \varphi - 4 \theta \nabq \varphi.
\end{align}
Next, we define
\begin{align}\label{eq:7+}
\bj_{\varphi,\bx}:= - \theta \nabx \varphi,
\end{align}
and then we set 
\begin{align}\label{eq:je} 
\bj_e - \int_D \cU_e (\bq) \bj_{\varphi,\bx}= -\kappa(\theta) \nabx \theta,
\end{align}
where $\kappa(\theta)>0$ is the temperature-dependent coefficient of heat conductivity; therefore, 
\[ \bj_e = -\kappa(\theta) \nabx\theta - \theta \nabx \int_D \cU_e(\bq) \varphi \dd \bq.\]
Concerning the first term on the right-hand side of \eqref{eq:et}, because $\bII : \bDD(\bv) = \mbox{div}\, \bv = 0$, 
the addition of a scalar multiple of the identity matrix $\bII \in \mathbb{R}^{d \times d}$ to the expression in the square brackets has no influence on the value of the right-hand side of the equality. Bearing this in mind, we recall that the deviatoric (traceless) part $\bSS_\delta$ of a matrix $\bSS \in \mathbb{R}^{d \times d}$ is defined by $\bSS_\delta := \bSS - \frac{1}{d} \mathrm{tr}(\bSS) \bII$, and we set 
\[ \left[\mathbb{T} - \int_D (\boldsymbol{F} \otimes \bq)\varphi \dd \bq \right]_\delta = 2 \nu(\theta) \bDD(\bv),
\] 
where $\nu(\theta)>0$ is the temperature-dependent shear viscosity coefficient. Consequently, 
\[2 \nu(\theta) \bDD(\bv) = \mathbb{T} - \int_D (\boldsymbol{F} \otimes \bq)\varphi \dd \bq  - \frac{1}{d} \mathrm{tr}\left(\mathbb{T} - \int_D (\boldsymbol{F} \otimes \bq)\varphi \dd \bq \right)\bII; \]
that is, 
\[ \bTT =  \frac{1}{d} \mathrm{tr}\left(\mathbb{T} - \int_D (\boldsymbol{F} \otimes \bq)\varphi \dd \bq \right) \bII + 2 \nu(\theta) \bDD(\bv) + \int_D (\boldsymbol{F} \otimes \bq)\varphi \dd \bq.\]
Motivated by this form, we define the Cauchy stress tensor $\bTT$ by
\begin{align*} 
\bTT := - \mathrm{p} \bII+ 2\nu(\theta) \bDD(\bv) - 2 k_\mathrm{B} \theta n_{\mathrm{P}} \bII + \int_D (\boldsymbol{F} \otimes \bq) \varphi  \dd \bq,
\end{align*}
where the sum of the first two terms is the usual Navier--Stokes--Fourier viscous stress tensor, and the sum of the last two terms is Kramers' polymeric stress tensor; $\mathrm{p}=\mathrm{p}(t,\bx)$ is the pressure (i.e., the constitutively undetermined part of the spherical stress), $k_{\mathrm{B}}>0$ is the Boltzmann constant, and $n_{\mathrm{P}}=n_{\mathrm{P}}(t,\bx)$ is the polymer number density, defined by
\begin{align}\label{eq:np} 
n_{\mathrm{P}}(t,\bx):=\int_D \varphi(t,\bx,\bq) \dd \bq.
\end{align}
Scaling the Boltzmann constant to 1, we therefore (re)define
\begin{align}\label{eq:T0} 
\bTT := - \mathrm{p} \bII+ 2\nu(\theta) \bDD(\bv) - 2 \theta n_{\mathrm{P}} \bII + \int_D (\boldsymbol{F} \otimes \bq) \varphi  \dd \bq.
\end{align}
Finally, in view of the second term on the left-hand side of \eqref{eq:et}, we define the entropy flux $\bj_\eta$ as
\begin{align} \label{eq:entf}
\bj_\eta:= \frac{\bj_e}{\theta}  - \frac{1}{\theta}\int_D \big[\cU_e(\bq) + \theta \cU_\eta(\bq) + \theta(\log \varphi + 1)\big] \bj_{\varphi,\bx}\dd \bq  = - \frac{\kappa(\theta) \nabx \varphi}{\theta} + \theta \nabx \int_D [\cU_\eta \varphi +  \varphi \log \varphi] \dd  \bq. 
\end{align}
Thereby we arrive at the desired form of the evolution equation for the specific entropy:
\begin{align}\label{eq:entr}
\begin{aligned}
&\partial_t \eta + \divx(\bv \eta) + \divx\left(\theta \nabx \int_D \cU_\eta \varphi \dd \bq + \theta \nabx \int_D \varphi \log \varphi \dd \bq \right) - \divx\left(\frac{\kappa(\theta) \nabx \theta}{\theta}\right)\\
&\quad = \frac{2\nu(\theta)|\bDD(\bv)|^2}{\theta} + \frac{\kappa(\theta) |\nabx \theta|^2}{\theta^2} + \theta \int_D \frac{|\nabx \varphi|^2}{\varphi}\dd \bq \\  & \qquad\qquad+ \int_D \frac{4}{\theta \varphi}\left| \theta \nabq\varphi + \theta \varphi \nabq\,\cU_\eta + \varphi \nabq \,\cU_e\right|^2\! \dd \bq :=\xi \ge 0,
\end{aligned}
\end{align}
which, by recalling the definition of $\eta$ (cf. \eqref{eq:eta}), integrating over $\Omega$, and using the boundary conditions $\bv \cdot \boldsymbol{n}_{\bx}|_{\partial \Omega}=0$, \eqref{eq:mbct}$_1$ and \eqref{eq:mbcx} then implies that
\begin{equation}\label{eq:pepa}
\begin{split}
&\frac{\dd}{\dd t}\int_\Omega \bigg[-\log \theta + \int_D \big[\cU_\eta(\bq) \varphi + \varphi \log \varphi\big] \dd \bq \bigg] \dd \bx + \int_\Omega \bigg[\frac{2\nu(\theta)|\bDD(\bv)|^2}{\theta} + \frac{\kappa(\theta) |\nabx \theta|^2}{\theta^2}\bigg]\dd \bx\\ 
& \qquad \qquad + \int_\Omega \bigg[\theta \int_D \frac{|\nabx \varphi|^2}{\varphi}\dd \bq + \int_D \frac{4}{\theta \varphi}\left| \theta \nabq\varphi + \theta \varphi \nabq\,\cU_\eta + \varphi \nabq \,\cU_e\right|^2\! \dd \bq\bigg] \dd \bx = 0.
\end{split}
\end{equation}
%
%

In order to obtain the required evolution equation for the absolute temperature $\theta$, we note
that, by taking the material derivative of \eqref{eq:e} followed by insertion of equations \eqref{eq:nond}$_{3,4}$ into the resulting equality, we obtain
\begin{align*}
\dot\theta &= \dot e -  \int_D \cU_e(\bq) \dot \varphi \dd \bq\\
& = \mathbb{T} : \mathbb{D}(\bv)- \mbox{div}_{\bx}\, \bj_e + \int_D \cU_e(\bq)\big[\mbox{div}_{\bx} \bj_{\varphi,\bx} + \mbox{div}_{\bq}\big((\nabx \bv)\bq \varphi + \bj_{\varphi,\bq}\big)\big] \dd \bq.
\end{align*}
It remains to substitute the defining expressions for $\bTT$, $\bj_e$, $\bj_{\varphi,\bx}$ and $\bj_{\varphi,\bq}$ stated in \eqref{eq:T0}, \eqref{eq:je}, \eqref{eq:7+} and \eqref{eq:7} above into the right-hand side of this equality to infer that the evolution equation for the temperature is the following:
\begin{align}
\begin{aligned}
\dot\theta - \divx(\kappa(\theta)\nabx \theta)&= 2\nu(\theta)|\bDD(\bv)|^2\\
&\qquad + 4\int_D \left(\theta\nabq \varphi \cdot \nabq\, \cU_e  + \theta (\nabq\, \cU_\eta \cdot \nabq\, \cU_e)\varphi + |\nabq\, \cU_e|^2 \varphi  \right)\! \dd \bq \\\label{eq:m4++}
& \qquad +  \bDD(\bv) : \int_D \theta (\nabq\,\cU_\eta \otimes \bq) \varphi \dd \bq.
\end{aligned}
\end{align}

Having completed the derivation of the model, we close this section by collecting the governing equations into a system of nonlinear partial differential equations and supplementing them with boundary and initial conditions.
\bigskip

\noindent
\textbf{Statement of the initial-boundary-value problem.} The initial-boundary-value problem that we will be concerned with in the rest of the paper is the following:
\begin{alignat}{2}
&\divx \bv =0\quad &&\mbox{in $(0,T) \times \Omega$}, \label{eq:m2}\\
&\partial_t \bv + \divx(\bv \otimes \bv) \nonumber\\
&\quad =  \divx \left(- \mathrm{p} \bII+ 2\nu(\theta) \bDD(\bv) - 2 \theta n_{\mathrm{P}} \bII + \int_D (\boldsymbol{F} \otimes \bq) \varphi  \dd \bq\right) + \bff \quad && \mbox{in $(0,T) \times \Omega$}, \label{eq:m1}\\
\partial_t \theta &+ \bv \cdot \nabx \theta - \divx(\kappa(\theta)\nabx \theta)= 2\nu(\theta)|\bDD(\bv)|^2 \nonumber\\
&\qquad + 4\int_D \left(\theta\nabq \varphi \cdot \nabq\, \cU_e  + \theta (\nabq\, \cU_\eta \cdot \nabq\, \cU_e)\varphi + |\nabq\, \cU_e|^2 \varphi  \right)\! \dd \bq \nonumber\\\label{eq:m4}
& \qquad +  \bDD(\bv) : \int_D \theta (\nabq\,\cU_\eta \otimes \bq) \varphi \dd \bq, \quad && \mbox{in $(0,T) \times \Omega$},\\
&\partial_t \varphi + \divx\left(\bv \varphi - \theta \nabx \varphi\right) \nonumber \\\label{eq:m3}
&\qquad + \divq\left((\nabx \bv) \bq \varphi  - 4\boldsymbol{F}\varphi - 4 \theta \nabq \varphi\right) =0\quad &&\mbox{in $(0,T) \times \Omega \times D$},
\end{alignat}
with 
\begin{align}\label{eq:pnd}
n_{\mathrm{P}}:= \int_D \varphi \dd \bq \quad \mbox{and}\quad  \bF:=\nabq\,\cU_e  + \theta \nabq\, \cU_\eta,
\end{align}
This system of partial differential equations will be considered subject to the following boundary and initial conditions, where the hypotheses imposed on the initial data are motivated by the requirement that the equations \eqref{eq:eta},\eqref{eq:E+}, \eqref{eq:en}, \eqref{eq:et}, involving $\bv$, $\theta$ and $\varphi$, be meaningful:

\begin{itemize}
\item \textit{For the velocity field} $\bv$:
\begin{align}\label{eq:mbcv}  
\bv = \boldsymbol{0} \quad \mbox{on $(0,T) \times \partial\Omega$}
\quad \mbox{and}\quad \bv(0,\bx) = \bv_0(\bx)\quad \mbox{for $\bx \in \Omega$},
\end{align}
with $\bv_0 \in L^2(\Omega;\mathbb{R}^d)$, $\divx \bv_0 = 0$ on $\Omega$ in the sense of distributions, and $\bv \cdot \boldsymbol{n}_{\bx}|_{\partial \Omega} = 0$ in $W^{-1/2,2}(\partial \Omega)$;

\item \textit{For the absolute temperature} $\theta$:
\begin{align} \label{eq:mbct}\kappa(\theta) \nabx \theta \cdot \boldsymbol{n}_{\bx}= 0\quad \mbox{on $(0,T) \times \partial \Omega$}\quad\mbox{and}\quad \theta(0,\bx) = \theta_0(\bx) \quad \mbox{for $\bx \in \Omega$};
\end{align}
where $\theta_0 \in L^1(\Omega)$ and there exists a positive constant $\theta_{\mathrm{min}}$ such that $\theta_0(\bx) \geq \theta_{\mathrm{min}}$, $\bx \in \Omega$;
\item \textit{For the probability density function} $\varphi$:
\begin{align}\label{eq:mbcx} 
\theta \nabx \varphi \cdot \boldsymbol{n}_{\bx} = 0 \quad \mbox{on $(0,T) \times \partial \Omega \times D$}, \end{align}
where $\boldsymbol{n}_{\bx}$ is the unit outward normal vector to $\partial\Omega$,
\begin{align}\label{eq:mbcq} ((\nabx \bv) \bq  \varphi  - 4\boldsymbol{F}\varphi - 4 \theta \nabq \varphi)\cdot \boldsymbol{n}_{\bq} =0\quad \mbox{on $(0,T) \times \Omega \times \partial D$},
\end{align}
where $\boldsymbol{n}_{\bq}$ is the unit outward normal vector to $\partial D$, and we impose the initial condition
\begin{align}\label{eq:mic}
\varphi(0,\bx,\bq) = \varphi_0(\bx,\bq)\quad \mbox{for $(\bx,\bq) \in \Omega \times D$},
\end{align}
where 
\[\varphi_0 \in L^1 \log L^1(\Omega \times D; \mathbb{R}_{\geq 0}) \cap L^\infty(\Omega;L^1(D;\mathbb{R}_{ \geq 0}))\quad \mbox{and}\quad \mathcal{U}_e \varphi_0 \in L^1(\Omega \times D; \mathbb{R}_{\geq 0}).\]
\end{itemize}

\medskip
\noindent
\textbf{Assumptions on the material functions.}
We shall assume that 
\begin{align}\label{eq:5}
\mbox{$U_e  \in C^2([0,b/2); \mathbb{R}_{\geq 0}) \cap L^1((0,b/2))$ and $U_\eta \in C^1([0,b/2]; \mathbb{R}_{\geq 0})$},  \hspace{6mm}\nonumber\\
\mbox{$U_e$ is  monotonically increasing and such that} \hspace{22mm}  \nonumber 
\\
\lim_{s \rightarrow (b/2)_{-}} U_e(s) = +\infty,\;
\lim_{s \rightarrow (b/2)_{-}} U'_e(s) = +\infty\; \mbox{and}  \; 0 \leq U''_e(s) \leq c_e [U_e'(s)]^2\; \mbox{for all $s \in [0,b/2)$},\\
\mbox{where $c_e>0$ is a positive constant, independent of $s \in [0,b/2)$}.  \hspace{5mm} \nonumber
\end{align}
%

We note in passing that in the case of the classical FENE (finitely extensible nonlinear elastic) potential proposed by Warner \cite{Warner}, defined by 
\[ \cU_e (\bq):= \frac{Hb}{2} \log \left(1 - \frac{|\bq|^2}{b}\right)^{-1}, \quad \mbox{$|\bq| < \sqrt{b}$},\]
corresponding to 
\[U_e(s):=\frac{Hb}{2}\log(1-2s/b)^{-1}\]
for $s \in [0,b/2)$, where $H>0$ is the spring-constant, hypothesis \eqref{eq:5} is satisfied. The same is true of variants of the FENE potential based on alternative approximations of the inverse $L^{-1}$ of the Langevin function $x \mapsto L(x):=\coth(x)-(1/x)$ (cf. \cite{Ammar} and \cite{Jed}). Another example of a function that satisfies \eqref{eq:5} is the FENE-like potential
\[ \cU_e(\bq):= \frac{Hb}{2} \left[\left(1 - \frac{|\bq|^2}{b}\right)^{-r} -1\right], \quad \mbox{$|\bq| < \sqrt{b}$},\]
with $r \in (0,1)$, corresponding to 
\[U_e(s):=\frac{Hb}{2}[(1-2s/b)^{-r}-1]
\]
for $s \in [0,b/2)$. 

At the very end of the next section (cf. Step 15 there) we shall demand that $U_e$ satisfies a stronger assumption than \eqref{eq:5}; specifically, it will be required there that
\begin{align}\label{eq:5+}
\mbox{$U_e  \in C^2([0,b/2); \mathbb{R}_{\geq 0}) \cap L^1((0,b/2))$ and $U_\eta \in C^1([0,b/2]; \mathbb{R}_{\geq 0})$},  \hspace{6mm}\nonumber\\
\mbox{$U_e$ is  monotonically increasing and such that} \hspace{23mm}  \nonumber
\\
\hspace{-2mm}\lim_{s \rightarrow (b/2)_{-}} U_e(s) = +\infty,\;
\lim_{s \rightarrow (b/2)_{-}} U'_e(s) = +\infty\; \mbox{and}  \; 0 \leq U''_e(s) \leq c_e(s) [U_e'(s)]^2\; \mbox{for all $s \in [0,b/2)$},\\ 
\mbox{where $c_e \in C([0,b/2);\mathbb{R}_{> 0})$, with $\lim_{s \to b/2_{-}} c_e(s)=0$}.\nonumber\hspace{15mm}
\end{align}

Clearly, if \eqref{eq:5+} holds, then so does \eqref{eq:5}. In the case of the classical FENE potential $c_e = 2/Hb$ in \eqref{eq:5}, while \eqref{eq:5+} does not hold. On the other hand, the FENE-like potential stated above satisfies \eqref{eq:5+} with $c_e(s):=(2/Hb)(1-2s/b)^r$ and, a fortiori, \eqref{eq:5} for all $r \in (0,1)$. As will become clear in the next section (cf. Step 15 there), the imposition of the stronger requirement \eqref{eq:5+} becomes necessary because the evolution equation for the absolute temperature $\theta$in the model generally does not provide a uniform upper bound on $\theta$ over $(0,T) \times \Omega$ for an initial datum in $L^1(\Omega; \mathbb{R}_{>0})$. However, until Step 15 of the next section, we shall only require the weaker assumption \eqref{eq:5} to hold. 

The classical FENE potential and its variants are physically motivated: their purpose is to emulate the finite extensibility of polymer molecules; they do so by ensuring that the magnitude of the spring force $\bq \in D \mapsto \boldsymbol{F}(\bq)$ tends to $+\infty$ as $\bq$ approaches the boundary $\partial D$ of the \textit{bounded} configuration space domain $D$. In the context of the model considered here, this singular behaviour of $\boldsymbol{F}$ is encoded in our assumptions on the elastic spring potential $U_e$ (cf. \eqref{eq:5} and \eqref{eq:5+}). Unsurprisingly, the singular behaviour of $U_e$, and thereby also of $\cU_e(\bq):=U_e(\frac{1}{2}|\bq|^2)$, complicates the analysis of the model and has a significant impact on the behaviour of the probability density function $\varphi$ as $\bq$ approaches $\partial D$. This becomes perhaps most evident by observing that in order to ensure that the integral
\[ \int_D |\nabq\, \cU_e|^2 \varphi \dd \bq\]
appearing in the temperature equation
is finite under our assumptions on $U_e$, the function $\varphi$ and its normal derivative have to vanish as $\bq$ approaches $\partial D$. We shall make this assertion rigorous in Lemma \ref{lem:1} below.\footnote{To give a simple illustration of this point, consider the FENE potential $U_e(s):= (Hb/2) \log (1- (2s/b))^{-1}$ and suppose that $\phi \in C^1([0,b/2];\mathbb{R}_{\geq 0}))$. 
Let $g(s):=(\phi(s) - \phi(b/2))/(s-b/2)$ for $s \in [0,b/2)$  with $g(b/2) := \phi'(b/2)$, and note that $g \in C([0,b/2])\cap C^1([0,b/2))$. Suppose further that $g(s) \leq 0$ for all $s \in [0,b/2]$. Then, to ensure that
\[\infty > \int_0^{b/2}|U_e'(s)|^2 \phi(s) \dd s  = 
H^2 \int_0^{b/2}
\frac{\phi(s)}{\big(1-\frac{2s}{b}\big)^2}\dd s =
\frac{H^2 b^2}{4} \left[\phi\bigg(\frac{b}{2}\bigg) \lim_{\sigma \to (b/2)_{-}} \int_0^\sigma \frac{\dd s}{\big(s-\frac{b}{2}\big)^2} + \lim_{\sigma \to (b/2)_{-}} \int_0^\sigma\frac{g(s)}{s- \frac{b}{2}} \dd s\right],  \]
it is necessary that $\phi(b/2)=0$ and $g(b/2)=0$, i.e., both $\phi$ and $\phi'$ must vanish at $s=b/2$.}

Our assumptions on $U_\eta$ imply that there exists a positive constant $C$ such that
\begin{alignat}{2}\label{eq:6}
 |U_\eta'(s)| &\leq C &&\quad  \mbox{for all $s \in [0,\frac{b}{2}]$}. 
\end{alignat}
For example, the classical Hookean spring-potential $\cU_\eta(\bq):= \frac{1}{2}|\bq|^2$ with $\mbox{$|\bq| \leq  \sqrt{b}$}$,
corresponding to $U_\eta(s) = s$ for $s \in [0,b/2]$ satisfies the hypotheses imposed on $U_\eta$.

The coefficient of heat conductivity $\kappa=\kappa(\theta)$ depends on the absolute temperature $\theta$. We shall assume in what follows that $\kappa \in C([0,\infty); \mathbb{R}_{>0})$ and 
\begin{align}\label{eq:9}
 C (1 + \theta^\beta) \leq \kappa(\theta) \leq \hat{C} (1 + \theta^\beta)\quad \mbox{for all $\theta \in [0,\infty)$, where $\beta>\frac{5}{6}$,}
\end{align}
and $C$ and $\hat{C}$ are positive constants, independent of $\theta$. The motivation for the assumed lower bound on the exponent $\beta$ appearing in \eqref{eq:9} will be made clear in the next section in the derivation of the a~priori bound on solutions to this coupled system of nonlinear partial differential equations.  

Concerning the temperature-dependent shear viscosity coefficient $\nu(\theta)$  that appears in the balance of linear momentum equation \eqref{eq:m1} we shall assume that 
\[\nu \in BC([0,\infty);\mathbb{R}_{>0}) \quad \mbox{and} 
\quad
\nu(\theta) \geq C \quad \mbox{for all $\theta \in [0,\infty)$},\] 
where $C$ is a positive constant independent of $\theta$. The density of body forces $\boldsymbol{f}$ appearing in \eqref{eq:m1} will be assumed to belong to $L^2(0,T;L^2(\Omega; \mathbb{R}^d))$. Up to the end of Step 14 in Section \ref{sec:2} the spring-potentials $U_e$ and $U_\eta$ (cf.~\eqref{eq:4a}) are assumed to satisfy \eqref{eq:5}; thereafter \eqref{eq:5+} is assumed. As we have yet to prove that the assumed positivity of $\theta_0$ implies that $\theta(t,\bx)>0$ for all $(t,\bx) \in [0,T] \times \overline{\Omega}$, we shall assume that the material functions $\kappa$ and $\nu$ are even functions of $\theta$; i.e., that $\kappa(\theta)=\kappa(|\theta|)$ and $\nu(\theta) = \nu(|\theta|)$ for all $\theta \in \mathbb{R}$. 
%

The rest of the paper is structured as follows. 

In the next section, we shall derive various formal a~priori bounds. By the word \textit{formal} we mean that we shall assume for the moment that the problem has a smooth solution. The first group of these estimates follows the guideline provided by the thermodynamic approach presented above. After checking that $\varphi$ is nonnegative and $\theta$ is strictly positive, starting from the governing equations \eqref{eq:m2}--\eqref{eq:m3} and knowing what the correct forms of the internal energy $e$ (see \eqref{eq:e}) and the entropy $\eta$ (see \eqref{eq:eta}) should be, we rederive the evolution equations for $e$, $E$ and $\eta$, giving us two fundamental estimates, see \eqref{eq:16} and \eqref{eq:17} below. The second group of estimates consists of improvements that stem from the assumption \eqref{eq:9} on the heat conductivity coefficient, and better properties of $\varphi$ in the interior of $(0,T)\times \Omega \times D$, 
or up to the boundary provided that $\theta$ is bounded above.

We shall then, in Section~\ref{sec:3}, state our definition of weak solution. Interestingly, the concept of solution that we consider here is different from any of the concepts used for the incompressible/compressible Navier--Stokes--Fourier system. There, one has at least three natural options \cite{fm06, BFM, FeiNov17}: (i) to use the equation for $E$, (ii) to use the equation for $e$, and (iii) to use the equation for $\eta$. While these formulations, when supplemented by conservation of the global total energy, are equivalent at the level of classical (regular enough) solutions, they are,---even if we replace the equality signs in the equations for $e$ and $\eta$ with inequalities,---different at the level of variational (weak) solutions. It transpires that, for the problem considered here, we cannot use any of these three formulations. This is due to the impossibility of ensuring integrability of the terms $\theta \nabx \int_D \cU_e(\bq) \varphi \dd \bq$ and $\int_D \big[\cU_\eta(\bq) \varphi + \theta \nabx \varphi \log \varphi\big] \dd \bq$ appearing in the flux terms of the equations for $e$ (and hence also for $E$) and for $\eta$, respectively. As a result of these difficulties and in view of the available estimates, we are thus forced to work with the equation for the temperature \eqref{eq:m4} in a renormalized sense, and have to replace the `=' sign by the `$\geq$' sign therein. 

Finally, also in Section \ref{sec:3}, we provide a rigorous analysis of the model, leading to the key contribution of the paper stated in Theorem \ref{thm:1}: the weak compactness of sequences of solutions to the model \eqref{eq:m1}--\eqref{eq:9}, which are shown to converge to a global-in-time large-data weak solution that satisfies the energy inequality \eqref{eq:17} below.

\section{Formal a priori estimates}\label{sec:2}

We begin by performing some formal preparatory calculations which are required for the derivation of the desired a priori estimates. We shall assume throughout this section that $\bv$, $\varphi$ and $\theta$ are continuous functions of their arguments,  defined on $[0,T] \times \overline \Omega$, $[0,T] \times \overline {\Omega} \times \overline{D}$ and $[0,T] \times \overline \Omega$, respectively, 
and that $\bv$, $\varphi$ and $\theta$ are sufficiently many times continuously differentiable so as to satisfy the partial differential equations that feature in the model, as well as the associated boundary and initial conditions, in the usual pointwise sense. We shall assume that the initial datum $\varphi_0$ for the Fokker--Planck equation is nonnegative  a.e.~on $\overline \Omega \times \overline D$ and that the initial datum $\theta_0$ for the temperature equation is strictly positive a.e.~on $\overline\Omega$; we shall then prove that the corresponding solution $\varphi$ to the Fokker--Planck equation is  nonnegative a.e.~on $[0,T] \times \overline{\Omega} \times \overline D$, and that the solution $\theta$ to the temperature equation is strictly positive a.e.~on $[0,T] \times \overline \Omega$.

More precisely, in this section we assume that $\bv \in C^{1,2}([0,T] \times \overline\Omega;\mathbb{R}^d)$, with $\bv|_{[0,T] \times \partial\Omega}=\boldsymbol{0}$, $\divx \bv = 0$ on $[0,T] \times \overline \Omega$, $\theta \in C^{1,2}([0,T] \times \overline \Omega)$ and $\varphi \in C^{1,2,2}([0,T] \times \overline \Omega \times \overline D)$. 
Motivated by the form of the fourth term on the right-hand side of the evolution equation \eqref{eq:entr} for the specific entropy $\eta$ and the final term of the first integral over $D$ appearing on the right-hand side of the evolution equation \eqref{eq:m4} for the absolute temperature, we shall further suppose that
\[ \int_D \frac{|\nabq \varphi(t,\bx,\bq)|^2}{\varphi(t,\bx,\bq)} \dd \bq < \infty\qquad \mbox{and}\qquad \int_D [\cU_e(\bq)|^2 \varphi(t,\bx,\bq) < \infty \]
for all $(t,\bx) \in [0,T] \times \overline{\Omega}$. 
 
In connection with the imposition of the boundary condition \eqref{eq:mbcq}  an important remark is in order: because the spring force $\boldsymbol{F}= \nabq\,\cU_e + \theta \nabq\,\cU_\eta$ appearing in \eqref{eq:mbcq} involves the gradient of $\cU_e$, and the first and second partial derivatives of $\cU_e$, as well as $\cU_e$ itself, blow up on $\partial D$ because of \eqref{eq:5}, for the boundary condition \eqref{eq:mbcq} to be meaningful, it is necessary that $\varphi|_{(0,T) \times \Omega \times \partial D}=0$. Since we have already imposed a boundary condition on $(0,T) \times \Omega \times \partial D$ in the form of \eqref{eq:mbcq}, we cannot of course impose a second boundary condition on $(0,T) \times \Omega \times \partial D$; the reason why $\varphi$ should vanish on $(0,T) \times \Omega \times \partial D$ therefore requires clarification. In this respect, we note that because of \eqref{eq:5} we have\footnote{From \eqref{eq:5} we have $U'_e(s) \leq U'_e(0) + c_e \int_0^s [U'_e(\sigma)]^2 \dd \sigma$ for all $s \in (0,b/2)$ and then, again by \eqref{eq:5}, it follows that $\int_0^{b/2} [U'_e(\sigma)]^2 \dd \sigma = + \infty$. Consequently, 
 \begin{align*}
 &\int_D |\nabla_{\bq}\, \mathcal{U}_e|^2 \dd \bq = \int_D |U_e'(|\bq|^2/2)|^2 |\bq|^2 \dd \bq \geq \int_{b/2 \leq |\bq|^2 <b} |U_e'(|\bq|^2/2)|^2 |\bq|^2 \dd \bq \geq \frac{b}{2}
 \int_{b/2
 \leq |\bq|^2 <b} |U_e'(|\bq|^2/2)|^2 \dd \bq \\& \geq  C(b) \int_{b/4}^{b/2} [U'_e(\sigma)]^2 \dd \sigma = C(b)\!\int_0^{b/2} [U'_e(\sigma)]^2 \dd \sigma - C(b)\! \int_0^{b/4} [U'_e(\sigma)]^2 \dd \sigma \geq 
 C(b) \! \int_0^{b/2} [U'_e(\sigma)]^2 \dd \sigma - C_1(b) =  +\infty.
 \end{align*}}
\[ \int_D |\nabq\, \cU_e|^2 \dd \bq = +\infty,\]
and therefore the integral
\begin{align}\label{eq:int} 
\int_D |\nabq\, \cU_e|^2 \varphi \dd \bq
\end{align}
appearing on the right-hand side of the temperature equation \eqref{eq:m4} will only be finite if the (nonnegative) function $\varphi$ vanishes on $\partial D$; i.e., 
\begin{align}\label{eq:bc0} 
\varphi = 0 \quad \mbox{on $(0,T) \times \Omega \times \partial D$}.
\end{align}
This homogeneous Dirichlet boundary condition on $\varphi$ is therefore imposed \textit{implicitly} (rather than explicitly as a genuine boundary condition) through the requirement that the integral \eqref{eq:int} involving the singular potential $\cU_e$ be finite. In fact, a stronger assertion than \eqref{eq:bc0} holds, which we shall formulate in the next lemma, and which will be used to justify the absence of various boundary integrals over $\partial D$ that would otherwise arise during the course of partial integration with respect to $\bq$. As solutions to the system of partial differential equations under consideration, which we are formally manipulating at this stage, will be replaced by solutions that are required to satisfy a certain energy (in)equality to be derived later in this section, the presence of the singular potential $\cU_e$ in that energy (in)equality will imply that both $\varphi$ and $\varphi (\nabq \,\cU_e \cdot \boldsymbol{n}_{\bq})$ vanish on $(0,T) \times \Omega \times \partial D$.

\begin{lemma}\label{lem:1}
Let $Q:=(0,T) \times \Omega$ and suppose that $g \in L^1(Q;\mathbb{R}_{> 0})$ and $\varphi \in L^1(Q;W^{1,1}(D;\mathbb{R}_{\geq 0}))$. Suppose further that $U \in C^2([0,\frac{b}{2});\mathbb{R}_{\geq 0})$, $U'$ is nonnegative, $U(s) \to +\infty$ and $U'(s) \to +\infty$ as $s \to (b/2)_{-}$, and that there exists a positive constant $c$ such that 
\begin{align}\label{eq:U21}
|U''(s)| \leq c[U'(s)]^2\quad \forall\, s \in [0,b/2).
\end{align}
Suppose, finally, that the following bound holds: 
\begin{equation}\label{AP1}
\int_{Q\times D} g(t,\bx) \left(\frac{|\nabla_{\bq} \varphi(t,\bx,\bq)|^2}{\varphi(t,\bx,\bq)} + \varphi(t,\bx,\bq)(1+ [U'(|\bq|^2/2)]^2)\right) \!\dd \bq \dd \bx \dd t\le C.
\end{equation}
Then, 
\begin{align}
g \varphi &\in L^1(Q; W^{1,1}_0(D)) \quad \mbox{and}\label{TR1}\\
g \varphi U' &\in L^1(Q; W^{1,1}_0(D)).\label{TR2}
\end{align}
In particular, $\varphi(t,\bx,\cdot)$ and $\varphi(t,\bx,\cdot) U'(|\cdot|^2/2)$ have zero trace in $L^1(\partial D)$ for a.e.~$(t,\bx) \in Q$.
\end{lemma}

\noindent
\textit{Proof} We begin by showing that $g \varphi$ and $g\varphi U'$ both belong to $L^1(Q; W^{1,1}(D))$. Thanks to the assumed nonnegativity of $g$, $\varphi$ and $U'$, we have from \eqref{AP1} that $g\varphi \in L^1(Q \times D) = L^1(Q;L^1(D))$, and 
\begin{align*}
0 &\leq \int_{Q\times D}g(t,\bx) \varphi(t,\bx,\bq)  U'(|\bq|^2/2) \dd \bq \dd \bx \dd t\\
&\le \frac{1}{2}\int_{Q\times D}g(t,\bx) \varphi(t,\bx,\bq) \left(1+ [U'(|\bq|^2/2)]^2\right)\!\dd \bq \dd \bx \dd t \le C,
\end{align*}
whereby $g\varphi U' \in L^1(Q \times D) = L^1(Q; L^1(D))$. Further, we also have that 
\begin{align*}
&\int_{Q \times D}|\nabq (g(t,\bx)\varphi(t,\bx,\bq))| + |\nabq(g(t,\bx)\varphi(t,\bx,\bq) U'(|\bq|^2/2))| \dd \bq \dd \bx \dd t\\
&\quad \leq C\int_Q g(t,\bx) \int_D |\nabq \varphi(t,\bx,\bq)|\big(1 + U'(|\bq|^2/2)\big) + \varphi(t,\bx,\bq) |U''(|\bq|^2/2)| \dd \bq \dd \bx \dd t\\
&\quad \leq C \int_Q g(t,\bx) \int_D \left( \frac{|\nabq \varphi(t,\bx,\bq)|^2}{\varphi(t,\bx,\bq)} + \varphi(t,\bx,\bq)\big(1+ [U'(|\bq|^2/2)]^2\big)  \right) \dd \bq \dd \bx \dd t < \infty,
\end{align*}
where the constant $C$ depends on the diameter $\sqrt{b}$ of $D$ and the constant $c$ from the assumed bound \eqref{eq:U21} in the statement of the lemma. Here in the transition from the second line to the third line we made use of the elementary calculation
\begin{align*}
|\nabq \varphi(t,\bx,\bq)|(1+ U'(|\bq|^2/2) &= \frac{|\nabq \varphi(t,\bx,\bq)|}{\sqrt{\varphi(t,\bx,\bq)}} \sqrt{\varphi(t,\bx,\bq)} (1+ U'(|\bq|^2/2))\\
&\leq \frac{1}{2}  \frac{|\nabq \varphi(t,\bx,\bq)|^2}{\varphi(t,\bx,\bq)} + \varphi(t,\bx,\bq) \big(1 +[U'(|\bq|^2/2)]^2\big).
\end{align*}

Thus, $g\varphi$ and $g\varphi U'$ both belong to $L^1(Q; W^{1,1}(D;\mathbb{R}_{\geq 0}))$. Hence, by the trace theorem, we infer that the traces of $g\varphi$ and $g\varphi U'$ on $\partial D$ are well-defined and belong to $L^1(Q; L^1(\partial D;\mathbb{R}_{\geq 0}))$. 

To complete the proofs of \eqref{TR1} and \eqref{TR2} it remains to show that the traces of $g\varphi$ and $g\varphi U'$ are equal to $0$ a.e.~on $Q \times \partial D$. We shall first show this for $g\varphi U'$, and then proceed analogously in the simpler situation of $g \varphi$.

Let us therefore show that the trace of $g\varphi U'$ is equal to $0$ a.e.~on $Q \times \partial D$. To this end, for $(t,\bx) \in Q$ we introduce
$$
u(t,\bx,r):=\int_{\partial B_r} \varphi (t,\bx, \bq) U'(r^2/2)\dd S(\bq), \quad 0<r<\sqrt{b}, 
$$
where $B_r:=B(\mathbf{0},r)$. Thanks to the fact that $g \varphi U' \in L^1(Q; W^{1,1}(D;\mathbb{R}_{\geq 0}))$, it follows from the trace theorem that $0 \leq \int_Q g(t,\bx) \int_{\partial B_r} \varphi(t,\bx,\bq) U'(r^2/2)\dd S(\bq) \dd \bx \dd t = \int_Q g(t,\bx) u(t,\bx,r) \dd \bx \dd t <\infty$ for all $r \in (0,\sqrt{b})$. Hence, the integrand of the last integral must be finite for a.e.~$(t,\bx) \in Q$ and for all $r \in (0,\sqrt{b})$. As, by hypothesis, $g(t,\bx)>0$ for a.e.~$(t,\bx) \in Q$, it then follows that $u(t,\bx,r)$ must be finite for a.e.~$(t,\bx) \in Q$ and for all $r \in (0,\sqrt{b})$. This means that $u(t,\bx,r)$ is well-defined for a.e.~$(t,\bx)\in Q$ and for all $r\in (0,\sqrt{b})$. Then, from the definition of $u$ and \eqref{AP1}, it follows that
\begin{equation}\label{TR5}
\int_0^{\sqrt{b}} U'(r^2/2) u(t,\bx,r) \dd r = \int_{D} \varphi(t,\bx,\bq) [U'(|\bq|^2/2)]^2\dd \bq < \infty
\end{equation}
for almost all $(t,\bx)\in Q$. In addition, since $\varphi U' \in W^{1,1}(D;\mathbb{R}_{\geq 0})$ for almost all $(t,\bx)\in Q$, we have $u(t,\bx,\cdot)\in \mathcal{C}([0,\sqrt{b}];\mathbb{R}_{\geq 0})$ for almost all $(t,\bx)\in Q$. Clearly, $u$ is a nonnegative function of its arguments. Assume for a moment that, for some $(t,\bx) \in Q$, $\lim_{r\to(\sqrt{b})_{-}} u(t,\bx,r)>0$. Then, thanks to the continuity of $u$ with respect to its third argument, there exist $\varepsilon>0$ and $r_0 \in (0,\sqrt{b})$ such that $u(t,\bx,r) \ge \varepsilon$ for all $r \in [r_0,\sqrt{b}]$. Consequently,
\begin{align*}
+\infty> \int_0^{\sqrt{b}} U'(r^2/2) u(t,\bx,r) \dd r &\ge \varepsilon
 \int_{r_0}^{\sqrt{b}} U'(r^2/2)\dd r =
\frac{\varepsilon}{\sqrt{2}}  \int_{r_0^2/2}^{b/2} U'(s) \frac{\dd s}{\sqrt{s}}\\
& \ge\frac{\varepsilon}{\sqrt{b}} \left(\lim_{s\to (b/2)_{-}}U(s) -U(r_0^2/2) \right)=+\infty.
\end{align*}
Thus, we have arrived at a contradiction, which then implies that, contrary to our supposition that $\lim_{r \to (\sqrt{b})_{-}} u(t,\bx,r)>0$ for some $(t,\bx) \in Q$, in fact,  
\[
\lim_{r \to (\sqrt{b})_{-}} u(t,\bx,r)=0 \quad \mbox{for a.e.~$(t,\bx) \in Q$}.
\]
This implies that, for a.e.~$(t,\bx) \in Q$, 
\[
\lim_{r\to (\sqrt{b})_{-}}
\int_{\partial B_r} \varphi (t,\bx, \bq) U'(r^2/2)\dd S(\bq)=0.
\]
Thanks to the assumed nonnegativity of $\varphi$ and $U'$ this means that $\varphi(t,\bx,\bq) U'(|\bq|^2/2)=0$ for a.e.~$(t,\bx,\bq) \in Q\times \partial D$. That completes the proof of \eqref{TR2}. 
To complete the proof of \eqref{TR1}, we need to show that the trace of $g\varphi$ is equal to $0$ a.e.~on $Q \times \partial D$. Similarly as in the case of $g\varphi U'$ above, we now define
\[ u(t,\bx,r):= \int_{\partial B_r} \varphi(t,\bx,\bq) \dd S(\bq), \quad 0 <r < \sqrt{b},\]
where, again $B_r:=B(\boldsymbol{0},r)$. Instead of \eqref{TR5}, we then have
\begin{equation*}
\int_0^{\sqrt{b}} U'(r^2/2) u(t,\bx,r) \dd r = 
\int_{D} \varphi(t,\bx,\bq) U'(|\bq|^2/2) \leq \frac{1}{2} \int_{D} \varphi(t,\bx,\bq) (1 + [U'(|\bq|^2/2)]^2)\dd \bq < \infty. 
\end{equation*}
By an argument based on contradiction, identical to the one in the case of \eqref{TR2}, we deduce that
\[
\lim_{r\to (\sqrt{b})_{-}} u(t,\bx,r)=0 \quad \mbox{for a.e.~$(t,\bx) \in Q$}.
\]
This implies that, for a.e.~$(t,\bx) \in Q$, 
\[
\lim_{r\to (\sqrt{b})_{-}}
\int_{\partial B_r} \varphi (t,\bx, \bq) \dd S(\bq)=0.
\]
Thanks to the assumed nonnegativity of $\varphi$ this means that $\varphi(t,\bx,\bq)=0$ for a.e.~$(t,\bx,\bq) \in Q\times \partial D$. That completes the proof of \eqref{TR1}.
\hfill $\Box$

\begin{remark}\label{rem:4}
In order to illustrate the relevance of Lemma \ref{lem:1} suppose that $\theta$ and $\varphi$ are such that the inequality \eqref{eq:up}, which we shall prove below, holds, i.e., that, for each positive real number $\theta_{\mathrm{max}}$,
\[
\int_{Q_{\theta_{\mathrm{max}}}}\int_{D} \left(\frac{\theta|\nabq\varphi|^2}{\varphi} + \frac{\varphi |\nabq\, \cU_e|^2}{\theta}\right)\! \dd \bq \leq C(\theta_{\mathrm{max}}),
\]
where $Q_{\theta_{\mathrm{max}}}:=\{(t,\bx) \in Q\,:\, 0 \leq \theta(t,\bx) \leq \theta_{\mathrm{max}}\}$ and $C(\theta_{\mathrm{max}})$ is a positive constant depending on $\theta_{\mathrm{max}}$. It then follows from Lemma \ref{lem:1}  with $Q$ there replaced by $Q_{\theta_{\rm max}}$ and $g$ chosen to be identically equal to  $1$, that the terms 
$\varphi (\nabx \bv)\bq \cdot \boldsymbol{n}_{\bq}$ and $\varphi (\nabq\, \cU_e)\cdot \boldsymbol{n}_{\bq}$, involved in the zero normal flux boundary condition \eqref{eq:mbcq}, i.e., in the boundary condition
\[ 0 =  \left((\nabx \bv)\bq \varphi + \bj_{\varphi,\bq}\right)\cdot \boldsymbol{n}_{\bq} =  \big(\varphi(\nabx \bv)\bq  - 4\varphi(\nabq\, \cU_e + \theta \nabq\,\cU_\eta) - 4 \theta \nabq \varphi\big)\cdot \boldsymbol{n}_{\bq}, \]
that was imposed for $\varphi$ on $(0,T) \times \Omega \times \partial D$, both vanish on  $Q_{\theta_{\mathrm{max}}} \times \partial D$, and therefore by passing to the limit $\theta_{\textrm{max}} \to +\infty$, also on $Q \times \partial D=(0,T) \times \Omega \times \partial D$. In fact, the term $\varphi(\theta \nabq\,\cU_\eta)$ also vanishes on $(0,T) \times \Omega \times \partial D$, because by hypothesis $\cU_\eta \in C^1(\overline{D})$ and the function $\varphi$ vanishes on $(0,T) \times \Omega \times \partial D$. Thus, the boundary condition \eqref{eq:mbcq} reduces to $\theta \nabla_{\bq} \varphi \cdot \boldsymbol{n}_{\bq}=0$ on $(0,T) \times \Omega \times \partial D$, which is of a similar form as the boundary condition \eqref{eq:mbcx} imposed for $\varphi$ on $(0,T) \times \partial \Omega \times D$. For related considerations in the isothermal setting, we refer Zhang \&
Zhang \cite{MR2221211}, Liu \& Liu \cite{MR2407125} and Masmoudi \cite{MR2456183}; see in particular Theorem 1.1 in \cite{MR2407125} and Remark 3.6 in \cite{MR2456183}.
\end{remark}

Motivated by Lemma \ref{lem:1}, in the formal calculations that follow we shall therefore assume that $\varphi$ vanishes on $(0,T) \times \Omega \times \partial D$, and that $\varphi(\cdot,\cdot,\bq)$ decays to $0$ sufficiently fast as $\bq \to \partial D$ to ensure that integrals over $\partial D$ are well defined; moreover, if in the course of our calculations boundary integrals over $\partial D$ happen to disappear following partial integration over $D$, then it will be understood that this is either because of the boundary condition \eqref{eq:mbcq} or because $\varphi(\cdot,\cdot,\bq)$ decays to $0$ sufficiently fast as $\bq \to \partial D$ to annihilate the boundary integral over $\partial D$ in question. In particular, we shall suppose that
\begin{align}\label{eq:1a}
\varphi \in L^1(Q; W^{1,1}_0(D; \mathbb{R})) \quad \mbox{and}
\quad \varphi\, U_e'(|\bq|^2/2) \in L^1(Q; W^{1,1}_0(D; \mathbb{R})).
\end{align}
Because, by hypothesis, $\varphi \in C^{1,2,2}([0,T] \times \overline\Omega \times \overline D))$ and $\bq \mapsto U_e'(|\bq|^2/2) \in C^1(D)$, it follows that $\varphi\, U_e'(|\bq|^2/2) \in C([0,T] \times \overline\Omega \times D))$ and $\lim_{\bq \to \partial D} \varphi\, U_e'(|\bq|^2/2) = 0$. Therefore, $\varphi\, U_e'(|\bq|^2/2)$ is uniformly continuous on $[0,T] \times \overline\Omega \times \overline D$, and thereby it is also bounded on $[0,T] \times \overline\Omega \times \overline D$; that is,
\begin{align*}
\max_{(t,\bx,\bq) \in [0,T] \times \overline\Omega \times \overline D} |\varphi(t,\bx,\bq)| U_e'(|\bq|^2/2) =:C_0 < \infty.
\end{align*}

We are now ready to embark on the derivation of the required (formal) a priori bounds. For the sake of clarity of the exposition, we split the process into sixteen separate steps. 

As $\varphi$ is a probability density function, it is required to be nonnegative. However, the nonnegativity of $\varphi$ is a property that should be deduced from the nonnegativity of the initial datum $\varphi_0$ for the Fokker--Planck equation (using a parabolic minimum principle, for example,) instead of assuming, \textit{a priori}, that it holds. Unfortunately, the parabolicity of the Fokker--Planck equation isn't automatically guaranteed, as it hinges on the positivity of the absolute temperature $\theta$. The absolute temperature is, for physical reasons, positive; however, this property should again not be assumed, but should be inferred from the positivity of the initial temperature $\theta_0$. However, to deduce the positivity of $\theta$ from the positivity of $\theta_0$, one needs nonnegativity of $\varphi$.

To break this vicious circle, for $z \in \mathbb{R}$ let $[z]_{+}:= \max\{z,0\}~(\geq 0)$,  $[z]_{-}:=\min\{z,0\}~(\leq 0)$, and replace for the moment the initial-boundary-value problem \eqref{eq:m3}, \eqref{eq:mbcx}, \eqref{eq:mbcq}, \eqref{eq:mic} for the Fokker--Planck equation with the following: 
\begin{align}\label{eq:m3+}
\partial_t \varphi &+ \divx\left(\bv \varphi - [\theta]_{+} \nabx \varphi\right) \nonumber \\
&+ \divq\left((\nabx \bv) \bq \varphi  - 4\boldsymbol{F}[\varphi]_{+} - 4 [\theta]_{+} \nabq \varphi\right) =0\quad \mbox{in $(0,T) \times \Omega \times D$},&&
\end{align}
with $\bF:=\nabq\,\cU_e  + \theta \nabq\, \cU_\eta$,
subject to the following boundary and initial conditions:
\begin{align}\label{eq:mbcx+} 
[\theta]_{+} \nabx \varphi \cdot \boldsymbol{n}_{\bx} = 0 \quad \mbox{on $(0,T) \times \partial \Omega \times D$}, \end{align}
where $\boldsymbol{n}_{\bx}$ is the unit outward normal vector to $\partial\Omega$, and
\begin{align}\label{eq:mbcq+} ((\nabx \bv) \bq  \varphi  - 4\boldsymbol{F}[\varphi]_{+} - 4 [\theta]_{+} \nabq \varphi)\cdot \boldsymbol{n}_{\bq} =0\quad \mbox{on $(0,T) \times \Omega \times \partial D$},
\end{align}
where $\boldsymbol{n}_{\bq}:=\bq/|\bq|$ is the unit outward normal vector to $\partial D$, and we impose the initial condition
\begin{align}\label{eq:mic+}
\varphi(0,\bx,\bq) = \varphi_0(\bx,\bq)\quad \mbox{for $(\bx,\bq) \in \Omega \times D$},
\end{align}
where $\varphi_0\geq 0$ a.e.~on $\overline\Omega \times D$. 
From the perspective of mathematical modelling, this initial-boundary-value problem is equivalent to \eqref{eq:m3}, \eqref{eq:mbcx}, \eqref{eq:mbcq}, \eqref{eq:mic} because, as noted above, the probability density function $\varphi$ is by definition nonnegative and the absolute temperature $\theta$ is, for physical reasons, positive.

In our initial step, Step 0, of the analysis that follows we shall prove that the assumed nonnegativity of $\varphi_0$ a.e.~on $\overline \Omega \times \overline D$ implies nonnegativity of $\varphi$, viewed as a solution of the problem \eqref{eq:m3+}--\eqref{eq:mic+}, a.e.~on $[0,T] \times \overline\Omega \times \overline{D}$. Using the nonnegativity of $\varphi$ thus proved, we shall then prove in Step 1 that the positivity of the initial temperature $\theta_0$ a.e.~on $\overline \Omega$ guarantees positivity of $\theta$ a.e.~on $[0,T]\times \overline\Omega$. Having proved the nonnegativity of $\varphi$ and the positivity of $\theta$, the initial-boundary-value problem \eqref{eq:m3+}--\eqref{eq:mic+}  will then collapse to its original form \eqref{eq:m3}, \eqref{eq:mbcx}, \eqref{eq:mbcq}, \eqref{eq:mic}, and in the remaining steps of the analysis we shall therefore continue working with that original form. 

\bigskip
\noindent
\textbf{Step 0.} In this step we shall prove that
the assumed nonnegativity of $\varphi_0$ implies the nonnegativity of $\varphi$, viewed as a solution of the problem \eqref{eq:m3+}--\eqref{eq:mic+}, a.e.~on $[0,T] \times \overline\Omega \times \overline{D}$. To this end, we multiply the partial differential equation \eqref{eq:m3+} by $[\varphi]_{-}$, integrate the resulting equality over $\Omega \times D$, and perform partial integration in all terms that involve $\divx$ and $\divq$ using the homogeneous boundary conditions \eqref{eq:mbcx+} and \eqref{eq:mbcq+}. As $\divx\, \bv = 0$ and $\divq ((\nabx \bv) \bq) =0$, it follows by decomposing $\varphi$ as $[\varphi]_{+} + [\varphi]_{-}$ that, for all $t \in (0,T)$,
\[ \frac{1}{2}
\frac{\dd}{\dd t} \int_{\Omega \times D} |[\varphi]_{-}|^2 \dd \bq \dd \bx + 
\int_{\Omega \times D} [\theta]_{+}|\nabx [\varphi]_{-}|^2 \dd \bq \dd \bx + 4 
\int_{\Omega \times D} [\theta]_{+}|\nabq [\varphi]_{-}|^2 \dd \bq \dd \bx = 0,
\]
thanks to the fact that the interiors of the supports $\mathrm{supp}( [\varphi]_{+})$ and $\mathrm{supp}([\varphi]_{-})$
of the functions $[\varphi]_{+}$ and $[\varphi]_{-}$ are disjoint sets. Hence, 
\[ \int_{\Omega \times D} |[\varphi(t,\bx, \bq)]_{-}|^2 \dd \bq \dd \bx \leq \int_{\Omega \times D}[\varphi(0,\bx, \bq)]_{-}|^2 \dd \bq \dd \bx = \int_{\Omega \times D}[\varphi_0(\bx, \bq)]_{-}|^2 \dd \bq \dd \bx = 0\]
for all $t \in [0,T]$. Consequently $\varphi \geq 0$ a.e. on $[0,T] \times \overline{\Omega} \times \overline{D}$.

\bigskip

\noindent\textbf{Step 1.} Having proved the nonnegativity of $\varphi$, we next establish a \textit{minimum principle} for the absolute temperature, $\theta$. To this end, we choose $\them \in (0,1)$ such that $\them \leq \mbox{ess.inf}_{x \in \Omega} \,\theta_0(x)$ and consider the function $H$ defined by 
\begin{align} \label{eq:H}
H(t,\bx):= \theta(t,\bx) - \them\mathrm{e}^{-\alpha t}\quad \mbox{for $t \in [0,T]$ and $\bx \in \Omega$},
\end{align}
with $\alpha>0$ to be fixed below.  For future reference, we also note that, thanks to our assumptions concerning the functions $\cU_e$ and $\cU_\eta$ (cf.~equations \eqref{eq:5} and \eqref{eq:6}), we have the following:
\begin{align*} 
|\nabq\, \cU_\eta| |\bq| = U'_\eta\left(\frac{|\bq|^2} {2}\right)|\bq|^2 \leq C\quad \mbox{for all $\bq \in D$}
\end{align*}
and
\begin{align*} 
0 \leq \nabq \,\cU_\eta  \cdot \nabq\, \cU_e &= |\bq|^2 U_\eta'\left(\frac{|\bq|^2}{2}\right)U_e'\left(\frac{|\bq|^2}{2}\right)\leq C|\bq|^2  \left(1 + \left[U_e'\left(\frac{|\bq|^2}{2}\right)\right]^2\right)\\
& \leq C(1 + |\nabq \,\cU_e|^2)\quad \mbox{for all $\bq \in D$}.
\end{align*}
Furthermore, because 
\begin{align*}
\Delta_{\bq}\,\cU_e = d U_e'\left(\frac{|\bq|^2}{2}\right) + |\bq|^2 U_e''\left(\frac{|\bq|^2}{2}\right), 
\end{align*}
and $U_e \in C^2([0,b/4])$,  it follows from \eqref{eq:5} that 
\begin{align*}
0 \leq \Delta_{\bq}\,\cU_e \leq C\quad \mbox{for $|\bq|^2\leq \frac{b}{2}$},
\end{align*}
where $C$ is a positive constant, independent of $\bq$, and
\begin{align*}
0 &\leq \Delta_{\bq}\,\cU_e  \leq \frac{d}{2}\left(1+\frac{2}{b}|\bq|^2\left|U_e'\left(\frac{|\bq|^2}{2}\right)\right|^2\right) + |\bq|^2 U_e''\left(\frac{|\bq|^2}{2}\right)\\
 &= \frac{d}{2}\left(1+\frac{2}{b}|\nabq\,\cU_e|^2\right) + |\bq|^2 \left[U_e'\left(\frac{|\bq|^2}{2}\right)\right]^2\left\{\left[U_e'\left(\frac{|\bq|^2}{2}\right)\right]^{-2}   U_e''\left(\frac{|\bq|^2}{2}\right)\right\}\\
& = \frac{d}{2}\left(1+\frac{2}{b}|\nabq\,\cU_e|^2\right) + |\nabq\,\cU_e|^2 \left\{\left[U_e'\left(\frac{|\bq|^2}{2}\right)\right]^{-2}   U_e''\left(\frac{|\bq|^2}{2}\right)\right\}
\quad \mbox{for $\frac{b}{2} \leq |\bq|^2 < b$}.
\end{align*}
Thanks to \eqref{eq:5}, the (nonnegative) expression appearing in the curly brackets in the last line is,  for all $\bq \in D$ such that $\frac{b}{2}\leq |\bq|^2 < b$,  bounded by a positive constant independent of $\bq$. In summary, then, writing $C$ for a generic positive constant independent of $\bq$, we have shown that
\begin{alignat}{2}\label{eq:Ubd}
\begin{aligned}
|\nabq\, \cU_\eta| |\bq| &\leq C \quad &&\mbox{for all $\bq \in D$},\\
0 \leq \nabq \,\cU_\eta \cdot \nabq\, \cU_e &\leq C(1 + |\nabq\, \cU_e|^2)  \quad &&\mbox{for all $\bq \in D$},\\
0 \leq \Delta_{\bq}\,\cU_e &\leq C(1 + |\nabq\, \cU_e|^2)\quad &&\mbox{for all $\bq \in D$}.
\end{aligned}
\end{alignat}
With these preparations, we return to the temperature equation \eqref{eq:m4} subject to the boundary and initial conditions \eqref{eq:mbct}, multiply equation \eqref{eq:m4} by $(H(t,\bx))_{-}:= \min\{H(t,\bx),0\} ~ ( \leq 0)$, where $H$ is the function defined by \eqref{eq:H}, integrate the resulting equality over $\bx \in \Omega$ and perform partial integration over $D$ in the first term of the first integral on the right-hand side of \eqref{eq:m4}. Hence, 
\begin{align}\label{eq:tr}
\begin{aligned}
&\int_\Omega \partial_t \theta\, H_{-} \dd \bx + \int_{\{ \bx \in \Omega\,:\, H(t,\bx) \leq 0\}} \kappa(\theta)|\nabx \theta|^2 \dd \bx + \int_\Omega 2\nu(\theta) |H_{-}| |\bDD(\bv)|^2 \dd \bx\\
&\qquad + 4\int_\Omega |H_{-}| \int_D |\nabq\, \cU_e|^2 \varphi \dd \bq
+ 4\int_\Omega \theta |H_{-}| \int_{\partial D}\varphi \nabq\,\cU_e \cdot {\boldsymbol{n}}_{\bq} \dd S(\bq)
\\
&\quad = 4 \int_\Omega
\theta \,|H_{-}| \int_D (\Delta_{\bq} \,\cU_{e})\varphi \dd \bq \dd \bx
+  4\int_\Omega \theta\, H_{-} \int_D (\nabq\, \cU_\eta \cdot \nabq\, \cU_e)\varphi 
\dd \bq \dd \bx\\
&\qquad + \int_\Omega  \theta \, H_{-} \,\bDD(\bv) : \int_D (\nabq\,\cU_\eta \otimes \bq) \varphi \dd \bq =: \mathrm{T}_1 + \mathrm{T}_2 + \mathrm{T}_3,
\end{aligned}
\end{align}
where we have used the trivial equality $H_{-} = - |H_{-}|$. We proceed by bounding the three terms on the right-hand side of this equality. Clearly,
$\theta |H_{-}| = (H + \theta_{\mathrm{min}}\mathrm{e}^{-\alpha t})|H_{-}| \leq \theta_{\mathrm{min}}
\mathrm{e}^{-\alpha t}
|H_{-}| \leq |H_{-}|$ for all $(t,\bx) \in [0,T]\times \Omega$.
Therefore, by \eqref{eq:Ubd}$_3$, 
\begin{align*}
\mathrm{T}_1:=&4 \int_\Omega
\theta \,|H_{-}| \int_D (\Delta_{\bq} \,\cU_{e})\varphi \dd \bq \dd \bx \leq 4 \int_{\{x \in \Omega\,:\, \theta \geq 0\}}
\theta \,|H_{-}| \int_D (\Delta_{\bq} \,\cU_{e})\varphi \dd \bq \dd \bx \\
&\quad\leq 4C  \int_{\{x \in \Omega\,:\, \theta \geq 0\}}
\theta \,|H_{-}| \int_D (1 + |\nabq\,\cU_{e}|^2)\varphi \dd \bq \dd \bx
\\
&\quad \leq  4 C \theta_{\rm min} \int_{\{x \in \Omega\,:\, \theta \geq 0\}}
 \,|H_{-}| \, \mathrm{e}^{-\alpha t} \int_D \varphi \dd \bq \dd \bx + 4C \theta_{\rm min}  \int_{\{x \in \Omega\,:\, \theta \geq 0\}}
\,|H_{-}| \int_D |\nabq\,\cU_{e}|^2\varphi \dd \bq \dd \bx\\
&\quad \leq  4 C_1 \theta_{\rm min} \int_{\Omega}
 \,|H_{-}| \, \mathrm{e}^{-\alpha t} \dd \bx + 4C \theta_{\rm min}  \int_{\Omega}
\,|H_{-}| \int_D |\nabq\,\cU_{e}|^2\varphi \dd \bq \dd \bx, 
\end{align*}
where $C$ here is the constant from inequality \eqref{eq:Ubd}$_3$ and
\[ C_1 := C\, \max_{(t,\bx) \in [0,T] \times \overline{\Omega}} \int_D \varphi \dd \bq.\]

Next, we consider the second term on the right-hand side of the inequality \eqref{eq:tr}. Clearly, 
\begin{align*}
\mathrm{T}_2:=4\int_\Omega \theta\, H_{-} \int_D (\nabq\, \cU_\eta \cdot \nabq\, \cU_e)\varphi 
\dd \bq \dd \bx \leq \frac{1}{\theta_{\rm min}} \int_{\Omega} |H_{-}|\, \theta^2 \int_D |\nabq\,\cU_\eta|^2 \varphi \dd \bq \dd \bx\\
+ 4 \theta_{\rm min} \int_{\Omega } |H_{-}|\, \int_D |\nabq\,\cU_e|^2 \varphi \dd \bq \dd \bx\\
\leq \frac{C_2}{\theta_{\rm min}} \int_{\Omega} |H_{-}|\, \theta^2  \dd \bx + 4 \theta_{\rm min} \int_{\Omega} |H_{-}| \int_D |\nabq\,\cU_e|^2 \varphi \dd \bq \dd \bx, 
\end{align*}
where
\[ C_2:= \max_{(t,\bx) \in [0,T] \times \overline{\Omega}} \int_D |\nabq \, \cU_\eta|^2 \varphi \dd \bq.\]

Finally, we consider the third term on the right-hand side of \eqref{eq:tr}:
\begin{align*}
\mathrm{T}_3:=\int_\Omega  \theta \, H_{-} \,\bDD(\bv) : \int_D (\nabq\,\cU_\eta \otimes \bq) \varphi \dd \bq \leq   \frac{C_3}{\theta_{\rm min}} \int_\Omega |H_{-}|\, \theta^2 \dd \bx + \theta_{\rm \min} \int_\Omega 2\nu(\theta) |H_{-}| |\bDD(\bv)|^2 \dd \bx,   
\end{align*}
where 
\[ C_3:= \frac{1}{8} |D| \max_{(t,\bx) \in [0,T] \times \overline{\Omega}}
\frac{1}{\nu(\theta)} \int_D |(\nabq\,\cU_\eta \otimes \bq))\varphi|^2 \dd \bq.\]

We then substitute these bounds on $\mathrm{T}_1, \mathrm{T}_2$ and
$\mathrm{T}_3$
in \eqref{eq:tr} and choose $\theta_{\rm min} \in (0,1)$ so small that $(C+1) \theta_{\rm min} \leq 1$. This enables us to absorb the second term of the bound on $\mathrm{T}_1$  and the second term of the bound on $\mathrm{T}_2$ into the fourth term on the left-hand side of \eqref{eq:tr}; it also enables us to absorb the second term in the bound on $\mathrm{T}_3$ into the third term on the left-hand side of \eqref{eq:tr}. We note further that
the final term on the left-hand side vanishes thanks to \eqref{eq:1a}.

Hence, recalling the definition \eqref{eq:H} of the function $H$, we have from \eqref{eq:tr}
the inequality
\[ \int_\Omega \partial_t \theta \left(\theta -\them \mathrm{e}^{-\alpha t}\right)_{-} \!\dd \bx \leq  4C_1 \theta_{\rm min} \int_\Omega |H_{-}|\, \mathrm{e}^{-\alpha t} \dd \bx + \frac{C_2 + C_3}{\theta_{\rm min}} \int_\Omega |H_{-}|\, \theta^2 \dd \bx.\]
Equivalently, we have the following inequality:
\begin{align*}
&\frac{1}{2}\frac{\dd}{\dd t} \int_\Omega \left[ \left(\theta -\them \mathrm{e}^{-\alpha t}\right)_{-}\right]^2 \dd \bx - \alpha \them \int_\Omega \left(\theta -\them \mathrm{e}^{-\alpha t}\right)_{-} \mathrm{e}^{-\alpha t} \dd \bx 
\\
&\qquad \leq  4C_1 \theta_{\rm min} \int_\Omega |H_{-}|\, \mathrm{e}^{-\alpha t} \dd \bx + \frac{C_2 + C_3}{\theta_{\rm min}} \int_\Omega |H_{-}|\, \theta^2 \dd \bx,
\end{align*}
whereby, because
\[ -\left(\theta -\them \mathrm{e}^{-\alpha t}\right)_{-} = - H_{-} = |H_{-}|, \]
also
\begin{align*} \frac{1}{2} \frac{\dd}{\dd t} \int_\Omega(H_{-})^2 \dd \bx  \leq 
4C_1 \theta_{\rm min} \int_\Omega |H_{-}|\, \mathrm{e}^{-\alpha t} \dd \bx + \frac{C_2 + C_3}{\theta_{\rm min}} \int_\Omega |H_{-}|\, \theta^2 \dd \bx
- \alpha \them \int_\Omega |H_{-}|\,\mathrm{e}^{-\alpha t} \dd \bx\\
=  
(4C_1 - \alpha)\, \theta_{\rm min} \int_\Omega |H_{-}|\, \mathrm{e}^{-\alpha t} \dd \bx + \frac{C_2 + C_3}{\theta_{\rm min}} \int_\Omega |H_{-}|\, \theta^2 \dd \bx.
\end{align*}
We shall now focus our attention on the last term on the right-hand side of this inequality. 

If $\theta(t,\bx) \geq 0$
at a point $(t,\bx) \in [0,T] \times \overline{\Omega}$, then
\[0 \leq \theta(t,\bx) = \theta(t,\bx)  - \theta_{\rm min} \mathrm{e}^{-\alpha t} + \theta_{\rm min} \mathrm{e}^{-\alpha t} = H(t,\bx)  + \theta_{\rm min} \mathrm{e}^{-\alpha t}.\]
Consequently, because $|H_{-}|\,H  \leq 0$,
\[ |H_{-}(t,\bx)|\, \theta(t,\bx) \leq \theta_{\rm min} \mathrm{e}^{-\alpha t}\, |H_{-}(t,\bx)|,\]
whereby also
\[ |H_{-}(t,\bx)|\, \theta^2(t,\bx) \leq \theta_{\rm min} \mathrm{e}^{-\alpha t}\, |H_{-}(t,\bx)|\, \theta(t,\bx) \leq \theta_{\rm min}^2 \mathrm{e}^{-2\alpha t}\, |H_{-}(t,\bx)| \leq  \theta_{\rm min}^2 \mathrm{e}^{-\alpha t}\, |H_{-}(t,\bx)|,\]
with a strict inequality when $H(t,\bx)<0$ and equality when $H(t,\bx) \geq 0$.

If, on the other hand $\theta(t,\bx) < 0$, then letting
\[ C_\ast:= \frac{1}{\theta_{\rm min}}\mathrm{e}^{\alpha T} \max_{(t,\bx) \in [0,T] \times \overline \Omega} |\theta(t,\bx)|^2,\]
we have
\[ \theta^2(t,\bx) \leq \max_{(t,\bx) \in [0,T] \times \overline \Omega} |\theta(t,\bx)|^2 \leq C_\ast \theta_{\rm min} \mathrm{e}^{-\alpha t},\]
whereby also
\[ \theta^2(t,\bx) + C_\ast \theta(t,\bx) < C_\ast \theta_{\rm min} \mathrm{e}^{-\alpha t}.\]
Thus, 
\[ \theta^2(t,\bx) < C_\ast (\theta_{\rm min} \mathrm{e}^{-\alpha t} - \theta(t,\bx)) = C_{\ast}(-H(t,\bx)) = C_\ast\,|H_{-}(t,\bx)|,\]
and therefore also
\[ |H_{-}(t,\bx)|\,\theta^2(t,\bx) \leq C_\ast |H_{-}(t,\bx)|^2,\]
with a strict inequality when $H(t,\bx)<0$ and equality when $H(t,\bx) \geq 0$.
Hence, for each $t \in (0,T)$,
\begin{align*} 
\frac{1}{2} \frac{\dd}{\dd t} \int_\Omega|H_{-}(t,\bx)|^2 \dd \bx  &\leq 
(4C_1 - \alpha)\, \theta_{\rm min} \int_\Omega |H_{-}(t,\bx)|\, \mathrm{e}^{-\alpha t} \dd \bx + \frac{C_2 + C_3}{\theta_{\rm min}} \int_\Omega |H_{-}(t,\bx)|\, \theta^2(t,\bx) \dd \bx\\
&=
(4C_1 - \alpha)\, \theta_{\rm min} \int_\Omega |H_{-}(t,\bx)|\, \mathrm{e}^{-\alpha t} \dd \bx\\
&\qquad + \frac{C_2 + C_3}{\theta_{\rm min}} \int_{\{\bx \in \Omega\,:\, \theta(t,\bx) \geq 0\}} |H_{-}(t,\bx)|\, \theta^2(t,\bx) \dd \bx \\
&\qquad + \frac{C_2 + C_3}{\theta_{\rm min}} \int_{\{\bx \in \Omega\,:\, \theta(t,\bx) < 0\}} |H_{-}(t,\bx)|\, \theta^2(t,\bx) \dd \bx\\
&\leq 
(4C_1 - \alpha)\, \theta_{\rm min} \int_\Omega |H_{-}(t,\bx)|\, \mathrm{e}^{-\alpha t} \dd \bx \\
&\qquad + (C_2 + C_3)\,\theta_{\rm min}\, \int_{\{\bx \in \Omega\,:\, \theta(t,\bx) \geq 0\}} |H_{-}(t,\bx)|\, \mathrm{e}^{-\alpha t} \dd \bx \\
&\qquad + \frac{C_2 + C_3}{\theta_{\rm min}} C_\ast \int_{\{\bx \in \Omega\,:\, \theta(t,\bx) < 0\}} |H_{-}(t,\bx)|^2 \dd \bx\\
& \leq (4C_1 + C_2 + C_3 - \alpha)\, \theta_{\rm min} \int_\Omega |H_{-}(t,\bx)|\, \mathrm{e}^{-\alpha t} \dd \bx\\
&\qquad + \frac{C_2 + C_3}{\theta_{\rm min}} C_\ast \int_{\Omega} |H_{-}(t,\bx)|^2 \dd \bx.
\end{align*}

Now we fix the value of the constant $\alpha$ by setting $\alpha = 2(4C_1 + C_2 + C_3)>0$. Thus, 
\begin{align*}
    \frac{1}{2} \frac{\dd}{\dd t} \int_\Omega |H_{-}|^2 \dd \bx  &\leq - \frac{\alpha}{2}\, \theta_{\rm min} \int_\Omega |H_{-}|\,\mathrm{e}^{-\alpha t} \dd \bx  +  \frac{C_2 + C_3}{\theta_{\rm min}} C_\ast \int_{\Omega} |H_{-}|^2 \dd \bx. \\
    & \leq \frac{C_2 + C_3}{\theta_{\rm min}} C_\ast \int_{\Omega} |H_{-}|^2 \dd \bx.
\end{align*}
As $|H_{-}(0,\bx)|^2 = ((\theta(0,\bx) - \them)_{-})^2 = ((\theta_0(\bx) - \them)_{-})^2= 0$ for all $\bx \in \overline{\Omega}$, it follows by Gronwall's lemma that 
\[ \int_\Omega|H_{-}(t,\bx)|^2 \dd \bx \leq 0\]
for all $t \in [0,T]$. Thus, we deduce that
\[ H_{-}(t,\bx) = 0 \quad \mbox{for all $(t,\bx) \in \overline Q$}.\]
In other words, $H \geq 0$ on $\overline Q$, i.e., 
\begin{align}\label{eq:20} 
\theta(t,\bx) \geq\them \mathrm{e}^{-\alpha t} \geq\them \mathrm{e}^{-\alpha T} \,\,(>0) \quad \mbox{for all $(t,\bx) \in \overline Q$}, 
\end{align}
with $\them \in (0,1)$ and $\alpha=2(4C_1 + C_2 + C_3)$, where $C_1, C_2, C_3 >0$ are as above. Thus, we have established a strictly positive lower bound on $\theta$ over $[0,T] \times \overline\Omega$.

\bigskip\noindent\textbf{Step 2.}
Next, we bound the \textit{polymer number density} $n_{\mathrm{P}}=\int_{D}\varphi \dd \bq$ (cf.~\eqref{eq:np} and \eqref{eq:pnd}). We integrate the Fokker--Planck equation \eqref{eq:m3} over $D$ and use the assumed zero normal flux boundary condition \eqref{eq:mbcq} on $(0,T) \times \Omega \times \partial D$ to deduce that
\begin{alignat}{2}\label{eq:fpi}
\partial_t \left(\int_D \varphi \dd \bq\right)  + \divx\left(\bv \left(\int_D \varphi \dd \bq\right)\right) - \divx\left(\theta \nabx \left(\int_D \varphi \dd \bq\right)\right) &=0\quad \mbox{in $(0,T) \times \Omega \times D$},
\end{alignat}
subject to the following zero normal flux boundary condition for the function $\int_D \varphi \dd \bq$ (which arises from integrating over the set $D$ the zero normal flux boundary condition \eqref{eq:mbcx} for $\varphi$, imposed on $(0,T) \times \partial \Omega \times D$):
\begin{align}\label{eq:fpi+} 
\left(\bv \left(\int_D \varphi \dd \bq\right) - \theta \nabx \left(\int_D \varphi \dd \bq\right) \right) \cdot \boldsymbol{n}_{\bx} = 0 \quad \mbox{on $(0,T) \times \partial \Omega \times D$}
\end{align}
and the initial condition
\[ \left(\int_D \varphi \dd \bq\right)(0,\bx) = \int_D \varphi_0(\bx,\bq) \dd \bq\quad \mbox{for $\bx \in \Omega$},\]
obtained by integrating the initial condition \eqref{eq:mic} imposed on $\varphi$ over the set $D$.
By the parabolic maximum principle applied to the function $(t,\bx) \in [0,T] \times \overline{\Omega} \mapsto \int_D \varphi(t,\bx,\bq) \dd \bq$ (recall from Step 1 that $\theta$ is strictly positive on $[0,T] \times \overline \Omega$, which guarantees parabolicity of the partial differential equation \eqref{eq:fpi} as an evolution equation for the function $(t,\bx) \to \int_D \varphi(t,\bx,\bq) \dd \bq$), we then have, for all $t \in [0,T]$, the inequality
\begin{align}\label{eq:11}
\underset{x \in \Omega}{\mbox{ess.inf}}\int_D \varphi_0(\bx,\bq) \dd \bq \leq \int_D \varphi(t,\bx,\bq) \dd \bq \leq
\underset{x \in \Omega}{\mbox{ess.sup}} \int_D \varphi_0(\bx,\bq) \dd \bq.
\end{align}
Thanks to the nonnegativity of $\varphi$ (cf. Step 0), it follows from \eqref{eq:11} that 
\begin{align}\label{eq:varphi}
    \|\varphi\|_{L^\infty(Q; L^1(D))} \leq \|\varphi_0\|_{L^\infty(\Omega; L^1(D))}
\end{align}
where, as previously, $Q:=(0,T) \times \Omega$.

\bigskip\noindent\textbf{Step 3.}
As a first step towards deriving evolution equations for the internal energy $e$ and the specific entropy $\eta$, defined, respectively, by 
\[e := \theta + \int_{D} \cU_e(\bq) \varphi \dd \bq \quad \mbox{and}\quad\eta := \log \theta - \int_D [\,\cU_\eta(\bq) \varphi  + \varphi \log \varphi ]\dd \bq,\]
(cf.~\eqref{eq:e} and \eqref{eq:eta}) recall from \eqref{eq:m4} that the evolution equation for $\theta$ is
\begin{align}\label{eq:12}
\begin{aligned}
\partial_t \theta + \divx(\bv \theta) &- \divx(\kappa(\theta)\nabx \theta) \\
& = 2\nu(\theta)|\bDD(\bv)|^2 + 4 \int_D \left(\theta\nabq \varphi \cdot \nabq\, \cU_e  + \theta (\nabq\, \cU_\eta \cdot \nabq\, \cU_e)\varphi + |\nabq\, \cU_e|^2 \varphi  \right)\! \dd \bq \\
& \quad +  \bDD(\bv) : \int_D \theta (\nabq\,\cU_\eta \otimes \bq) \varphi \dd \bq.
\end{aligned}
\end{align}
Dividing \eqref{eq:12} by $\theta$ (note that Step 1 guarantees that $\theta$ is bounded below by a positive constant throughout $[0,T] \times \overline{\Omega}$) then gives
\begin{align}\label{eq:13}
\begin{aligned}
 &\partial_t \log \theta + \divx(\bv \log \theta) - \divx(\kappa(\theta) \nabx \log \theta)   = \frac{2\nu(\theta)|\bDD(\bv)|^2}{\theta} + \frac{\kappa(\theta)|\nabx \theta|^2}{\theta^2} \\
 &\quad+ \bDD(\bv) : \int_D (\nabq\,\cU_\eta \otimes \bq) \varphi \dd \bq
  + 4\int_D \left(\nabq \varphi \cdot \nabq\, \cU_e  + \varphi (\nabq\, \cU_\eta \cdot \nabq\, \cU_e) + \frac{\varphi |\nabq\, \cU_e|^2}{\theta} \right)\! \dd \bq. 
\end{aligned}
\end{align}
\bigskip\noindent\textbf{Step 4.}
As a second step towards deriving an evolution equation for the internal energy we multiply the Fokker--Planck equation \eqref{eq:m3} by $\cU_e(\bq)$ and integrate the resulting equality over $D$. Partial integration with respect to $\bq$ in the final term then gives 
\begin{align}\label{eq:11a}
\begin{aligned}
\partial_t\left(\int_D \cU_e \varphi \dd \bq\right) + \divx\left(\bv \int_D \cU_e \varphi \dd \bq\right) &- \divx  \left( \theta \nabx \int_D \cU_e \varphi \dd \bq\right)  \\
& - \int_D \left((\nabx \bv)\bq \varphi + \bj_{\varphi,\bq}\right) \cdot \nabq\, \cU_e(\bq) \dd \bq = 0,
\end{aligned}
\end{align}
where
\[ \bj_{\varphi,\bq} := - 4\boldsymbol{F}\varphi - 4 \theta \nabq \varphi = -4\varphi(\nabq\, \cU_e + \theta\,\nabq\,\cU_\eta) - 4 \theta \nabq \varphi.\]
Because $e := \theta + \int_D \cU_e(\bq) \varphi \dd \bq$, summing \eqref{eq:12} and \eqref{eq:11a} and recalling the definition \eqref{eq:T0} of the Cauchy stress tensor, we find that $e$ satisfies the partial differential equation \eqref{eq:nond}$_3$; that is, 
\[\partial_t e + \divx(\bv e) + \divx \left(  - \kappa(\theta) \nabx \theta - \theta \nabx \int_D \cU_e(\bq) \varphi \dd \bq \right) =\bTT : \bDD(\bv) \quad \mbox{in $(0,T) \times \Omega$}, \]
where the expression in parentheses in the last term on the left-hand side is what was defined in \eqref{eq:je} as the
energy flux $\boldsymbol{j}_e = \boldsymbol{j}_e(t,\bx)$.
\bigskip\noindent\textbf{Step 5.}
As a second step towards deriving an evolution equation for the specific entropy, we return to the Fokker--Planck equation \eqref{eq:m3}, but this time we multiply it by $\cU_\eta$ and again integrate the resulting equality over $D$ and perform partial integration; hence, we deduce that
\begin{align}\label{eq:14}
&\partial_t \left(\int_D \cU_\eta \varphi \dd \bq\right) + \divx \left(\bv \int_D \cU_\eta \varphi \dd \bq\right) - \divx \left(\theta \,\nabx \int_D \cU_\eta \varphi \dd \bq \right)\nonumber\\
&\quad =\int_D \bDD(\bv) : (\bq \otimes \nabq\, \cU_\eta)\varphi \dd \bq - 4\int_D (\bF \varphi + \theta \nabq \varphi) \cdot \nabq\, \cU_\eta \dd \bq \nonumber\\
&\quad = \int_D \left[\bDD(\bv) : (\bq \otimes \nabq\, \cU_\eta)\varphi - 4\varphi (\nabq\, \cU_e \cdot \nabq\, \cU_\eta) -4 \theta \varphi |\nabq\, \cU_\eta|^2 -4 \theta (\nabq \varphi \cdot \nabq\,\cU_\eta)\right]\!\dd \bq, 
\end{align}
by recalling the definition \eqref{eq:F} of the spring force $\bF$.

\bigskip\noindent\textbf{Step 6.}
As a third step towards deriving an evolution equation for the specific entropy, we multiply the Fokker--Planck equation by $(1+ \log (\varphi + \sigma)) \zeta$, where $\sigma>0$ and $\zeta \in C^\infty_0(Q)$, integrate the resulting equality over $Q\times D$, perform partial integration with respect to $t$ and $\bx$ using that $\zeta$ is smooth with compact support, and also with respect to $\bq$ using the zero normal flux boundary conditions imposed on $\varphi$. We then pass to the limit $\sigma \to 0_{+}$ in the resulting equality, term by term, using the monotone/dominated convergence theorem. Having passed to the limit $\sigma \to 0$, we reverse the partial integration with respect to $t$ and $\bx$ using the assumed smoothness of $\varphi$. Since the resulting equality holds for all test functions $\zeta \in C^\infty_0(Q)$, it follows that
\begin{align}\label{eq:15}
&\partial_t \left(\int_D \varphi \log \varphi \dd \bq\right) + \divx \left(\bv \int_D \varphi \log \varphi \dd \bq \right) - \divx\left(\theta \nabx \int_D \varphi \log \varphi \dd \bq \right)\nonumber\\
&\quad=\int_D \left[\nabx \bv : (\bq \otimes \nabq \varphi) - 4\nabq\,\cU_e \cdot \nabq \varphi - 4\theta \nabq\, \cU_\eta \cdot \nabq\varphi  - 4\frac{\theta |\nabq \varphi|^2}{\varphi} - \frac{\theta |\nabx \varphi|^2}{\varphi}\right]\! \dd \bq \nonumber\\
&\quad=\int_D \left[- 4\nabq\,\cU_e \cdot \nabq \varphi - 4\theta \nabq\, \cU_\eta \cdot \nabq\varphi  - 4\frac{\theta |\nabq \varphi|^2}{\varphi} - \frac{\theta |\nabx \varphi|^2}{\varphi}\right]\! \dd \bq, 
\end{align}
where in the transition to the last line we made use of the fact that 
\[ \int_D \nabx \bv : (\bq \otimes \nabq \varphi) \dd \bq = \nabx \bv : \int_D (\bq \otimes \nabq \varphi) \dd \bq = - (\nabx \bv : \bII) \int_D \varphi \dd \bq = - (\divx \bv) \int_D \varphi \dd \bq = 0, \]
by partial integration with respect to $\bq$, using the implicitly imposed boundary condition \eqref{eq:bc0}, and recalling the divergence-free property \eqref{eq:m2} of the velocity field $\bv$.

\bigskip\noindent\textbf{Step 7.}
We shall now combine Steps 3, 5 and 6 to infer the evolution equation for the specific entropy $\eta$. Subtracting \eqref{eq:14} and \eqref{eq:15} from \eqref{eq:13} and recalling the definition \eqref{eq:eta} of the specific entropy $\eta$, it follows that
\begin{align}\label{eq:ent}
&\partial_t \eta + \divx(\bv \eta) + \divx\left(\theta \int_D \cU_\eta \nabx \varphi \dd \bq + \theta \nabx \int_D \varphi \log \varphi \dd \bq \right) - \divx\left(\frac{\kappa(\theta) \nabx \theta}{\theta}\right)\nonumber\\
&\quad = \frac{2\nu(\theta)|\bDD(\bv)|^2}{\theta} + \frac{\kappa(\theta) |\nabx \theta|^2}{\theta^2} + \theta \int_D \frac{|\nabx \varphi|^2}{\varphi}\dd \bq + 4 \theta \int_D \frac{|\nabq \varphi|^2}{\varphi}\dd \bq\nonumber\\
&\qquad + 4\int_D \left[\frac{|\nabq \,\cU_e|^2 \varphi}{\theta} + 2(\nabq \varphi \cdot \nabq\, \cU_e) + 2\varphi (\nabq\, \cU_\eta \cdot \nabq\, \cU_e) + 2\theta (\nabq \varphi \cdot \nabq\, \cU_\eta)  + \theta \varphi |\nabq\, \cU_\eta|^2  \right]\! \dd \bq\nonumber\\
&\quad = \frac{2\nu(\theta)|\bDD(\bv)|^2}{\theta} + \frac{\kappa(\theta) |\nabx \theta|^2}{\theta^2} + \theta \int_D \frac{|\nabx \varphi|^2}{\varphi}\dd \bq + \int_D \frac{4}{\theta \varphi}\left| \theta \nabq\varphi + \theta \varphi \nabq\,\cU_\eta + \varphi \nabq \,\cU_e\right|^2\! \dd \bq.
\end{align}

\smallskip
\noindent\textbf{Step 8.} Next, we shall establish the balance of total energy. We take the scalar product of \eqref{eq:m1} with $\bv$, integrate the resulting equality over $(0,t) \times \Omega$, where $t \in (0,T)$, and perform partial integration with respect to $\bx$ using \eqref{eq:mbcv}. We then combine this with \eqref{eq:12} integrated over $(0,t) \times \Omega$ after performing partial integration using the boundary condition \eqref{eq:mbct}, and with \eqref{eq:11a} integrated over $(0,t)\times \Omega$ after performing partial integration using the boundary condition \eqref{eq:mbcx}, to deduce that 
\begin{align}\label{eq:16+}
\begin{aligned}
&\int_\Omega \frac{1}{2} |\bv(t)|^2 \dd \bx  + \int_\Omega \theta(t) \dd \bx + \int_{\Omega \times D}  \cU_e \varphi(t) \dd \bq \dd \bx \\
&\qquad = \int_\Omega \frac{1}{2} |\bv_0|^2 \dd \bx  + \int_\Omega \theta_0 \dd \bx + \int_{\Omega \times D}  \cU_e \varphi_0 \dd \bq \dd \bx + \int_0^t \int_\Omega \bff(s) \cdot \bv(s)\dd \bx \dd s.
\end{aligned}
\end{align}
Looking back at the definition \eqref{eq:e} of the internal energy $e=e(t,\bx)$ and of the total energy $E$ defined in \eqref{eq:E+}, i.e., 
\[ E(t,\bx) = \frac{1}{2} |\bv(t,\bx)|^2  + e(t,\bx) = \frac{1}{2} |\bv(t,\bx)|^2   + \theta(t,\bx) + \int_{D}  \cU_e(\bq) \varphi(t,\bx,\bq) \dd \bq,\]
one sees that the process leading to \eqref{eq:16+} coincides with the integration of the evolution equation \eqref{eq:E} for the total energy $E$, first over $\Omega$ in conjunction with the use of the divergence theorem and the relevant boundary conditions to deduce \eqref{eq:en}, and then the integration of \eqref{eq:en} over $t$, which yields
\[ \int_\Omega E(t,\bx) \dd \bx = \int_\Omega E(0,\bx) \dd \bx+ \int_0^t \int_{\Omega} \boldsymbol{f}(s,\bx) \cdot \bv(s,\bx) \dd \bx \dd s.\]
Comparing this with \eqref{eq:16+}, we see that the equality \eqref{eq:16+} derived above expresses the balance of total energy, in fact.

\smallskip

Applying Gronwall's inequality to \eqref{eq:16+}, we have
\begin{align} \label{eq:16}
\begin{aligned}
&\underset{t \in (0,T)}{\mbox{ess.sup}}\left(\frac{1}{2}\|\bv(t)\|^2_{L^2(\Omega)} + \|\theta(t)\|_{L^1(\Omega)} + \|\,\cU_e \varphi \|_{L^1(\Omega \times D)}\right)\\
&\quad \leq \mathrm{e}^{T}\left(\frac{1}{2}\|\bv_0\|^2_{L^2(\Omega)} + \|\theta_0\|_{L^1(\Omega)} + \|\,\cU_e \varphi_0 \|_{L^1(\Omega \times D)} + \frac{1}{2}\|\bff\|^2_{L^2(Q)}\right) \leq C.
\end{aligned}
\end{align}
Here and hereafter $C$ signifies a generic positive constant that depends only on the data. 

\bigskip\noindent\textbf{Step 9.}
In this crucial step, we derive the formal energy equality that lies at the heart of the compactness argument presented in Section \ref{sec:3}. We multiply \eqref{eq:ent} by $-1$ and integrate the resulting equality over $\Omega$. Using the boundary condition \eqref{eq:mbcv} to eliminate the boundary integral originating from the second term on the left-hand side of \eqref{eq:ent} upon integration over $\Omega$, the boundary condition \eqref{eq:mbcx} to eliminate the boundary integral that comes from integrating the third term on the left-hand side of \eqref{eq:ent} over $\Omega$, and the boundary condition \eqref{eq:mbct} to eliminate the boundary integral that stems from integrating the fourth term on the left-hand side of \eqref{eq:ent} over $\Omega$, we have
\begin{align}
\begin{aligned}\label{eq:ent1}
&\frac{\dd}{\dd t}\int_\Omega -\eta \dd \bx + \int_\Omega \bigg[\frac{2\nu(\theta)|\bDD(\bv)|^2}{\theta} + \frac{\kappa(\theta) |\nabx \theta|^2}{\theta^2}\bigg]\dd \bx\\ 
&\quad + \int_\Omega \bigg[\theta \int_D \frac{|\nabx \varphi|^2}{\varphi}\dd \bq + \int_D \frac{4}{\theta \varphi}\left| \theta \nabq\varphi + \theta \varphi \nabq\,\cU_\eta + \varphi \nabq \,\cU_e\right|^2\! \dd \bq\bigg] \dd \bx = 0.
\end{aligned}
\end{align}
We note that because the second and third terms on the left-hand side of \eqref{eq:ent1} are nonnegative, the \textit{total entropy}, $t \in [0,T] \mapsto \int_\Omega \eta(t,\bx) \dd \bx$, is an increasing function of $t$, as is physically expected.

Thanks to the definition \eqref{eq:eta} of the specific entropy $\eta$ we have
\[ -\eta = - \log \theta + \int_D [\,\cU_\eta \varphi  + \varphi \log\varphi] \dd \bq.\]
By substituting this into \eqref{eq:ent1} it follows that
\begin{align*}
&\frac{\dd}{\dd t}\int_\Omega \bigg[-\log \theta + \int_D [\,\cU_\eta \varphi + \varphi \log \varphi] \dd \bq \bigg] \dd \bx + \int_\Omega \bigg[\frac{2\nu(\theta)|\bDD(\bv)|^2}{\theta} + \frac{\kappa(\theta) |\nabx \theta|^2}{\theta^2}\bigg]\dd \bx\\ 
&\qquad + \int_\Omega \bigg[\theta \int_D \frac{|\nabx \varphi|^2}{\varphi}\dd \bq + \int_D \frac{4}{\theta \varphi}\left| \theta \nabq\varphi + \theta \varphi \nabq\,\cU_\eta + \varphi \nabq \,\cU_e\right|^2\! \dd \bq\bigg] \dd \bx = 0.
\end{align*}
By integrating equation \eqref{eq:fpi} over the domain $\Omega$ and using the boundary condition \eqref{eq:fpi+}, it follows that $\frac{\dd}{\dd t}\int_{\Omega \times D} (- \varphi + 1) \dd \bq \dd \bx =0$; therefore,
\begin{align*}
&\frac{\dd}{\dd t}\int_\Omega \bigg[-\log \theta + \int_D \big[\,\cU_\eta \varphi + [\varphi ( \log \varphi-1) + 1 ]\big] \dd \bq \bigg] \dd \bx + \int_\Omega \bigg[\frac{2\nu(\theta)|\bDD(\bv)|^2}{\theta} + \frac{\kappa(\theta) |\nabx \theta|^2}{\theta^2}\bigg]\dd \bx\\ 
&\qquad + \int_\Omega \bigg[\theta \int_D \frac{|\nabx \varphi|^2}{\varphi}\dd \bq + \int_D \frac{4}{\theta \varphi}\left| \theta \nabq\varphi + \theta \varphi \nabq\,\cU_\eta + \varphi \nabq \,\cU_e\right|^2\! \dd \bq\bigg] \dd \bx = 0.
\end{align*}
Defining 
\begin{align}\label{eq:HF}
\mathcal{H}(s):=s - 1 - \log s\quad \mbox{for $s>0$} \qquad \mbox{and}\qquad \mathcal{F}(s):= \left\{\begin{array}{ll} s (\log s - 1) + 1 & \mbox{for $s > 0$,}\\ 1 & \mbox{for $s=0$,}
\end{array} \right.
\end{align}
we then have
\begin{align*}
&\frac{\dd}{\dd t}\int_\Omega \bigg[\mathcal{H}(\theta) + \int_D [\,\cU_\eta\varphi + \mathcal{F}(\varphi)]\dd \bq  \bigg] \dd \bx + \! \int_\Omega \bigg[\frac{2\nu(\theta)|\bDD(\bv)|^2}{\theta} + \frac{\kappa(\theta) |\nabx \theta|^2}{\theta^2}\bigg]\dd \bx\\ 
&\qquad + \int_\Omega \bigg[\theta \int_D \frac{|\nabx \varphi|^2}{\varphi}\dd \bq + \int_D \frac{4}{\theta \varphi}\left| \theta \nabq\varphi + \theta \varphi \nabq\,\cU_\eta + \varphi \nabq \,\cU_e\right|^2\! \dd \bq\bigg] \dd \bx \\
&\quad 
= \frac{\dd}{\dd t}\int_\Omega (\theta - 1) \dd \bx = \frac{\dd}{\dd t}\int_\Omega \theta \dd \bx .
\end{align*}
Note that $s \mapsto \mathcal{H}(s)$ and $s\mapsto \mathcal{F}(s)$ are strictly convex nonnegative functions of $s$, which only vanish at $s=1$. By integrating this equality with respect to $t$, it follows that
\begin{align*}
&\int_\Omega \bigg[\mathcal{H}(\theta(t)) + \int_D [\,\cU_\eta \varphi(t) + \mathcal{F}(\varphi(t))]\dd \bq  \bigg] \dd \bx + \int_0^t \! \int_\Omega \bigg[\frac{2\nu(\theta)|\bDD(\bv)|^2}{\theta} + \frac{\kappa(\theta) |\nabx \theta|^2}{\theta^2}\bigg]\dd \bx \dd s\\ 
&\qquad + \int_0^t \int_\Omega \bigg[\theta \int_D \frac{|\nabx \varphi|^2}{\varphi}\dd \bq + \int_D \frac{4}{\theta \varphi}\left| \theta \nabq\varphi + \theta \varphi \nabq\,\cU_\eta + \varphi \nabq \,\cU_e\right|^2\! \dd \bq\bigg] \dd \bx \dd s\\
&
= \int_\Omega \bigg[\mathcal{H}(\theta_0) + \int_D [\,\cU_\eta \varphi_0 + \mathcal{F}(\varphi_0)]\dd \bq  \bigg] \dd \bx+ \int_\Omega  \theta(t) \dd \bx - \int_\Omega  \theta_0 \dd \bx.
\end{align*}
Adding this to the balance of total energy stated in \eqref{eq:16+} then yields the following  \underline{\textit{energy equality}}:
%
\begin{align}
\begin{aligned}\label{eq:energy}
&\int_\Omega \bigg[\frac{1}{2}|\bv(t)|^2 +  \mathcal{H}(\theta(t)) \bigg] \dd \bx  + \int_{\Omega \times D} [\,\cU_e \varphi(t) + \cU_{\eta} \varphi(t) + \mathcal{F}(\varphi(t))]\dd \bq \dd \bx \\
&\qquad+ \int_0^t \! \int_\Omega \bigg[\frac{2\nu(\theta)|\bDD(\bv)|^2}{\theta} + \frac{\kappa(\theta) |\nabx \theta|^2}{\theta^2}\bigg]\dd \bx \dd s\\ 
&\qquad + \int_0^t \int_\Omega \bigg[\theta \int_D \frac{|\nabx \varphi|^2}{\varphi}\dd \bq + \int_D \frac{4}{\theta \varphi}\left| \theta \nabq\varphi + \theta \varphi \nabq\,\cU_\eta + \varphi \nabq \,\cU_e\right|^2\! \dd \bq\bigg] \dd \bx \dd s\\
& =  \int_\Omega \bigg[\frac{1}{2}|\bv_0|^2 +  \mathcal{H}(\theta_0) \bigg] \dd \bx  + \int_{\Omega \times D} [\,\cU_e \varphi_0 + \cU_{\eta} \varphi_0 + \mathcal{F}(\varphi_0)]\dd \bq \dd \bx + \int_0^t \int_\Omega \bff(s) \cdot \bv(s)\dd \bx \dd s.
\end{aligned}
\end{align}
This formal energy equality will play a crucial role in Section \ref{sec:3}. In particular, once a solution is known to satisfy \eqref{eq:energy} (or a weaker version of \eqref{eq:energy}, with the equality sign replaced by the symbol $\leq$), it will automatically follow that one can meet the hypotheses of Lemma \ref{lem:1}, as a consequence of which the functions $\varphi$ and $\varphi \,\mathcal{U}_e$ will be forced to have zero trace on $(0,T)\times \Omega \times \partial D$. In other words, instead of being explicitly imposed, these homogeneous boundary conditions for $\varphi$ and $\varphi \, \mathcal{U}_e$  will be implied by the energy equality (or an energy inequality, where the equality sign in \eqref{eq:energy} is replaced by the symbol $\leq$). 

The energy equality \eqref{eq:energy} further implies that
\begin{align}\label{eq:17}
&\underset{t \in (0,T)}{\mbox{ess.sup}} \left(\frac{1}{2}\|\bv(t)\|^2_{L^2(\Omega)}+ \| \mathcal{H}(\theta(t))\|_{L^1(\Omega)} + \|\,\cU_e \varphi(t)\|_{L^1(\Omega \times D)}  + \|\,\cU_\eta \varphi(t)\|_{L^1(\Omega \times D)} + \|\mathcal{F}(\varphi(t))\|_{L^1(\Omega \times D)} \right)\nonumber\\
&\quad + \int_{0}^T \|\log \theta\|^2_{W^{1,2}(\Omega)} + \|\theta^{\beta/2}\|^2_{W^{1,2}(\Omega)} \dd t + \int_Q \frac{2\nu(\theta)|\bDD(\bv)|^2}{\theta} + \frac{\kappa(\theta) |\nabx \theta|^2}{\theta^2} \dd \bx \dd t \nonumber\\
&\quad + \int_Q \left[ \theta \int_D \frac{|\nabx \varphi|^2}{\varphi}\dd \bq + \int_D \frac{4}{\theta \varphi}\left| \theta \nabq\varphi + \theta \varphi \nabq\,\cU_\eta + \varphi \nabq \,\cU_e\right|^2\! \dd \bq \right]\! \dd \bx \dd t\nonumber\\
& \leq C \bigg(\,\frac{1}{2}\|\bv_0\|^2_{L^2(\Omega)} + \|\mathcal{H}(\theta_0)\|_{L^1(\Omega)}  +  \|\,\cU_e \varphi_0 \|_{L^1(\Omega \times D)} + \|\,\cU_\eta \varphi_0\|_{L^1(\Omega \times D)} \nonumber\\
&\qquad \qquad + \|\mathcal{F}(\varphi_0)\|_{L^1(\Omega \times D)} + \|\bff\|^2_{L^2(Q)}\bigg) \leq C,
\end{align}
where $\beta>\frac{5}{6}$, as was assumed in \eqref{eq:9}. In connection with the final integral on the left-hand side of \eqref{eq:17} we note that
\begin{align*}
\int_{Q \times D} \frac{1}{\theta \varphi}   \left| \theta \nabq\varphi + \theta \varphi \nabq\,\cU_\eta + \varphi \nabq \,\cU_e\right|^2\! \dd \bq \dd \bx \dd t   = \int_{Q \times D} \frac{1}{\theta \varphi} \left|\theta \, \mathrm{e}^{-\frac{\cU_e + \theta\, \cU_\eta}{\theta}} \nabq \left(\varphi \,\mathrm{e}^{\frac{\cU_e + \theta\, \cU_\eta}{\theta}}\right)\right|^2
\dd \bq \dd \bx \dd t.
\end{align*}
Hence, by defining the rescaled probability density function
\begin{align}\label{eq:res} \widehat\varphi:= \varphi \,\mathrm{e}^{\frac{\cU_e + \theta\, \cU_\eta}{\theta}},
\end{align}
we then have from \eqref{eq:17} that
\begin{align}\label{eq:18}
\begin{aligned}
&\int_{Q \times D} \frac{1}{\theta \varphi}   \left| \theta \nabq\varphi + \theta \varphi \nabq\,\cU_\eta + \varphi \nabq \,\cU_e\right|^2\! \dd \bq \dd \bx \dd t  
\\
& \quad = \int_{Q \times D}
\theta \, \mathrm{e}^{-\frac{\cU_e + \theta\, \cU_\eta}{\theta}} \frac{|\nabq \widehat{\varphi}|^2}{\widehat\varphi} \dd \bq \dd \bx \dd t 
\leq C.
\end{aligned}
\end{align}
Furthermore, recalling the definition of $\bj_{\varphi,\bq}$ from \eqref{eq:7}, we have from \eqref{eq:18} that
\begin{align}\label{eq:19}
\int_{Q \times D}\frac{1}{\theta \varphi} \,
 |\,\bj_{\varphi,\bq}|^2 \dd \bq \dd \bx \dd t \leq C.
\end{align}
We also note that by adding $\int_\Omega 1 \dd \bx = |\Omega|$ to both sides of the energy equality \eqref{eq:17}, the norm $\|\mathcal{H}(\theta(t))\|_{L^1(\Omega)}$ on the left-hand side of \eqref{eq:17} is replaced by $ \int_\Omega (\theta(t) - \log \theta(t))\dd \bx$ and $\|\mathcal{H}(\theta_0)\|_{L^1(\Omega)}$ on the right-hand side is replaced by $\int_\Omega (\theta_0 - \log \theta_0)\dd \bx$. Because 
\[|\theta_0 - \log \theta_0| \leq \theta_0 + |\log \theta_0| \leq \theta_0 + \max\{\theta_0, |\log \theta_{\mathrm{min}}|\}\] 
(cf. Step 1 for the choice of $\theta_{\mathrm{min}} \in (0,1)$), and 
\[\theta(t) - \log \theta(t) \geq \frac{1}{4}(\theta(t) + |\log(\theta(t))| \geq 0,\] 
it follows in particular that 
\[ \underset{t \in (0,T)}{\mbox{ess.sup}} \left(\|\theta(t)\|_{L^1(\Omega)} + \|\log \theta(t)\|_{L^1(\Omega)} \right)\leq C.\]

\smallskip\noindent\textbf{Step 10.}
We note at this point in relation to the integral 
\[ \int_D (\bF \otimes \bq) \varphi \dd \bq\]
appearing in the definition of the Cauchy stress tensor $\bTT$ that because the function $\theta$ is independent of the configuration space variable $\bq \in D$, recalling the definition \eqref{eq:res} of the rescaled probability density function $\widehat{\varphi}
$, this integral can be rewritten in an equivalent form as follows:
\begin{align}\label{eq:19a}
\begin{aligned}
\int_D (\bF \otimes \bq)\varphi \dd \bq &= \theta \int_D \mathrm{e}^{-\frac{\cU_e + \theta \,\cU_\eta}{\theta}} \left(\nabq \frac{\cU_e + \theta\, \cU_\eta}{\theta} \otimes \bq\right)\varphi\,  \mathrm{e}^{\frac{\cU_e + \theta \,\cU_\eta}{\theta}} \dd \bq\\
 &=-\theta \int_D \left(\nabq \mathrm{e}^{-\frac{\cU_e + \theta \,\cU_\eta}{\theta}}\right)\otimes \left(\bq\, \varphi\, \mathrm{e}^{\frac{\cU_e + \theta\, \cU_\eta}{\theta}}\right)\! \dd \bq\\
 &= d\, \theta \left(\int_D \mathrm{e}^{-\frac{\cU_e + \theta \,\cU_\eta}{\theta}} \varphi \, \mathrm{e}^{\frac{\cU_e + \theta\, \cU_\eta}{\theta}}      \dd \bq \right)\, \bII + \theta \int_D \mathrm{e}^{-\frac{\cU_e + \theta \,\cU_\eta}{\theta}} \left(\bq \otimes \nabq \left(\varphi\, \mathrm{e}^{\frac{\cU_e + \theta\, \cU_\eta}{\theta}}\right) \right)\!\dd \bq\\
 &=d\,\theta\,n_{\text{p}}\, \bII + \theta \int_D \mathrm{e}^{-\frac{\cU_e + \theta \,\cU_\eta}{\theta}} \left(\bq \otimes \nabq \widehat\varphi\right)\!\dd \bq\\
 &=d\,\theta\,n_{\text{p}}\, \bII +  \theta \int_D \sqrt{\widehat\varphi} \, \mathrm{e}^{-\frac{\cU_e +\theta \,\cU_\eta}{\theta}} \left(\bq \otimes \frac{\nabq \widehat\varphi}{\sqrt{\widehat\varphi}}\right)\!\dd \bq\\
 &=d\,\theta\,n_{\text{p}}\, \bII +  \int_D \sqrt{\theta}\, \mathrm{e}^{-\frac{\cU_e +\theta \,\cU_\eta}{2\theta}}\,\sqrt{\widehat\varphi} \left(\bq \otimes  
 \sqrt{\theta}\, \mathrm{e}^{-\frac{\cU_e +\theta \,\cU_\eta}{2\theta}}\, \frac{\nabq \widehat\varphi}{\sqrt{\widehat\varphi}}\right)\!\dd \bq\\
 &=d\,\theta\,n_{\text{p}}\, \bII +  \int_D \sqrt{\theta}\, \sqrt{\varphi} \left(\bq \otimes  
 \sqrt{\theta}\, \mathrm{e}^{-\frac{\cU_e +\theta \,\cU_\eta}{2\theta}}\, \frac{\nabq \widehat\varphi}{\sqrt{\widehat\varphi}}\right)\!\dd \bq.
 \end{aligned}
 \end{align}
 In the transition from the second displayed line to the third displayed line we performed partial integration with respect to $\bq$ and made use of the implicitly imposed homogeneous Dirichlet boundary condition \eqref{eq:bc0}. Since, by the Cauchy--Schwarz inequality and the bound \eqref{eq:18},
\begin{align*}
\int_Q \int_D \left|\sqrt{\theta}\, \sqrt{\varphi} \left(\bq \otimes  
 \sqrt{\theta}\, \mathrm{e}^{-\frac{\cU_e +\theta \,\cU_\eta}{2\theta}}\, \frac{\nabq \widehat\varphi}{\sqrt{\widehat\varphi}}\right)\right|\!\dd \bq \dd \bx \dd t\leq \left(\int_Q \theta \left(\int_D \varphi \dd \bq\right)\dd \bx \dd t\right)^{\frac{1}{2}}\\
\times \left(\int_{Q\times D}\left| \bq \otimes  
 \sqrt{\theta}\, \mathrm{e}^{-\frac{\cU_e +\theta \,\cU_\eta}{2\theta}}\, \frac{\nabq \widehat\varphi}{\sqrt{\widehat\varphi}} \right|^2\! \dd \bq \dd \bx \dd t\right)^{\frac{1}{2}}\\
 \leq \sqrt{b} \left(\int_Q \theta \, n_{\text{p}} \dd \bx \dd t\right)^{\frac{1}{2}} \left(\int_{Q \times D} 
 \theta\, \mathrm{e}^{-\frac{\cU_e +\theta \,\cU_\eta}{\theta}} \frac{|\nabq \widehat\varphi|^2}{\widehat\varphi}\,\dd \bq \dd \bx \dd t\right)^{\frac{1}{2}}
 \leq C \sqrt{b} \left(\int_Q \theta \, n_{\text{p}} \dd \bx \dd t\right)^{\frac{1}{2}}, 
\end{align*}
where, as before, $n_{\text{p}}(t,\bx) :=\int_D \varphi(t,\bx, \bq) \dd \bq$ and, by \eqref{eq:varphi}, 
\[ \left(\int_Q \theta \, n_{\text{p}} \dd \bx \dd t\right)^{\frac{1}{2}} \leq \|\varphi_0\|_{L^\infty(Q; L^1(D))}^{\frac{1}{2}} \left(\int_Q \theta \dd \bx \dd t\right)^{\frac{1}{2}}, \]
it follows from the bound on the second term on the left-hand side of inequality \eqref{eq:16} that 
\[ \sqrt{\theta}\, \sqrt{\varphi} \left(\bq \otimes  
 \sqrt{\theta}\, \mathrm{e}^{-\frac{\cU_e +\theta \,\cU_\eta}{2\theta}}\, \frac{\nabq \widehat\varphi}{\sqrt{\widehat\varphi}}\right) \in L^1(0,T;L^1(\Omega \times D;\mathbb{R}^{d \times d})).\]
Similarly, 
\[ d\, \theta\,n_{\text{p}}\, \bII  \in L^\infty(0,T;L^1 (\Omega;\mathbb{R}^{d \times d})).\]
Thus, the polymeric contributions to the Cauchy stress tensor $\bTT$, that is, the tensors 
\[ d\,\theta\, n_{\text{p}\,} \bII \quad \mbox{and}\quad \int_D (\bF \otimes \bq)\varphi \dd \bq,\]
belong to $L^1(0,T;L^1 (\Omega;\mathbb{R}^{d \times d}))$, and therefore they are well-defined under the stated assumptions on the data in the sense that they have finite Frobenius norm, a.e.~on $Q$. We note in passing that the right-hand side of \eqref{eq:19a} can be further transformed as follows:
\begin{align}\label{eq:19b}
\int_D (\bF \otimes \bq)\varphi \dd \bq = 
d\,\theta\,n_{\text{p}}\, \bII + \theta \int_D \left( \nabq \log \widehat{\varphi} \otimes \bq \right)\varphi\dd \bq.
\end{align}

\bigskip\noindent\textbf{Step 11.}
To proceed, we require a bound on $\theta$ in the norm of $L^r(Q)= L^r(0,T;L^r(\Omega))$, for a suitable $r>1$. We begin by noting that the bounds \eqref{eq:16} and \eqref{eq:17} yield
\begin{align}\label{eq:21}
\|\theta\|_{L^\infty(0,T; L^1(\Omega))} \leq C \quad \mbox{and}\quad \|\theta^{\beta/2}\|_{L^2(0,T;W^{1,2}(\Omega))} \leq C, 
\end{align}
where $C$ is a positive constant that depends only on the data, and $\beta>\frac{5}{6}$.
As $\beta>\frac{5}{6}> \frac{1}{3}$, it follows that $1< \frac{2}{3}+\beta < 3\beta$, so we can interpolate between $L^1(\Omega)$ and $L^{3\beta}(\Omega)$ to obtain the following bound in $L^{\frac{2}{3} + \beta}(Q) = L^{\frac{2}{3} + \beta}(0,T;L^{\frac{2}{3} + \beta}(\Omega))$:
\begin{align}\label{eq:22}
\begin{aligned}
\int_0^T \|\theta\|_{L^{\frac{2}{3}+\beta}(\Omega)}^{\frac{2}{3}+\beta} \dd t &\leq \int_0^T \|\theta\|_{L^1(\Omega)}^{\frac{2}{3}}\|\theta\|_{L^{3\beta}(\Omega)}^\beta \dd t\\
& \leq \|\theta\|_{L^\infty(0,T; L^1(\Omega))} ^{\frac{2}{3}} \int_0^T \|\theta\|_{L^{3\beta}(\Omega)}^\beta \dd t\\
&=  \|\theta\|_{L^\infty(0,T; L^1(\Omega))}^{\frac{2}{3}} \int_0^T \|\theta^\frac{\beta}{2}\|_{L^{6}(\Omega)}^2 \dd t\\
& \leq C\|\theta\|_{L^\infty(0,T; L^1(\Omega))}^{\frac{2}{3}} \int_0^T \|\theta^\frac{\beta}{2}\|_{W^{1,2}(\Omega)}^2 \dd t \leq C,
\end{aligned}
\end{align}
where in the transition to the last inequality, we have used the bounds of \eqref{eq:21}.

\bigskip\noindent\textbf{Step 12.}
Next, we derive an additional bound on the velocity field $\bv$. To this end, we define
\[ p:= 2 - \frac{6}{5 + 3\beta}~.\]
Again, because $\beta>\frac{5}{6}>\frac{1}{3}$, it follows that $p \in (1,2)$. Hence, by Young's inequality, the bound \eqref{eq:17}, and the assumed strict positivity of the shear viscosity coefficient $\nu(\theta)$, 
\begin{align}\label{eq:23}
\begin{aligned}
\int_Q |\mathbb{D}(\bv)|^p \dd \bx \dd t &= \int_Q \left(\frac{|\mathbb{D}(\bv)|^2}{\theta}\right)^{\frac{p}{2}} \theta^{\frac{p}{2}} \dd \bx \dd t \\
&\leq \frac{p}{2} \int_Q \frac{|\mathbb{D}(\bv)|^2}{\theta} \dd \bx \dd t + \frac{2-p}{2} \int_Q \theta^{\frac{p}{2-p}} \dd \bx \dd t\\
&\leq C + \frac{3}{5+3\beta}\int_Q \theta^{\frac{2-\frac{6}{5+3\beta}}{2-2+\frac{6}{5+3\beta}}}\dd \bx \dd t \\
&= C + \frac{3}{5+3\beta}\int_Q \theta^{\frac{2}{3}+\beta} \dd \bx \dd t \leq C,
\end{aligned}
\end{align}
where the last inequality follows from the bound \eqref{eq:22}.

Our assumption that $\beta>\frac{5}{6}$ guarantees that, with $p$ as defined above, $\bv \in L^{\frac{5p}{3}}(Q;\mathbb{R}^d)$ with $\frac{5p}{3}>2$. Indeed, because $\bv \in L^p(0,T;W^{1,p}(\Omega;\mathbb{R}^d)) \cap L^\infty(0, T; L^2(\Omega;\mathbb{R}^d))$ by parabolic interpolation\footnote{Lemma (Parabolic Interpolation, cf.~Thm.~2.1 in \cite{DiBenedetto}). Suppose that $\Omega$ is a bounded open Lipschitz domain in $\mathbb{R}^d$, $d \geq 2$. Let $p$, $q \in [1, \infty)$, $s \in [p, \infty)$,  and suppose that 
\[ \gamma :=\frac{\frac{1}{p} - \frac{1}{s}}{\frac{1}{d} -\frac{1}{q}+\frac{1}{p}} \in [0,1].\]
Then, there exists a constant $c > 0$ such that
\[ \|v\|_{L^r(0,T;L^s(\Omega))} \leq  c \|\nabla v\|^\gamma 
_{L^q(Q)} \|v\|^{1-\gamma}_{L^\infty(0,T;L^p(\Omega))}\]
for all $v \in L^q(0,T;W^{1,q}(\Omega)) \cap L^\infty(0, T; L^p(\Omega))$,  with 
\[r:=\frac{s (q(p + d) - dp)}{(s - p)d} \in (1,\infty].\]
}
we have that 
\begin{align*}
\int_Q |\bv|^{\frac{5p}{3}} \dd \bx \dd t \leq C,
\end{align*}
where $\frac{5p}{3} = \frac{10}{3} \frac{2 + 3 \beta}{5 + 3\beta} > 2$ thanks to our assumption that $\beta>\frac{5}{6}$.

\bigskip\noindent\textbf{Step 13.}
We proceed by deriving a bound on $\mathbb{T}$. Because by hypothesis $\nu(\cdot)$ is bounded from below and from above by positive constants, \eqref{eq:23} implies that the term $\nu(\theta) \mathbb{D}(\bv)$ is bounded in $L^p(Q;\mathbb{R}^{d \times d})$ with $p$ as defined above. Next, we bound the polymeric part of the stress tensor. Recall from \eqref{eq:19a} that
\begin{align*}
   \int_D (\boldsymbol{F} \otimes \bq ) \varphi \dd \bq
   & = d\, \theta\, n_{\mathrm{P}}\,\bII + \int_D \sqrt{\theta}\, \sqrt{\varphi}\left(\bq \otimes \sqrt{\theta} \,\mathrm{e}^{-\frac{\cU_e + \theta \cU_\eta}{2\theta}}\, \frac{\nabq \widehat\varphi}{\sqrt{\widehat\varphi}}\right)\!\dd \bq\\
   & = d\, \theta\, n_{\mathrm{P}}\,\bII + \int_D
   \underbrace{\left(\frac{\nabq \widehat\varphi}{\sqrt{\widehat\varphi}}  
   \sqrt{ \theta \,\mathrm{e}^{-\frac{\cU_e + \theta \cU_\eta}{\theta}}} \otimes \bq \right)}_{z_1} 
   \underbrace{\sqrt{\widehat\varphi \, \theta\,  \mathrm{e}^{-\frac{\cU_e + \theta \cU_\eta}{\theta}}}}_{z_2}\dd \bq,
\end{align*}
where $\widehat\varphi$ is the rescaled probability density function defined in \eqref{eq:ren}. 
Thanks to \eqref{eq:18}, the term $z_1$ belongs to $L^2(Q \times D;\mathbb{R}^{d\times d})$. Also, 
\begin{align*} 
\int_{Q \times D} z_2^2 \dd \bx \dd \bq \dd t &= \int_{Q \times D} \varphi\, \theta \dd \bx \dd \bq \dd t \\ &= \int_Q \theta \left(\int_D \varphi \dd \bq\right)\! \dd \bx \dd t \\& = \int_Q \theta\, n_{\mathrm{P}} \dd \bx \dd t\\
& \leq \|n_{\mathrm{P}}\|_{L^\infty(Q)}\int_Q \theta \dd \bx \dd t \leq C,
\end{align*}
whereby $z_2 \in L^2(Q \times D)$. Hence, for a suitable $q \in (1,2)$, to be fixed below,
\begin{align*}
   \int_Q \left|\theta \int_D \mathrm{e}^{-\frac{\cU_e + \theta \cU_\eta}{\theta}} \bq \otimes \nabq \widehat\varphi \dd \bq\right|^q  \dd \bx \dd t  &= \int_Q\left| \int_D z_1 z_2 \dd \bq\right|^q \dd \bx \dd t\\ 
   & \leq \int_Q \left(\int_D z_1^2 \dd \bq\right)^{\frac{q}{2}} \left(\int_D z_2^2 \dd \bq\right)^{\frac{q}{2}} \dd \bx \dd t\\
   & = \int_Q \left(\int_D z_1^2 \dd \bq\right)^{\frac{q}{2}} \left(n_{\mathrm{P}}\,\theta \right)^{\frac{q}{2}} \dd \bx \dd t\\
   & \leq \left(\int_{Q \times D} z_1^2 \dd \bq \dd \bx \dd t\right)^{\frac{q}{2}} \left(\int_Q (n_{\mathrm{P}}\, \theta)^{\frac{q}{2-q}} \dd \bx \dd t\right)^{\frac{2-q}{2}}\\
   &\leq \left(\int_{Q \times D} z_1^2 \dd \bq \dd \bx \dd t\right)^{\frac{q}{2}} \|n_{\mathrm{P}}\|_{L^\infty(Q)} ^{\frac{q}{2}} \left(\int_Q \theta^{\frac{q}{2-q}} \dd \bx \dd t\right)^{\frac{2-q}{2}}\\
   & \leq C  \left(\int_Q \theta^{\frac{q}{2-q}} \dd \bx \dd t\right)^{\frac{2-q}{2}}.
\end{align*}
Because, thanks to \eqref{eq:22}, $\theta \in L^{\frac{2}{3} + \beta}(Q)$, we choose $q \in (1,2)$ by setting
\[ \frac{q}{2-q} = \frac{2}{3} + \beta.\]
Expressing $q$ in terms of $\beta$ from this equality, it follows that
\[ q = \frac{4 + 6\beta}{5 + 3\beta} = 2 - \frac{6}{5 + 3 \beta} = p.\]
Consequently,
\begin{align*}
\int_D (\bF \otimes \bq)\varphi \dd \bq \in L^p(Q;\mathbb{R}^{d \times d})\quad \mbox{with}\quad p=2 - \frac{6}{5+3\beta}, \quad \beta>\frac{5}{6}.
\end{align*}
Recall that
\[ \bTT := - \mathrm{p} \bII+ 2\nu(\theta) \bDD(\bv) - 2 \theta n_{\mathrm{P}} \bII + \int_D (\boldsymbol{F} \otimes \bq) \varphi  \dd \bq.\]
Hence, the sum of the second and fourth terms of $\bTT$ belongs to $L^p(Q;\mathbb{R}^{d \times d})$. That is, 
\begin{align*}
\int_Q |\mathbb{S}|^p \dd \bx \dd t = \int_Q |\bTT + \mathrm{p} \bII+ 2 \theta n_{\mathrm{P}} \bII|^p &\dd \bx \dd t \leq C\quad \mbox{for}\quad p=2 - \frac{6}{5 + 3 \beta}\quad\mbox{with} \quad \beta>\frac{5}{6}.
\end{align*}
Furthermore, since for $\beta> \frac{1}{3}$ (which holds in our case, since it was assumed that $\beta>\frac{5}{6}$) one has
\[ \frac{2}{3} + \beta > 2 - \frac{6}{5+3\beta} = p,\]
it follows that 
\[ 2\theta n_{\mathrm{P}} \bII \in L^{\frac{2}{3}+\beta}(Q;\mathbb{R}^{d \times d}) \subset L^p(Q;\mathbb{R}^{d \times d}),\]
and therefore also
\begin{align*}
\int_Q |\bTT + \mathrm{p} \bII|^p &\dd \bx \dd t \leq C\quad \mbox{for}\quad p=2 - \frac{6}{5 + 3 \beta}\quad\mbox{with}\quad \beta>\frac{5}{6}.
\end{align*}

\bigskip\noindent\textbf{Step 14.}
Next, we derive bounds on the probability density function $\varphi$. We begin by noting that, thanks to \eqref{eq:varphi}, \eqref{eq:16}, \eqref{eq:17} and \eqref{eq:19} we have that 
\begin{align}\label{eq:26}
\begin{aligned}
&\|\varphi\|_{L^\infty(Q;L^1(D))} +\|\cU_e \varphi\|_{L^\infty(0,T;L^1(\Omega \times D))} + \|\cU_\eta \varphi\|_{L^\infty(0,T;L^1(\Omega \times D))}\\ 
&\quad + \|\varphi \log\varphi\|_{L^\infty(0,T;L^1(\Omega \times D))} + \left\|\frac{\bj_{\varphi,\bq}}{\sqrt{\theta \varphi}}\right\|_{L^2(Q \times D)} + \left\|\sqrt{\frac{\theta}{\varphi}}\,\nabx \varphi\right\|_{L^2(Q \times D)} \leq C. 
\end{aligned}
\end{align}
Suppose that $\wtD$ is an arbitrary open subset of $D$ such that $\wtD\subset\joinrel\subset D$ (that is, the closure of $\wtD$ is contained in the open ball $D$). Let us define
\[ K_{\wtD}:= \max\{ \sup_{\bq \in \wtD} |\nabq\, \cU_e(\bq)|, \, \sup_{\bq \in \wtD} |\nabq\, \cU_\eta (\bq)|\} \,\, (< \infty) \]
and note the elementary inequality $|a+b+c|^2 \geq \frac{1}{2}|a|^2 - 2(|b|^2 + |c|^2)$ for $a, b, c \in \mathbb{R}^d$. Then
\begin{align} \label{eq:27}
\begin{aligned}
&\int_{Q \times \wtD} \frac{1}{\theta \varphi} \left|\theta \nabq \varphi + \theta \varphi (\nabq \,\cU_\eta) + \varphi(\nabq\, \cU_e)\right|^2 \dd \bq \dd \bx \dd t\\
& \qquad  \geq \frac{1}{2}\int_{Q \times \wtD} \frac{\theta |\nabq \varphi|^2}{\varphi} \dd \bq \dd \bx \dd t - 2 K_{\wtD}^2 \int_{Q \times \wtD} \varphi \left(\theta + \frac{1}{\theta}\right) \dd \bq \dd \bx \dd t.
\end{aligned}
\end{align}
Now, 
\begin{align*} 0 &\leq \int_{Q \times \wtD} \varphi \left(\theta + \frac{1}{\theta}\right) \dd \bq \dd \bx \dd t = \int_{Q} \left(\theta + \frac{1}{\theta}\right) \left(\int_{\wtD} \varphi \dd \bq \right) \dd \bx \dd t\\
&\leq \int_{Q} \left(\theta + \frac{1}{\theta}\right) \left(\int_{D} \varphi \dd \bq \right) \dd \bx \dd t = \int_{Q} \left(\theta + \frac{1}{\theta}\right) n_{\mathrm{P}}  \dd \bx \dd t\\
&\leq \|n_{\mathrm{P}}\|_{L^\infty(Q)} \int_{Q} \left(\theta + \frac{1}{\theta}\right) \dd \bx \dd t. 
\end{align*}
Using \eqref{eq:varphi}, \eqref{eq:20} and \eqref{eq:21} it follows that 
\[0 \leq \int_{Q \times \wtD} \varphi \left(\theta + \frac{1}{\theta}\right) \dd \bq \dd \bx \dd t \leq C.\]
Substituting this bound into the right-hand side of inequality \eqref{eq:27} and recalling that thanks to \eqref{eq:18} the expression on the left-hand side of \eqref{eq:27} is bounded by a constant, it then follows from \eqref{eq:17} that there exists a positive constant $C_{\wtD}$, dependent on $\wtD$, such that
\[ \frac{1}{2}\int_{Q \times \wtD} \frac{\theta |\nabq \varphi|^2}{\varphi} \dd \bq \dd \bx \dd t \leq C_{\wtD}.\]
In addition, it trivially follows from the bound on the last term on the left-hand side of \eqref{eq:26} that 
\[ \int_{Q \times \wtD} \frac{\theta |\nabx \varphi|^2}{\varphi} \dd \bq \dd \bx \dd t \leq C.\]
Thus, by combining the last two bounds, we have the following bound:
\begin{align}\label{eq:28}
\int_{Q \times \wtD} \theta \left(|\nabx \sqrt{\varphi}|^2 + |\nabq \sqrt{\varphi}|^2\right)\! \dd \bq \dd \bx \dd t \leq C_{\wtD}.
\end{align}
Thanks to \eqref{eq:20}, the function $\theta$ (which is independent of $\bq \in \wtD$) is bounded below by a positive constant on $\overline Q$. In addition, \eqref{eq:varphi} implies that
\[ \int_{Q \times \wtD} |\sqrt{\varphi}|^2 \dd \bq \dd \bx \dd t \leq C.\]
Combining this with \eqref{eq:28}, it follows that  
\begin{align}\label{eq:28a}
\|\sqrt{\varphi}\|_{L^2(0,T;W^{1,2}(\Omega \times \wtD))} \leq C_{\wtD}.
\end{align}
By parabolic interpolation between this inequality and the inequality
\begin{align*}
\|\sqrt{\varphi}\|_{L^\infty(0,T;L^2(\Omega \times \wtD))} \leq C_{\wtD}
\end{align*}
implied by \eqref{eq:varphi}, it follows that
\begin{align}\label{eq:29} 
\int_{Q \times \wtD} |\varphi|^{1+\delta} \dd \bx \dd t \leq C_{\wtD}\quad \mbox{for some $\delta>0$}.
\end{align}
A simple calculation yields $1+\delta= 1 + \frac{1}{d} \in (1,2)$ for $d=2,3$, independently of $\wtD$.

\bigskip\noindent\textbf{Step 15.}
In this final step, we derive bounds on $\varphi$ up to the boundary of $D$. To do this, suppose that $G \in L^\infty([0,\infty);\mathbb{R}_{\geq 0})$. It follows from the bound stated in \eqref{eq:18} that
\begin{align}\label{eq:Gt}
\begin{aligned}
C &\geq \int_{Q \times D} \frac{G(\theta)}{\theta \varphi}\left| \theta \nabq\varphi + \theta \varphi \nabq\, \cU_\eta + \varphi \nabq\, \cU_e\right|^2 \dd \bq \dd \bx \dd t\\
&= \int_{Q \times D}\bigg[ \frac{\theta |\nabq\varphi|^2}{\varphi} + \theta \varphi |\nabq\, \cU_\eta|^2 + \frac{\varphi |\nabq\, \cU_e|^2}{\theta} \\
&\qquad  \qquad + 2 \theta \nabq \varphi \cdot \nabq\, \cU_\eta + 2\nabq \varphi \cdot \nabq\,\cU_e  + 2 \varphi \nabq\, \cU_\eta \cdot \nabq\,\cU_e \bigg] G(\theta)  \dd \bq \dd \bx \dd t.
\end{aligned}
\end{align}
We proceed by considering, for $(t,\bx) \in Q$, the integral over $D$ of the expression appearing in the square brackets on the right-hand side of \eqref{eq:Gt}, that is, 
\begin{align*}
\mathcal{A}(t,\bx)&:=\int_{D}\bigg[ \frac{\theta |\nabq\varphi|^2}{\varphi} + \theta \varphi |\nabq\, \cU_\eta|^2 + \frac{\varphi |\nabq\, \cU_e|^2}{\theta} \\
&\qquad  \quad + 2 \theta \nabq \varphi \cdot \nabq\, \cU_\eta + 2\nabq \varphi \cdot \nabq\,\cU_e  + 2 \varphi \nabq\, \cU_\eta \cdot \nabq\,\cU_e \bigg] \dd \bq.   
\end{align*}
We shall use the elementary inequality $2a\cdot b \geq  -\frac{1}{2}|a|^2 - 2|b|^2$, with $a, b \in \mathbb{R}^d$, to partially absorb the fourth term into the first term and the sixth term into the third term. Indeed, because
\[ \int_D  2 \theta \nabq \varphi \cdot \nabq\, \cU_\eta \dd \bq \geq -\frac{1}{2}\int_D \frac{\theta|\nabq\varphi|^2}{\varphi} \dd \bq - 2 \int_D \theta \varphi |\nabq\,\cU_\eta|^2 \dd \bq\]
and
\[ \int_D  2 \varphi \nabq\, \cU_\eta \cdot \nabq\,\cU_e \dd \bq \geq -\frac{1}{2}\int_D \frac{\varphi|\nabq\, \cU_e|^2}{\theta} \dd \bq - 2 \int_D \theta \varphi |\nabq\,\cU_\eta|^2 \dd \bq, \]
it follows (by dropping the nonnegative second term appearing in the square brackets in the definition of $\mathcal{A}(t,\bx)$ above) that 
\[ \mathcal{A}(t,\bx) \geq \frac{1}{2} \int_D \left[ \frac{\theta|\nabq\varphi|^2}{\varphi} + \frac{\varphi|\nabq\, \cU_e|^2}{\theta} \right] \dd \bq - 4 \int_D \theta \varphi |\nabq\,\cU_\eta|^2 \dd \bq + 2 \int_D \nabq \varphi \cdot \nabq\,\cU_e \dd \bq.
\]
Next, we perform partial integration in the last term on the right-hand side of this inequality using the implicitly imposed boundary condition $\varphi|_{Q \times \partial D}=0$ (cf.~\eqref{eq:bc0}) to deduce that 
\[ \mathcal{A}(t,\bx) \geq \frac{1}{2} \int_D \left[ \frac{\theta|\nabq\varphi|^2}{\varphi} + \frac{\varphi|\nabq\, \cU_e|^2}{\theta} \right] \dd \bq - 4 \int_D \theta \varphi |\nabq\,\cU_\eta|^2 \dd \bq - 2 \int_D \varphi \cdot \Delta_{\bq}\,\cU_e \dd \bq.
\]
As $|\nabq\,\cU_\eta(\bq)| \leq C$ for all $\bq \in \overline D$, and $\theta$ is independent of $\bq$, it follows from \eqref{eq:varphi} that 
\begin{align*} 
\mathcal{A}(t,\bx) &\geq  - C\theta + \frac{1}{2}\int_D \left[ \frac{\theta|\nabq\varphi|^2}{\varphi} + \frac{\varphi|\nabq\, \cU_e|^2}{\theta} \right] \dd \bq - 2 \int_D \varphi \cdot \Delta_{\bq}\,\cU_e \dd \bq\\
&= - C\theta + \frac{1}{2}\int_D \frac{\theta|\nabq\varphi|^2}{\varphi}
\dd \bq + \frac{1}{2}\int_D \frac{\varphi}{\theta}\left(|\nabq\, \cU_e|^2 - 4 \theta \Delta_{\bq}\, \cU_e \right)\! \dd \bq.
\end{align*}
Now, let $\theta_{\mathrm{max}} ~(>\theta_{\mathrm{min}})$ be an arbitrary positive constant, and let 
\[Q_{\theta_{\mathrm{max}}}:= \{(t,\bx) \in \overline Q\,:\, 0 \leq \theta(t,\bx) \leq \theta_{\mathrm{max}}\}.\]
Importantly, as indicated in the previous section, at this point we need to demand that $U_e$ satisfies \eqref{eq:5+} rather than the weaker hypothesis \eqref{eq:5}. Then, because by \eqref{eq:5+}
\[ \lim_{|\bq|\to (\sqrt{b})_{-}} \frac{ \Delta_{\bq}\,\cU_e(\bq)}{|\nabq\, \cU_e(\bq)|^2} =0, \]
it follows that there exists a $\delta= \delta(\theta_{\mathrm{max}}) \in (0,1)$ sufficiently small such that 
\[ |\nabq\, \cU_e|^2 - 4 \theta_{\mathrm{max}} \Delta_{\bq}\, \cU_e = |\nabq\, \cU_e|^2 \left(1 - 4 \theta_{\mathrm{max}} \frac{ \Delta_{\bq}\,\cU_e}{ |\nabq\, \cU_e|^2}\right) \geq \frac{1}{2} |\nabq\, \cU_e|^2 \quad \mbox{for all $\bq \in D \setminus (1-\delta)D$}. \]
Hence, for any $(t,\bx) \in Q_{\theta_{\mathrm{max}}}$, we have that
\begin{align*} 
\mathcal{A}(t,\bx) &\geq - C\theta + \frac{1}{2}\int_D \frac{\theta|\nabq\varphi|^2}{\varphi}
\dd \bq + \frac{1}{2}\int_{(1-\delta)D} \frac{\varphi}{\theta}\left(|\nabq\, \cU_e|^2 - 4 \theta \Delta_{\bq}\, \cU_e \right)\! \dd \bq\\
&\quad\, +  \frac{1}{2}\int_{D\setminus (1-\delta)D} \frac{\varphi}{\theta}\left(|\nabq\, \cU_e|^2 - 4 \theta \Delta_{\bq}\, \cU_e \right)\! \dd \bq\\
 &\geq - C\theta + \frac{1}{2}\int_D \frac{\theta|\nabq\varphi|^2}{\varphi}
\dd \bq + \frac{1}{2}\int_{(1-\delta)D} \frac{\varphi |\nabq\, \cU_e|^2}{\theta} \dd \bq  - 2 \int_{(1-\delta)D} \varphi \Delta_{\bq}\, \cU_e \dd \bq\\
&\quad\, +  \frac{1}{4}\int_{D\setminus (1-\delta)D} \frac{\varphi |\nabq\, \cU_e|^2}{\theta}  \dd \bq\\
& \geq - C\theta + \frac{1}{2}\int_D \frac{\theta|\nabq\varphi|^2}{\varphi}
\dd \bq + \frac{1}{4}\int_{D} \frac{\varphi |\nabq\, \cU_e|^2}{\theta}  \dd \bq  - 2  \int_{(1-\delta)D} \varphi \Delta_{\bq}\, \cU_e \dd \bq\\
& \geq - C\theta + \frac{1}{4} \int_D \left(\frac{\theta|\nabq\varphi|^2}{\varphi} + \frac{\varphi |\nabq\, \cU_e|^2}{\theta}\right)\! \dd \bq  - 2  K(\theta_{\mathrm{max}})\int_{(1-\delta)D} \varphi \dd \bq, 
\end{align*}
where 
\[K(\theta_{\mathrm{max}}):= \max_{\bq \in (1-\delta)\overline{D}}  \Delta_{\bq}\, \cU_e(\bq) < \infty,\] 
thanks to the uniform continuity of $\Delta_{\bq}\,\cU_e$ on compact subsets of $D$. Because, by \eqref{eq:varphi},
\[ 0 \leq \int_{(1-\delta)D} \varphi \dd \bq \leq \int_{D} \varphi \dd \bq \leq \|\varphi_0\|_{L^\infty(\Omega;L^1(D))}\]
and the constant $K(\theta_{\mathrm{max}})$ is positive, we have that 
\begin{align*}
 \mathcal{A}(t,\bx) &\geq -(C\theta +2K(\theta_{\mathrm{max}}) \|\varphi_0\|_{L^\infty(\Omega;L^1(D))}) + \frac{1}{4} \int_D \left(\frac{\theta|\nabq\varphi|^2}{\varphi} + \frac{\varphi |\nabq\, \cU_e|^2}{\theta}\right)\! \dd \bq
\\ & \geq  -(C \theta_{\mathrm{max}} +2K(\theta_{\mathrm{max}}) \|\varphi_0\|_{L^\infty(\Omega;L^1(D))})+ \frac{1}{4} \int_D \left(\frac{\theta|\nabq\varphi|^2}{\varphi} + \frac{\varphi |\nabq\, \cU_e|^2}{\theta}\right)\! \dd \bq 
\end{align*}
for all $(t,\bx) \in Q_{\theta_{\mathrm{max}}}$. In other words, 
\begin{align}\label{eq:ac}
\frac{1}{4} \int_D \left(\frac{\theta|\nabq\varphi|^2}{\varphi} + \frac{\varphi |\nabq\, \cU_e|^2}{\theta}\right)\! \dd \bq \leq \mathcal{A}(t,\bx) + {C}(\theta_{\mathrm{max}}), 
\end{align}
where ${C}(\theta_{\mathrm{max}})$ is a positive constant which depends on $\theta_{\mathrm{max}}$. We return with this bound to the inequality \eqref{eq:Gt} and take 
\[ G(\theta):=1\!\!1_{[0,\theta_{\mathrm{max}}]}(\theta) \quad \mbox{for $\theta \in [0,\infty)$}\]
there. Clearly, $G(\theta(t,\bx)) =
1\!\!1_{[0,\theta_{\mathrm{max}}]}(\theta(t,\bx)) = 1\!\!1_{Q_{\theta_{\mathrm{max}}}}(t,\bx)$ for all $(t,\bx) \in Q$, and therefore \eqref{eq:Gt} and \eqref{eq:ac} imply that
\begin{align*} 
\int_{Q \times D} 1\!\!1_{Q_{\theta_{\mathrm{max}}}} \left(\frac{\theta|\nabq\varphi|^2}{\varphi} + \frac{\varphi |\nabq\, \cU_e|^2}{\theta}\right)\! \dd \bq \dd \bx \dd t &= 
\int_{Q \times D} G(\theta) \left(\frac{\theta|\nabq\varphi|^2}{\varphi} + \frac{\varphi |\nabq\, \cU_e|^2}{\theta}\right)\! \dd \bq \dd \bx \dd t \nonumber \\
&\leq C+ C(\theta_{\mathrm{max}})|Q|,
\end{align*}
In other words, 
\begin{align}\label{eq:up}
\int_{Q_{\theta_{\mathrm{max}}} \times D}  \left(\frac{\theta|\nabq\varphi|^2}{\varphi} + \frac{\varphi |\nabq\, \cU_e|^2}{\theta}\right)\! \dd \bq \dd \bx \dd t  \leq C+ C(\theta_{\mathrm{max}})|Q|.
\end{align}

Applying Lemma \ref{lem:1} with $Q$ there replaced by $Q_{\theta_{\mathrm{max}}}$ and $g(t,\bx) := 1$ for all $(t,\bx) \in Q_{\theta_{\mathrm{max}}}$, we deduce that $\varphi(t,\bx,\cdot)$ and $\varphi(t,\bx,\cdot) U_e'(|\cdot|^2/2)$ have well-defined traces in $L^1(\partial D)$ which are equal to zero for a.e.~$(t,\bx) \in Q_{\theta_{\mathrm{max}}}$.
Because, by \eqref{eq:16}, 
\[\|\theta\|_{L^\infty(0,T;L^1(\Omega))} \leq \mathrm{e}^{T}\left(\frac{1}{2}\|\bv_0\|^2_{L^2(\Omega)} + \|\theta_0\|_{L^1(\Omega)} + \|\,\cU_e \varphi_0 \|_{L^1(\Omega \times D)} + \frac{1}{2}\|\bff\|^2_{L^2(Q)}\right) < \infty,\]
whereby $\theta$ is finite a.e.~on $Q$, in the limit of $\theta_{\mathrm{max}}\to +\infty$ the increasing nested family of sets $\{Q_{\theta_{\mathrm{max}}}\}$  exhausts $Q$. Therefore, $\varphi(t,\bx,\cdot)$ has a well-defined trace $\varphi(t,\bx,\cdot)|_{\partial D} \in L^1(\partial D)$ for a.e.~$(t,\bx) \in Q$ and $\varphi|_{\partial D}$ is equal to $0$ a.e.~on $Q \times\partial D$; also,  $\varphi(t,\bx,\cdot) \nabq\, \cU_e(\cdot) \cdot \boldsymbol{n}_{\bq}|_{\partial D} \in L^1(\partial D)$ for a.e.~$(t,\bx) \in Q$, and  it, too, vanishes on $\partial D$ for a.e.~$(t,\bx) \in Q$. We draw the reader's attention to the fact that, although $\nabq\, \cU_e \in C^1(D;\mathbb{R}^d)$, the second assertion in the previous sentence does not directly follow from the fact that $\varphi|_{\partial D} = 0$, as $\nabq\,\cU_e(\bq)\cdot\boldsymbol{n}_{\bq} = U_e'(|\bq|^2/2)|\bq|= \sqrt{b}\, U_e'(|\bq|^2/2) = \infty$ when $\bq \in \partial D$.

\section{Weak sequential stability}\label{sec:3}
We begin by stating the definition of a weak solution to the problem under consideration. 

\begin{definition}[Weak solution]
Suppose that $\Omega \subset \mathbb{R}^d$ is a bounded, open, simply connected set in $\mathbb{R}^d$, $d \in \{2,3\}$, with Lipschitz continuous boundary $\partial\Omega$, and $T>0$. Assume that the initial data $\bv_0$, $\theta_0$,  $\varphi_0$  and the source term $\bff$ satisfy the following properties:
\begin{itemize}
    \item $\bv_0 \in L^2(\Omega;\mathbb{R}^d)$, $\mathrm{div}_{\bx} \bv_0 = 0$ in the sense of distributions on $\Omega$ and $\bv_0 \cdot \boldsymbol{n}_{\bx}|_{\partial \Omega} = 0$ in $W^{-1/2,2}(\partial \Omega)$; 
\item  $\theta_0 \in L^1(\Omega)$, $\theta_0 \geq \theta_{\mathrm{min}}>0$;
\item $\varphi_0 \in L^1 \log L^1(\Omega \times D) \cap L^\infty(\Omega;L^1(D))$, $\mathcal{U}_e \varphi_0 \in L^1(\Omega \times D)$, $\varphi_0 \geq 0$;
\item $\bff \in L^2(Q)$.
\end{itemize}
Let
\[ \bF:= \nabq\, \cU_e + \theta \nabq\, \cU_\eta \qquad \mbox{and} 
\qquad
\bj_{\varphi,\bq} := - 4\boldsymbol{F}\varphi - 4\theta \nabq \varphi.
\]
Suppose further that $\beta>5/6$ and $p:= 2 - \frac{6}{5+3\beta}$. We shall say that a triple of functions $(\bv,\theta,\varphi)$ is a weak solution to the problem under consideration provided that:
\begin{itemize}
\item $\bv \in L^\infty(0,T;L^2_{0,\mathrm{div}}(\Omega;\mathbb{R}^d)) \cap 
L^p(0,T; W^{1,p}_{0,\mathrm{div}}(\Omega; \mathbb{R}^d))$,  $\partial_t \bv \in L^s(0,T;W^{-1,s'}_{0,\mathrm{div}}(\Omega;\mathbb{R}^d))$ for some $s>1$, where $s':=s/(s-1)$ is the H\"older conjugate of $s$;


\item $\theta \in L^\infty(0,T;L^1(\Omega;\mathbb{R}_{> 0})) \cap L^{\frac{2}{3} + \beta}(Q)$, $\theta^{\beta/2} \in L^2(0,T;W^{1,2}(\Omega))$ and
$\log \theta \in L^\infty(0,T;L^1(\Omega))
\cap L^2(0,T;W^{1,2}(\Omega))$;
%
\item 
\[ \frac{\nu(\theta)|\mathbb{D}(\bv)|^2}{\theta} \in L^1(Q), \quad \frac{\kappa(\theta)|\nabx \theta|^2}{\theta^2} \in L^1(Q),\]
\[ \frac{\theta}{\varphi}|\nabx \varphi|^2 \in L^1(Q\times D),\quad \frac{1}{\theta \varphi}|\theta \nabq \varphi + \theta \varphi \nabq\, \mathcal{U}_\eta + \varphi \nabq \,\mathcal{U}_e|^2 \in L^1(Q \times D),\]
\[ 2 \nu(\theta) \mathbb{D}(\bv) - 2 \theta n_{\mathrm{P}} \mathbb{I} + \int_D (\boldsymbol{F} \otimes \bq) \varphi \dd \bq \in L^p(Q;\mathbb{R}^{d \times d});\]
\item
$\varphi \!\in\! L^\infty(0,T; L^1 \log L^1(\Omega \times D;\mathbb{R}_{\geq 0}))$, $\mathcal{U}_e \varphi \in L^\infty(0,T; L^1(\Omega \times D))$, $\partial_t \varphi \!\in  \![C([0,T];C^1(\overline{\Omega \times D}))]^*$;
\end{itemize}
furthermore, the velocity field $\bv$ satisfies the balance of linear momentum equation in the following weak sense:
\begin{align*}
- \int_Q \bv\cdot \partial_t \bw \dd \bx \dd t- \int_Q (\bv \otimes \bv) : \nabx \bw \dd \bx \dd t + \int_Q \bSS:\mathbb{D}(\bw)  \dd \bx \dd t = \int_\Omega \bv_0 \cdot \bw(0) \dd \bx + \int_Q \bff \cdot \bw \dd \bx \dd t
\end{align*}
for all test functions $\bw \in C^1_0([0,T); W^{1,s}_{0,\mathrm{div}}(\Omega;\mathbb{R}^d))$, where 
\[ \bSS = 2\nu(\theta) \bDD(\bv) + \int_D (\boldsymbol{F} \otimes \bq) \varphi  \dd \bq;  \]
for all positive integers $k$ the temperature $\theta$ satisfies the following renormalized variational inequality:
\begin{align*}
&-\int_Q T_{k}(\theta) \,\partial_t \psi \dd \bx \dd t - \int_Q \bv\, T_{k}(\theta) \cdot \nabx \psi \dd \bx \dd t + \int_Q \kappa(\theta) \nabx T_{k}(\theta)\cdot \nabx \psi \dd \bx \dd \\
& \quad \geq \int_\Omega T_{k}(\theta_0)\,\psi(0) \dd \bx   + \int_Q 2\nu(\theta) T_{k}'(\theta)|\bDD(\bv)|^2\psi  \dd \bx \dd t\\
&\qquad +\int_Q \left(T_{k}'(\theta) \bDD(\bv): \int_D \theta (\nabq\, \cU_\eta \otimes \bq) \varphi \dd \bq \right) \psi \dd \bx \dd t \\
&\qquad+ \int_Q \left(T_{k}'(\theta) \int_D \theta\nabq \varphi \cdot \nabq\, \cU_e + \theta \varphi \nabq\, \cU_\eta \cdot \nabq\, \cU_e + |\nabq\, \cU_e|^2 \varphi \dd \bq\right) \psi \dd \bx \dd t
\end{align*}
for all test functions $\psi \in C^1_0([0,T); W^{1,z}(\Omega))$ with $z>d$ such that $\psi \geq 0$ on $\overline{Q}$, where 
\[ T_k(s):= \left\{ \begin{array}{cl} \frac{s}{|s|} \max\{k,|s|\} & \mbox{for $s \neq 0$},\\
0 & \mbox{for $s=0$};\end{array} \right. \]
and the probability density function $\varphi$ satisfies the Fokker--Planck equation in the following weak sense: 
\begin{align*}
& - \int_{Q \times D} \varphi\, \partial_t\zeta \dd \bq \dd \bx \dd t - \int_{Q \times D} \bv \varphi \cdot \nabx \zeta \dd \bq \dd \bx \dd t  +  \int_{Q \times D}\theta \nabx \varphi \cdot \nabx \zeta \dd \bq \dd \bx \dd t \\
&\quad - \int_{Q \times D} ((\nabx \bv) \bq \varphi + \bj_{\varphi,\bq})\cdot \nabq \zeta \dd \bq \dd \bx \dd t = \int_{\Omega \times D} \varphi_0\, \zeta(0) \dd \bq \dd \bx
\qquad \forall\, \zeta \in C^1_0([0,T);C^1(\overline{\Omega \times D})).
\end{align*}
Finally, a weak solution is understood to be energy-dissipative in the sense that for a.e.~$t \in (0,T)$ it satisfies the following energy inequality, with $\mathcal{H}$ and $\mathcal{F}$ as defined in \eqref{eq:HF}:
\begin{align}
\begin{aligned}\label{eq:energy+}
&\int_\Omega \bigg[\frac{1}{2}|\bv(t)|^2 +  \mathcal{H}(\theta(t)) \bigg] \dd \bx  + \int_{\Omega \times D} [\,\cU_e \varphi(t) + \cU_{\eta} \varphi(t) + \mathcal{F}(\varphi(t))]\dd \bq \dd \bx \\
&\qquad+ \int_0^t \! \int_\Omega \bigg[\frac{2\nu(\theta)|\bDD(\bv)|^2}{\theta} + \frac{\kappa(\theta) |\nabx \theta|^2}{\theta^2}\bigg]\dd \bx \dd s\\ 
&\qquad + \int_0^t \int_\Omega \bigg[\theta \int_D \frac{|\nabx \varphi|^2}{\varphi}\dd \bq + \int_D \frac{4}{\theta \varphi}\left| \theta \nabq\varphi + \theta \varphi \nabq\,\cU_\eta + \varphi \nabq \,\cU_e\right|^2\! \dd \bq\bigg] \dd \bx \dd s\\
& \leq  \int_\Omega \bigg[\frac{1}{2}|\bv_0|^2 +  \mathcal{H}(\theta_0) \bigg] \dd \bx  + \int_{\Omega \times D} [\,\cU_e \varphi_0 + \cU_{\eta} \varphi_0 + \mathcal{F}(\varphi_0)]\dd \bq \dd \bx + \int_0^t \int_\Omega \bff(s) \cdot \bv(s)\dd \bx \dd s.
\end{aligned}
\end{align}
%
\end{definition}
%
%
%

The central result of the paper is the following theorem, 
which guarantees weak sequential stability of the nonisothermal kinetic model under consideration.

\begin{theorem}\label{thm:1}
Under the stated assumptions on the data, sequences of smooth solutions to the initial-boundary-value problem \eqref{eq:m1}--\eqref{eq:mbct} for the system of nonlinear partial differential equations under consideration, which satisfy the energy equality \eqref{eq:energy}, converge to an energy-dissipative weak solution of the problem. 
\end{theorem}

The rest of the paper is devoted to the proof of this theorem. Because the proof is long, we draw the reader's attention to the fact that, even though the argument is seemingly interrupted by titles of subsections, the proof runs through several subsections and ends at the end of Section 3. 

\bigskip

\noindent
\textbf{Proof.} Suppose that we have a sequence of functions $((\bv^n, \theta^n,\varphi^n  ))_{n=1}^\infty$ with $\bv^n \in C^{1,2}([0,T] \times \overline\Omega;\mathbb{R}^d)$, $\theta^n \in C^{1,2}([0,T] \times \overline \Omega)$ and $\varphi^n \in C^{1,2,2}([0,T] \times \overline \Omega \times \overline D)$, satisfying the inequality \eqref{eq:energy+}, the equations \eqref{eq:m1}, \eqref{eq:m2}, 
\eqref{eq:m3+} and \eqref{eq:m4} subject to the boundary and initial conditions \eqref{eq:mbcv}, \eqref{eq:mbcx+}, \eqref{eq:mbcq+}, \eqref{eq:mic}, and \eqref{eq:mbct}, i.e., for all $n=1,2,\ldots$,
\[ \bv^n = \boldsymbol{0} \;\; \mbox{on $(0,T) \times \partial\Omega$},\quad \bv^n(0,\bx) = \bv_0^n(x)\;\; \mbox{for $\bx \in \Omega$}, \]
\[ \kappa(\theta^n) \nabx \theta^n \cdot \boldsymbol{n}_{\bx} = 0, \;\; \mbox{on $(0,T) \times \partial\Omega$},\quad 
\theta^n(0,\bx) = \theta_0^n(\bx)\;\; \mbox{for $\bx \in \Omega$},\]
\[ \left\{ \begin{array}{c}  \left(\bv^n \varphi^n - [\theta^n]_{+} \nabx \varphi^n\right) \cdot \boldsymbol{n}_{\bx} = 0 \;\; \mbox{on $(0,T) \times \partial \Omega \times D$}, \\
((\nabx \bv^n) \bq \varphi^n + \bj_{\varphi,\bq}^n)\cdot \boldsymbol{n}_{\bq} =0\;\; \mbox{on $(0,T) \times \Omega \times \partial D$}, \end{array} \right. 
\quad\!\! \varphi^n(0,\bx,\bq) = \varphi_0^n(\bx,\bq) \;\; \mbox{for $(\bx,\bq) \in \Omega \times D$},
\]
where $(\bv_0^n)_{n=1}^\infty$, $(\theta_0^n)_{n=1}^\infty$, $(\varphi_0^n)_{n=1}^\infty$ are sequences of  functions defined on $\overline\Omega$, $\overline\Omega$ and $\overline\Omega \times \overline D$, respectively, such that $\bv^n_0 \in C^2(\overline\Omega;\mathbb{R}^d)$,  $\theta^n_0 \in C^2(\overline\Omega)$, 
$\varphi^n_0 \in C^{2,2}(\overline\Omega \times \overline D)$,
$\divx \bv_0^n =0$ on $\Omega$ and $\bv_0^n \cdot \boldsymbol{n}_{\bx} = 0$ on $\partial \Omega$, $\theta_0^n \geq \theta_{\mathrm{min}}>0$ on $\overline\Omega$, and  $\varphi_0^n \geq 0$ on $\overline \Omega \times \overline D$ for all $n=1,2,\ldots $; with $\bv^n_0 \to \bv_0$ in $L^2(\Omega;\mathbb{R}^d)$,  $\theta^n_0 \to \theta_0$ in $L^1(\Omega)$, $\varphi_0^n \to \varphi_0$ in $L^1(\Omega \times D)$, $\cU_e \varphi_0^n \to \cU_e \varphi_0$ in $L^1(\Omega \times D)$, $\mathcal{F}(\varphi_0^n) \to \mathcal{F}(\varphi_0)$ in $L^1(\Omega \times D)$,  as $n \to \infty$, and
\begin{align*}
\bj_{\varphi,\bq}^n &:= - \boldsymbol{F}^n[\varphi^n]_{+} - [\theta^n]_{+} \nabq \varphi^n \qquad \mbox{with} \qquad 
\boldsymbol{F}^n := \nabq\, \cU_e + \theta^n \nabq \,\cU_\eta.
\end{align*}

As in Step 0 in the previous section, the assumed nonnegativity of $\varphi_0^n$ on $\overline \Omega \times \overline D$ implies the nonnegativity of $\varphi^n$ on $[0,T] \times \overline \Omega \times \overline D$ and, as in Step 1, the assumed positive lower bound on $\theta_0^n$ on $\overline \Omega$ implies the positivity of $\theta^n$ on $[0,T] \times \overline \Omega$. Consequently, $[\varphi^n]_{+}$ and $[\theta^n]_{+}$
appearing in the equations above can be replaced by
$\varphi^n$ and $\theta^n$; therefore, we shall remove the now redundant symbol $[\cdot]_{+}$ in what follows. 

Because $\bv^n=\mathbf{0}$ on $\partial\Omega$ and $\theta^n$ is strictly positive  on $[0,T] \times \overline{\Omega}$,  it follows that the zero normal flux boundary condition for $\varphi^n$ on $(0,T)\times \partial\Omega \times D$ collapses to the following homogeneous Neumann boundary condition: 
\[ \nabx \varphi^n \cdot \boldsymbol{n}_{\bx} = 0 \quad \mbox{on $(0,T) \times \partial \Omega \times D$.}\]
We note further that because 
\begin{align}\label{eq:bnj}
\begin{aligned}
\bj_{\varphi,\bq}^n &:= - \boldsymbol{F}^n\varphi^n - \theta^n \nabq \varphi^n = -(\nabq\, \cU_e + \theta^n \nabq \,\cU_\eta)\varphi^n -\theta^n \nabq \varphi^n\\
&= - \theta^n \left[\varphi^n\nabq\left( \frac{\,\cU_e + \theta^n \cU_\eta}{\theta^n}\right) + \nabq \varphi^n\right] =- \theta^n 
 \left(\nabq \left(\mathrm{e}^{\frac{\,\cU_e + \theta^n \cU_\eta}{\theta^n}}\varphi^n\right)\,\mathrm{e}^{-\frac{\cU_e +\theta^n \cU_\eta}{\theta^n}}\right),
\end{aligned}
\end{align}
the zero normal flux boundary condition \eqref{eq:mbcq} imposed for $\varphi^n$ on $(0,T) \times \Omega \times \partial D$ can be restated as follows:
\[ \left[(\nabx \bv^n)\bq \varphi^n -\theta^n \left(\nabq \left(\mathrm{e}^{\frac{\,\cU_e + \theta^n \cU_\eta}{\theta^n}}\varphi^n\right)\,\mathrm{e}^{-\frac{\cU_e +\theta^n \cU_\eta}{\theta^n}}\right)\right]\cdot \boldsymbol{n}_{\bq} = 0 \quad \mbox{on $(0,T) \times \Omega \times \partial D$}.\]

We shall use the bounds established in Steps 1--15 in Section \ref{sec:2} to pass to the respective limits in the sequences $(\bv^n)_{n=1}^\infty$, $(\theta^n)_{n=1}^\infty$ and $(\varphi^n)_{n=1}^\infty$.

\subsection{Limit passage in the balance of linear momentum equation}

The a priori bounds \eqref{eq:varphi}--\eqref{eq:26} and \eqref{eq:28} imply the existence of weakly convergent subsequences (not indicated) such that, as $n \to \infty$,
\begin{alignat}{2}
\begin{aligned}\label{eq:v}
\bv^n &\rightharpoonup^\ast \bv\quad && \mbox{weakly-$\ast$ in $L^\infty(0,T;L^2(\Omega;\mathbb{R}^d))$}, \\
\bv^n &\rightharpoonup\bv\quad && \mbox{weakly in $L^{\frac{5}{3}\frac{2+3\beta}{5+3\beta}}(Q;\mathbb{R}^d))$},\\
\nabx\bv^n &\rightharpoonup \nabx\bv\quad && \mbox{weakly in $L^{\frac{4+6\beta}{5+3\beta}}(Q;\mathbb{R}^{d\times d})$}, \\
\bSS^n &\rightharpoonup \overline{\bSS}\quad && \mbox{weakly in $L^p(Q;\mathbb{R}^{d\times d})$},\quad p:=2- \frac{6}{5+3\beta} = \frac{4+6\beta}{5+3\beta}, 
\end{aligned}
\end{alignat}
with $\beta>\frac{5}{6}$, where
\[ \bSS^n:=2\nu(\theta^n) \bDD(\bv^n) + \int_D (\boldsymbol{F^n} \otimes \bq) \varphi^n  \dd \bq = \nu(\theta^n) \bDD(\bv^n) + \int_D ((\nabq\, \cU_e + \theta^n \nabq \,\cU_\eta) \otimes \bq) \varphi^n  \dd \bq,  \]
and $\overline{\bSS}$ is a weak limit, to be identified. From the balance of linear momentum equation \eqref{eq:m1} satisfied by $\bv^n$ one has a uniform bound on $\partial_t \bv^n$ in $L^s(0,T;W^{-1,s'}_{0,\mathrm{div}}(\Omega;\mathbb{R}^d))$ for some $s>1$, where $s':=\frac{s}{s-1}$, which then implies the existence of a subsequence (not indicated) such that, as $n \to \infty$,
\[\partial_t\bv^n \rightharpoonup \partial_t\bv\quad  \mbox{weakly in $L^s(0,T;W^{-1,s'}_{0,\mathrm{div}}(\Omega;\mathbb{R}^d))$ for some $s>1$}.\]
Thanks to the compact embedding of $W^{1,\frac{4+6 \beta}{5+3\beta}}(\Omega;\mathbb{R}^d))$ in $L^{r}(\Omega;\mathbb{R}^d)$ with $1 \leq r < \frac{6(2+3\beta)}{11+3\beta}$ for $d=2,3$, and noting that $6>\frac{6(2+3\beta)}{11+3\beta}>\frac{10}{3}\frac{2+3 \beta}{5+3\beta}>2$ for $\beta>\frac{5}{6}$, the Aubin--Lions lemma implies that, as $n \to \infty$,
\begin{align}\label{eq:va} 
\bv^n \rightarrow \bv \quad \mbox{strongly in $L^r(Q;\mathbb{R}^d)$ for $1 \leq r \leq \frac{10}{3}\frac{2+3 \beta}{5+3\beta}$}
\end{align}
for $d=2,3$ and $\beta > \frac{5}{6}$, and in particular, as $n \to \infty$,
\begin{align}\label{eq:va2}\bv^n \rightarrow \bv \quad \mbox{strongly in $L^2(Q;\mathbb{R}^d)$ for $d=2,3$}.
\end{align}
These convergence results then enable us to pass to the limit in the sequence of balance of linear momentum equations as $n \rightarrow \infty$ to deduce that, for a.e.~$t \in (0,T)$, 
\begin{align}\label{eq:E1F}
- \int_Q \bv\cdot \partial_t \bw \dd \bx \dd t- \int_Q (\bv \otimes \bv) : \nabx \bw \dd \bx \dd t + \int_Q \overline{\bSS}:\mathbb{D}(\bw)  \dd \bx \dd t = \int_\Omega \bv_0 \cdot \bw(0) \dd \bx + \int_Q \bff \cdot \bw \dd \bx \dd t
\end{align}
for all test functions $\bw \in C^1_0([0,T); W^{1,s}_{0,\mathrm{div}}(\Omega;\mathbb{R}^d))$, 
for some $s>1$, where the weak limit $\overline\bSS$ is yet to be identified.

\subsection{Limit passage in the sequence of absolute temperatures}

For the sequence $(\theta^n)_{n=1}^\infty$, we have from \eqref{eq:22}, \eqref{eq:17} and \eqref{eq:21} that, as $n \to \infty$,
\begin{alignat}{2}\label{eq:**}
\begin{aligned}
\theta^n &\rightharpoonup \theta\quad && \mbox{weakly in $L^{\frac{2}{3}+\beta}(Q)$}, \\
\log\theta^n &\rightharpoonup \overline{\log \theta}\quad && \mbox{weakly in $L^2(0,T;W^{1,2}(\Omega))$}, \\
(\theta^n)^{\frac{\beta}{2}} &\rightharpoonup \overline{\theta^{\frac{\beta}{2}}}\quad && \mbox{weakly in $L^2(0,T;W^{1,2}(\Omega))$},
\end{aligned}
\end{alignat}
where $\overline{\log \theta}$ and $\overline{\theta^{\frac{\beta}{2}}}$
are weak limits, to be identified. 

For $k \geq 1$ we introduce the real-valued continuous piecewise-linear cutoff function $T_k \in C^{0,1}(\mathbb{R})$ defined by 
\[ T_k(s):= \left\{ \begin{array}{cl} \frac{s}{|s|} \max\{k,|s|\} & \mbox{for $s \neq 0$},\\
0 & \mbox{for $s=0$}. \end{array} \right. \]
Clearly, $T_k(s)=s$ for $s \in [-k,k]$ and $T_k(s) = k\,\mathrm{sign}(s)$ for $|s| \geq k$. Let $T_{k,\varepsilon} \in C^2(\mathbb{R})$ be a mollification of $T_k$, with $\varepsilon \in (0,1)$, such that 
\begin{alignat*}{2}
T_{k,\varepsilon}(s) &= s \quad && \mbox{for $|s| \leq k$},\\
T_{k,\varepsilon}(s) &= \frac{s}{|s|}\left(k + \frac{1}{2}\varepsilon\right) \quad && \mbox{for $|s| \geq k + \varepsilon$},\\
T_{k,\varepsilon}(s)&=-\frac{s}{2\varepsilon|s|}\left(s^2 -2\frac{s}{|s|}(k+\varepsilon)s + k^2\right)\quad && \mbox{for $k \leq |s| \leq k+ \varepsilon$}.
\end{alignat*}
Thus,
\[ 0 \leq T'_{k,\varepsilon}(s) \leq 1, \quad sT''_{k,\varepsilon}(s) \leq 0 \quad \mbox{and} \quad |T_{k,\varepsilon}''(s)| \leq \frac{1}{\varepsilon}\quad \mbox{for all $s \in \mathbb{R}$}.\]

We multiply by $T_{k,\varepsilon}'$ the temperature equation satisfied by $\theta^n$ to obtain a renormalized version of the equation; hence, 
\begin{align}\label{eq:ren}
\begin{aligned}
\partial_t T_{k,\varepsilon}(\theta^n) &+ \divx(\bv^n\, T_{k,\varepsilon}(\theta^n)) - \divx(\kappa(\theta^n) \nabx T_{k,\varepsilon}(\theta^n))
= 2\nu(\theta^n) T_{k,\varepsilon}'(\theta^n)|\bDD(\bv^n)|^2\\
&+T_{k,\varepsilon}'(\theta^n) \bDD(\bv^n): \int_D \theta^n (\nabq\, \cU_\eta \otimes \bq) \varphi^n \dd \bq - T_{k,\varepsilon}''(\theta^n)\kappa(\theta^n)|\nabx\theta^n|^2\\
&+ T_{k,\varepsilon}'(\theta^n) \int_D \theta^n \nabq \varphi^n \cdot \nabq\, \cU_e + \theta^n  \varphi^n \nabq\, \cU_\eta \cdot \nabq\, \cU_e + |\nabq\, \cU_e|^2 \varphi^n \dd \bq.
\end{aligned}
\end{align}
Since $T_{k,\varepsilon}'$ has compact support $[-k-\varepsilon,k+\varepsilon] \subset \mathbb{R}$, we deduce from \eqref{eq:ren} that 
\begin{align}
\label{eq:*} 
\|\partial_t T_{k,\varepsilon}(\theta^n)\|_{\mathcal{M}(0,T;(W^{1,z}(\Omega))^\ast)}\leq C(k,\varepsilon)\quad \mbox{for all $z \in (d,\infty)$},
\end{align}
where, for a Banach space $X$, $\mathcal{M}(0,T;X^*)$ denotes the continuous dual space of $C([0,T];X)$, that is, the Banach space of $X^*$-valued (finite, signed) Radon measures on $[0,T]$.

To show that the bound \eqref{eq:*} holds, we test \eqref{eq:ren} with $\chi \in C([0,T];W^{1,z}(\Omega))$, where $z>d$, and move the second and third term from the left-hand side of the resulting equality to the right-hand side. Hence, 
\begin{align*}
&\langle \partial_t T_{k,\varepsilon}(\theta^n),\chi \rangle = \int_Q \bv^n\, T_{k,\varepsilon}(\theta^n)\cdot \nabx \chi \dd \bx \dd t \\
&\quad - \int_Q \kappa(\theta^n) \nabx T_{k,\varepsilon}(\theta^n)\cdot \nabx \chi \dd \bx \dd t\\
&\quad + \int_Q 2\nu(\theta^n) T_{k,\varepsilon}'(\theta^n)|\bDD(\bv^n)|^2 \chi \dd \bx \dd t\\
&\quad + \int_Q T_{k,\varepsilon}'(\theta^n) \bDD(\bv^n): \left(\int_D \theta^n (\nabq\, \cU_\eta \otimes \bq) \varphi^n \dd \bq\right) \chi \dd \bx \dd t \\
&\quad - \int_Q  T_{k,\varepsilon}''(\theta^n)\kappa(\theta^n)|\nabx\theta^n|^2 \chi \dd \bx \dd t\\
&\quad + \int_Q T_{k,\varepsilon}'(\theta^n) \left(\int_D \theta^n \nabq \varphi^n \cdot \nabq\, \cU_e + \theta^n  \varphi^n \nabq\, \cU_\eta \cdot \nabq\, \cU_e + |\nabq\, \cU_e|^2 \varphi^n \dd \bq\right) \chi  \dd \bx \dd t\\
&=:\mathrm{T}_1 + \mathrm{T}_2 + \mathrm{T}_3 + \mathrm{T}_4+ \mathrm{T}_5 + \mathrm{T}_6.
\end{align*}
Then, by \eqref{eq:v}$_1$ and since $|T_{k,\varepsilon}(s)| \leq k+\frac{1}{2}\varepsilon< k+1$ for all $s \in \mathbb{R}$, we have 
\[ |\mathrm{T}_1| \leq \|\bv^n\|_{L^\infty(0,T;L^2(\Omega))}\|T_{k,\varepsilon}(\theta^n)\|_{L^\infty(Q)}\|\nabx \chi\|_{L^1(0,T;L^2(\Omega))} \leq C(k)\|\nabx \chi\|_{L^1(0,T;L^2(\Omega))}.\]
To bound the term $\mathrm{T}_2$, we note that $\kappa(\theta^n)\nabx T_{k,\varepsilon}(\theta^n) =\kappa(\theta^n) \theta^n T'_{k, \varepsilon}(\theta^n) \nabx \log \theta^n$. Thus, by \eqref{eq:**}$_2$ and because $T_{k,\varepsilon}'$ has compact support $[-k-\varepsilon,k+\varepsilon]$, we have the following:
\[ |\mathrm{T}_2| \leq C(k)\|\nabx \chi\|_{L^2(0,T;L^2(\Omega))}. \]
To deal with the term $\mathrm{T}_3$, we note that $\nu(\theta^n)T_{k,\varepsilon}'(\theta^n) |\mathbb{D}(\bv^n)|^2 =\left(T_{k,\varepsilon}'(\theta^n) \theta^n \right)\frac{\nu(\theta^n) |\mathbb{D}(\bv^n)|^2}{\theta^n}$. Thus, again using that $T_{k,\varepsilon}'$ has compact support in conjunction with the bound on the third term from the second line of \eqref{eq:17}, with $\bv$ and $\theta$ there replaced by $\bv^n$ and $\theta^n$, we have
\[ |\mathrm{T}_3| \leq C(k)\|\chi\|_{C([0,T];L^\infty(\Omega))}.\]
To bound the term $\mathrm{T}_4$, we write $T'_{k,\varepsilon}(\theta^n) \mathbb{D}(\bv^n) \theta^n = (T'_{k,\varepsilon}(\theta^n) (\theta^n)^{\frac{3}{2}}) (\mathbb{D}(\bv^n)/\sqrt{\theta^n})$, use that $T_{k,\varepsilon}'$ has compact support and $\nabq\,\cU_\eta \otimes \bq \in C(\overline{D};\mathbb{R}^{d \times d})$, together with the bound on the third term of the second line of \eqref{eq:17}, with $\bv$ and $\theta$ replaced by $\bv^n$ and $\theta^n$, and the bound on the first term on the left-hand side of inequality \eqref{eq:26}, with $\varphi$ replaced by $\varphi^n$; thus, 
\[ |\mathrm{T}_4| \leq C(k) \|\chi\|_{L^2(0,T;L^2(\Omega))}.\]
For the term $\mathrm{T}_5$, we rewrite $T_{k,\varepsilon}''(\theta^n)\kappa(\theta^n)|\nabx\theta^n|^2 = T_{k,\varepsilon}''(\theta^n)\kappa(\theta^n)(\theta^n)^2|\nabx\log \theta^n|^2$ and use \eqref{eq:**}$_2$ to find that
\[ |\mathrm{T}_5|\leq \frac{C(k)}{\varepsilon}\|\chi\|_{C([0,T];L^\infty(\Omega))}.\]
Finally, to deal with the term $\mathrm{T}_6$, we note that
\begin{align}\label{eq:*x}
\begin{aligned}
&T_{k,\varepsilon}'(\theta^n) \left(\int_D \theta^n \nabq \varphi^n \cdot \nabq\, \cU_e + \theta^n  \varphi^n \nabq\, \cU_\eta \cdot \nabq\, \cU_e + |\nabq\, \cU_e|^2 \varphi^n \dd \bq\right)\chi\\
&\quad = T_{k,\varepsilon}'(\theta^n) \theta^n \left(\int_D \frac{\sqrt{\varphi^n}\,\nabq\,\cU_e}{\sqrt{\theta^n}} \cdot \frac{\theta^n \nabq \varphi^n + \theta^n  \varphi^n \nabq\, \cU_\eta + \nabq\, \cU_e \varphi^n}{\sqrt{\theta^n}{\sqrt{\varphi^n}}}\dd \bq\right)\chi.
\end{aligned}
\end{align}
We recall the bound that appears
in the third line of \eqref{eq:17} with $\theta$ and $\varphi$ replaced by $\theta^n$ and $\varphi^n$, the bound \eqref{eq:up} with $\theta$ and $\varphi$ replaced by $\theta^n$ and $\varphi^n$ and $\theta_{\mathrm{max}}$ chosen to be equal to $k + \varepsilon$, and note that for $(t,\bx) \in Q \setminus Q_{\theta_{\mathrm{max}}} = Q\setminus Q_{k+\varepsilon}$ the expression appearing in the last line above is equal to zero because $T'_{k,\varepsilon}(\theta^n(t,\bx)) = 0$ for all such $(t,\bx)$, whereby the integral of the expression \eqref{eq:*x} over $Q$ is equal to its integral over $Q_{\theta_{\mathrm{max}}} = Q_{k+\varepsilon}$. Therefore, integrating \eqref{eq:*x} over $Q$ and applying H\"{o}lder's inequality, we find that 
\[ |\mathrm{T}_6| \leq C(k)\|\chi\|_{C([0,T];L^\infty(\Omega))}.\]
It remains to collect the bounds on the terms $\mathrm{T}_1,\ldots, \mathrm{T}_6$ and note that $C([0,T];W^{1,z}(\Omega))$ is, for $z>d$, continuously embedded in $L^1(0,T;W^{1,2}(\Omega))$, $L^2(0,T;W^{1,2}(\Omega))$, $C([0,T];L^\infty(\Omega))$, and $L^2(0,T;L^2(\Omega))$, and therefore we have
\[ |\langle \partial_t T_{k,\varepsilon}(\theta^n),\chi \rangle| \leq 
|\mathrm{T}_1| + |\mathrm{T}_2| + |\mathrm{T}_3| + |\mathrm{T}_4|+ |\mathrm{T}_5| + |\mathrm{T}_6| \leq C(k,\varepsilon) \|\chi\|_{C([0,T];W^{1,z}(\Omega))},\]
which directly implies \eqref{eq:*}.

In addition, 
\[ \int_Q |\nabx T_{k,\varepsilon}(\theta^n)|^2 \dd \bx \dd t \leq C(k) \int_Q \frac{|\nabx\theta^n|^2}{(\theta^n)^2} \dd \bx \dd t \leq C(k).\]
Hence, 
\begin{align}\label{eq:*b} 
\nabx T_{k,\varepsilon}(\theta^n) \rightharpoonup \overline{\nabx T_{k,\varepsilon}(\theta)}\quad \mbox{weakly in $L^2(0,T;L^2(\Omega;\mathbb{R}^d))$},
\end{align}
where $\overline{\nabx T_{k,\varepsilon}(\theta)}$ is a weak limit, to be identified. 

Now, \eqref{eq:*} and \eqref{eq:*b}, via a generalization of the Aubin--Lions lemma (cf.~Corollary~7.9 in \cite{Roubicek}), imply that as $n \rightarrow \infty$,
\begin{align}\label{eq:***} 
T_{k,\varepsilon}(\theta^n) \rightarrow \overline{T_{k,\varepsilon}(\theta)}\quad \mbox{strongly in $L^2(0,T;L^2(\Omega))$},
\end{align}
where $\overline{T_{k,\varepsilon}(\theta)}$ is a strong limit, to be identified. 

We note that 
\[ \int_Q |\theta^n - T_{k,\varepsilon}(\theta^n)|\dd \bx \dd t \leq \int_{\{(t,\bx) \in Q\, :\, \theta^n(t,\bx)  \geq k\}} |\theta^n| \dd \bx \dd t \leq \frac{C}{k^\alpha}\quad \mbox{with $\alpha:=\beta - \frac{1}{3}>0$},\]
where the second inequality is a consequence of \eqref{eq:22} applied to $\theta^n$ (recall that $\frac{2}{3}+\beta>1$). Therefore, by the triangle inequality, 
\[ \int_Q |\theta^m - \theta^n| \dd \bx \dd t \leq \int_Q |T_{k,\varepsilon}(\theta^m) - T_{k,\varepsilon}(\theta^n)| \dd \bx \dd t + \frac{C}{k^\alpha} \quad \mbox{for all $m,n \in \mathbb{N}$}.\]
Let us now choose a real number $\delta>0$ and then fix $k=k(\delta)$ such that
\[\frac{C}{k^\alpha} \leq \frac{\delta}{2},\]
where $C$ is the positive constant from the above inequality. According to \eqref{eq:***} the sequence $(T_{k,\varepsilon}(\theta^n))_{n=1}^\infty$ is (strongly) convergent in $L^2(0,T;L^2(\Omega)) = L^2(Q)$ and, therefore, because $Q=(0,T) \times \Omega$ is a bounded set, it is also strongly convergent in $L^1(Q)$. Thus, with $k=k(\delta)$ fixed as above, there exists a positive integer $N_0=N_0(\delta)$ such that 
\[ \int_Q |T_{k,\varepsilon}(\theta^m) - T_{k,\varepsilon}(\theta^n)| \dd \bx \dd t \leq \frac{\delta}{2}\quad \forall\,m,n \geq N_0.\]
Hence, we have shown that for any $\delta>0$ there exists a positive integer $N_0=N_0(\delta)$ such that 
\[ \int_Q |\theta^m - \theta^n| \dd \bx \dd t \leq \delta \quad \forall\, m,n \geq N_0.\]
In other words, $(\theta^n)_{n=1}^\infty$ is a Cauchy sequence in $L^1(Q)$, and therefore, thanks to the completeness of $L^1(Q)$, as $n \to \infty$,
\begin{align} \label{eq:str}
\theta^n \rightarrow \theta \quad \mbox{strongly in $L^1(Q)$}.
\end{align}
Hence, for a subsequence (not indicated), 
\begin{align}\label{eq:thp}
\theta^n \rightarrow \theta \quad \mbox{a.e.~on $Q$}.
\end{align}
This then implies that $ \log \theta^n \rightarrow \log \theta$ a.e.~on $Q$ and $(\theta^n)^{\frac{\beta}{2}} \rightarrow \theta^{\frac{\beta}{2}}$ a.e.~on $Q$. Thanks to the bound appearing in the second line of inequality \eqref{eq:17}, with $\theta$ replaced by $\theta^n$ there, we find that the sequences $(\log \theta^n)_{n=1}^\infty$ and $((\theta^n)^{\frac{\beta}{2}})_{n=1}^\infty$ are equi-integrable in $L^1(Q)$. Thus, by Vitali's convergence theorem, we have the following strong convergence results:
\begin{alignat*}{2}
\log\theta^n &\rightarrow \log \theta\quad && \mbox{strongly in $L^1(Q)$}, \\
(\theta^n)^{\frac{\beta}{2}} &\rightarrow \theta^{\frac{\beta}{2}}\quad && \mbox{strongly in $L^1(Q)$},
\end{alignat*}
as $n \to \infty$. Therefore, by the uniqueness of the weak limit, from \eqref{eq:**} we have that, as $n \to \infty$,
\begin{alignat}{2}
\begin{aligned}
\log\theta^n &\rightharpoonup {\log \theta}\quad && \mbox{weakly in $L^2(0,T;W^{1,2}(\Omega))$}, \label{eq:logfi}\\
(\theta^n)^{\frac{\beta}{2}} &\rightharpoonup {\theta^{\frac{\beta}{2}}}\quad && \mbox{weakly in $L^2(0,T;W^{1,2}(\Omega))$}. 
\end{aligned}
\end{alignat}
In addition, since $T_{k,\varepsilon} \in C^{0,1}(\mathbb{R})$ with Lipschitz constant equal to $1$, it follows from \eqref{eq:str} that, as $n \to \infty$,
\[ \|T_{k,\varepsilon}(\theta^n) - T_{k,\varepsilon}(\theta)\|_{L^1(Q)} \leq \|\theta^n - \theta\|_{L^1(Q)} \rightarrow 0,\]
and therefore, thanks to the uniqueness of the strong limit, $\overline{T_{k,\varepsilon}(\theta)} = T_{k,\varepsilon}(\theta)$. Thus, we have identified $\overline{T_{k,\varepsilon}(\theta)}$; that is, as $n \to \infty$,
\begin{align}\label{eq:*a} 
T_{k,\varepsilon}(\theta^n) \rightarrow T_{k,\varepsilon}(\theta)\quad \mbox{strongly in $L^2(0,T;L^2(\Omega))$}.
\end{align}
Now, \eqref{eq:*a}, \eqref{eq:*b} and the uniqueness of the weak limit yield that
$\overline{\nabx T_{k,\varepsilon}(\theta)}
= \nabx T_{k,\varepsilon}(\theta)$, and therefore
\begin{align*}
\nabx T_{k,\varepsilon}(\theta^n) \rightharpoonup \nabx T_{k,\varepsilon}(\theta) \quad \mbox{weakly in $L^2(0,T;L^2(\Omega;\mathbb{R}^d))$}.
\end{align*}
Further, we have from \eqref{eq:*} that 
\begin{align}\label{eq:*d}  
\partial_t T_{k,\varepsilon}(\theta^n) \rightharpoonup^\ast \partial_t T_{k,\varepsilon}(\theta) \quad \mbox{weakly-$\ast$ in $\mathcal{M}(0,T;(W^{1,z}(\Omega))^\ast)$ for all $z>d$}.
\end{align}
Since the (distributional) time derivative of $T_{k,\varepsilon}(\theta)$ is a (finite signed) Radon measure, which then implies that  $T_{k,\varepsilon}(\theta) \in BV (0,T;(W^{1,z}(\Omega))^\ast)$ for all $z>d$, it makes sense to define, for all $t \in (0,T)$,
\[ T_{k,\varepsilon}(\theta(t_+,\cdot)):= \lim_{s \to t_{+}} T_{k,\varepsilon}(\theta(s,\cdot)) \]
and
\[ T_{k,\varepsilon}(\theta(t_-,\cdot)):= \lim_{s \to t_{-}} T_{k,\varepsilon}(\theta(s,\cdot)), \]
where both limits are considered in the space $(W^{1,z}(\Omega))^\ast$, $z>d$.\footnote{See, for example, Proposition 4.2 in \cite{MR957087}, which is also available from
\href{https://hal.science/hal-01363799v1}{https://hal.science/hal-01363799v1}, concerning the existence of one-sided limits of functions contained in $BV(0,T;E)$, where $E$ is a complete metric space.}

In addition, we also have, for all $k\geq 1$ and $\varepsilon \in (0,1)$ that, as $n \to \infty$, 
\begin{align}\label{eq:*r}  
\kappa(\theta^n) \nabx T_{k,\varepsilon}(\theta^n) = \kappa(\theta^n)
T_{k,\varepsilon}'(\theta^n) \nabx \theta^n \rightharpoonup \kappa(\theta) T_{k,\varepsilon}'(\theta) \nabx\theta = \kappa(\theta) \nabx T_{k,\varepsilon}(\theta) \quad \mbox{weakly in $L^2(Q;\mathbb{R}^d)$}.
\end{align}

Using \eqref{eq:*d}, \eqref{eq:va} (with $r=2$), \eqref{eq:*a} and \eqref{eq:*r} we can now pass to the limit on the left-hand side of the weak formulation of \eqref{eq:ren} as $n \to \infty$. Since the right-hand side of \eqref{eq:ren} also involves $\varphi^n$, we need to establish suitable convergence results for $\varphi^n$ to be able to pass to the limit there. 

Therefore, in the next subsection we explore the convergence properties of the sequence of probability density functions $(\varphi^n)_{n=1}^\infty$.

\subsection{Limit passage in the sequence of probability density functions}

Thanks to \eqref{eq:19}, we have that
\[
\int_{Q \times D}\frac{1}{\theta^n \varphi^n} |\bj^n_\varphi|^2 \dd \bq \dd \bx \dd t \leq C.
\]
Therefore, as $n \to \infty$,
\begin{align}\label{eq:jbd}
\frac{1}{\sqrt{\theta^n \varphi^n}}\bj^n_\varphi \rightharpoonup \overline{\frac{1}{\sqrt{\theta \varphi}}\bj_{\varphi,\bq}}\quad \mbox{weakly in $L^2(Q \times D;\mathbb{R}^d)$},
\end{align}
where $\overline{\frac{1}{\sqrt{\theta \varphi}}\bj_{\varphi,\bq}}$ is a weak limit, to be identified. 

Further, from \eqref{eq:18}, we have
\begin{align}\label{eq:hat} 
\int_{Q \times D}
\theta^n \, \mathrm{e}^{-\frac{\cU_e + \theta^n\, \cU_\eta}{\theta^n}} |\nabq \sqrt{\varphi^{n}  \, \mathrm{e}^{\frac{\cU_e + \theta^n\, \cU_\eta}{\theta^n}}  }|^2\dd \bq \dd \bx \dd t = \frac{1}{4}\int_{Q \times D}
\theta^n \, \mathrm{e}^{-\frac{\cU_e + \theta^n\, \cU_\eta}{\theta^n}} \frac{|\nabq \widehat{\varphi^{n}}|^2}{\widehat{\varphi^n}} \dd \bq \dd \bx \dd t \leq C,
\end{align}
and therefore, as $n \to \infty$,
\begin{align}\label{eq:wL2} 
\sqrt{\theta^n \, \mathrm{e}^{-\frac{\cU_e + \theta^n\, \cU_\eta}{\theta^n}} } \nabq \sqrt{\varphi^{n}  \, \mathrm{e}^{\frac{\cU_e + \theta^n\, \cU_\eta}{\theta^n}}  } \rightharpoonup \overline{  \sqrt{\theta\, \mathrm{e}^{-\frac{\cU_e + \theta\, \cU_\eta}{\theta}} } \nabq \sqrt{\varphi  \, \mathrm{e}^{\frac{\cU_e + \theta\, \cU_\eta}{\theta}}}} \quad  \mbox{weakly in $L^2(Q \times D;\mathbb{R}^d)$},
\end{align}
where $\overline{  \sqrt{\theta\, \mathrm{e}^{-\frac{\cU_e + \theta\, \cU_\eta}{\theta}} } \nabq \sqrt{\varphi  \, \mathrm{e}^{\frac{\cU_e + \theta\, \cU_\eta}{\theta}}}}$ is a weak limit, to be identified. 

It follows from \eqref{eq:26} that, for all $n \geq 1$,
\begin{align}\label{eq:29a} 
\|\varphi^n \log \varphi^n\|_{L^\infty(0,T;L^1(\Omega \times D))} \leq C.
\end{align}
Thus, by de la Vallée Poussin's criterion,
\begin{align}\label{eq:29b} 
\varphi^n \rightharpoonup \varphi\quad \mbox{weakly in $L^p(0,T;L^1(\Omega \times D))$}\quad \mbox{for all $p \in [1,\infty)$},
\end{align}
as $n\to \infty$. It further follows from the last bound in \eqref{eq:28} that, for all $n \geq 1$,  
\[ \|\sqrt{\theta^n} \nabx \sqrt{\varphi^n}\|_{L^2(Q \times D)} \leq C.\]
Therefore, as $n \to \infty$,
\[ \sqrt{\theta^n} \nabx \sqrt{\varphi^n} \rightharpoonup \overline{\sqrt{\theta} \nabx \sqrt{\varphi}}\quad \mbox{weakly in $L^2(Q \times D;\mathbb{R}^d)$},\]
where $\overline{\sqrt{\theta} \nabx \sqrt{\varphi}}$ is a weak limit, to be identified. 

We recall that we have already shown above that  $\theta^n \to \theta$ a.e.~on $Q$ as $n \to \infty$. In order to identify the various unidentified weak limits listed above that involve the sequences $(\theta^n)_{n=1}^\infty$ and $(\varphi^n)_{n=1}^\infty$, it therefore remains to be shown that, as $n \to \infty$,
\[ \varphi^n \to \varphi \quad \mbox{a.e.~on $Q \times D$}.\]
We shall do so via a compensated compactness argument based on the div-curl lemma. To this end, however, we need to establish additional estimates on expressions entering in the Fokker--Planck equation satisfied by $\varphi^n$. 
For ease of writing, suppose that $d=3$ (for $d=2$ the argument is identical). We begin by noting that the Fokker--Planck equation satisfied by $\varphi^n$ is, in fact, the divergence with respect to the variables $(t,\bx,\bq) \in (0,T) \times \Omega \times D \subset \mathbb{R}^7$ of the vector function
\[ \boldsymbol{H}^n:=
(\varphi^n, \bv^n \varphi^n - \theta^n \nabx \varphi^n , (\nabx \bv^n) \bq \varphi^n + \bj_{\varphi,\bq}^n).\]
In terms of this notation, the Fokker--Planck equation satisfied by $\varphi^n$ is simply
\[ \mathrm{div}_{t,\bx,\bq}(\boldsymbol{H}^n) = 0\quad \mbox{for all $n \geq 1$}.\]
In a sequence of steps, we shall next bound $\varphi^n$ in the norm of $L^r(\Omega \times \wtD)$ for some $r \in (1,2)$, and each of $\bv^n\varphi^n$, $\theta^n\nabx\varphi^n$, $(\nabx \bv^n) \bq \varphi^n$ and $\bj_{\varphi,\bq}^n$ in the norm of $L^r(Q \times \wtD;\mathbb{R}^3)$ for some $r \in (1,2)$, independently of $n$, which will then imply that the sequence $\boldsymbol{H}^n$ is bounded in $L^r(Q \times \wtD;\mathbb{R}^7)$ for some $r \in (1,2)$.

According to \eqref{eq:29}, for any open subset $\wtD$ of $D$ such that $\wtD\subset\joinrel\subset D$, 
\begin{align}\label{eq:29c} \|\varphi^n\|_{L^{1+\delta} (Q \times \wtD)} \leq C_{\wtD}
\end{align}
for some $\delta>0$ (independent of $\wtD$) and all $n \geq 1$. Hence, 
\begin{align}\label{eq:30}
\varphi^n \rightharpoonup \varphi \quad \mbox{weakly in $L^{1+\delta}(Q \times \wtD)$}
\end{align}
as $n \to \infty$.

\bigskip

It follows from \eqref{eq:28} and the fact that $\theta$ is bounded below by a positive constant on $Q$ (cf. \eqref{eq:20}) that 
\[ \int_Q \|\nabq \sqrt{\varphi^n}\|^2_{L^2(\wtD)} \leq C_{\wtD}. \]

Because $\|\varphi^n\|_{L^\infty(Q;L^1(D))} \leq C$, also $\|\sqrt{\varphi^n}\|_{L^\infty(Q;L^2(D))} \leq C$, by which $\|\sqrt{\varphi^n}\|_{L^2(Q;W^{1,2}(\wtD))} \leq C_{\wtD}$. Hence, by Sobolev embedding, $\|\varphi^n\|_{L^1(Q;L^3(\wtD))} = \|\sqrt{\varphi^n}\|_{L^2(Q;L^6(\wtD))}^2 \leq C_{\wtD}$. 

Next, let $r \in (1,\frac{25p}{15p+12}]$ with $p=2 - \frac{6}{5+3\beta} = \frac{2(2+3\beta)}{5+3\beta}$.  We note that $\frac{25p}{15p+12} =\frac{25(2+3\beta)}{60+63\beta}>1$ for all $\beta>\frac{5}{6}$, and therefore the interval $(1,\frac{25p}{15p+12}]$ is not empty. For such $r$, we then define $\gamma:=\frac{5p}{6r}$ and $\lambda:=\frac{3-r}{2r}$. Hence, $r \gamma =\frac{5p}{6} = \frac{5}{3} \frac{2+3\beta}{5+3\beta}>1$ and 
\[0< (1-\lambda)r \frac{\gamma}{\gamma-1} \leq 1.\]
Thus, by function space interpolation and H\"older's inequality, 
\begin{align*}
\int_{Q \times \wtD} |\bv^n \varphi^n|^r \dd \bq \dd \bx \dd t &= \int_Q |\bv^n|^r \left(\int_{\wtD}|\varphi^n|^r \dd \bq \right)  \dd \bx \dd t = \int_Q |\bv^n|^r \|\varphi^n\|_{L^r(\wtD)}^r \dd \bx \dd t\\
&\leq \int_Q |\bv^n|^r \|\varphi^n\|_{L^1(\wtD)}^{\lambda r}\|\varphi^n\|_{L^3(\wtD)}^{(1-\lambda)r} \dd \bx \dd t\\
&\leq \|\varphi^n\|_{L^\infty(Q;L^1(D))}^{\lambda r} \left(\int_Q |\bv^n|^{r \gamma} \dd \bx \dd t \right)^{\frac{1}{\gamma}} \left(\int_Q \|\varphi^n\|_{L^3(\wtD)}^{(1-\lambda) r \frac{\gamma}{\gamma-1}} \dd \bx \dd t\right)^{\frac{\gamma-1}{\gamma}}.
\end{align*}
Because $0<(1-\lambda) r \frac{\gamma}{\gamma-1}\leq 1$ and the sequence $(\varphi^n)_{n=1}^\infty$ is bounded in $L^1(Q;L^3(\wtD))$, it then follows that
\begin{align}\label{eq:30b}
\int_{Q \times \wtD} |\bv^n \varphi^n|^r \dd \bq \dd \bx \dd t &\leq  \|\varphi^n\|_{L^\infty(Q;L^1(D))}^{\lambda r}\|\bv^n\|_{L^{\frac{5p}{6}}(Q)}^r \left(\int_Q \|\varphi^n\|_{L^3(\wtD)}^{(1-\lambda) r \frac{\gamma}{\gamma-1}} \dd \bx \dd t\right)^{\frac{\gamma-1}{\gamma}} \leq C_{\wtD}.
\end{align}

\bigskip

Next, let $r \in \left(1, \frac{2+2 \delta}{2+\delta}\right)$; thus $r<2$ and $1< \frac{r}{2-r}<1+\delta$. Then, it follows by \eqref{eq:28} that 
\begin{align*}
\int_{Q \times \wtD} |\theta^n \nabxq \varphi^n|^r \dd \bq \dd \bx \dd t &= 2^r \int_{Q \times \wtD} |\sqrt{\theta^n} \nabxq \sqrt{\varphi^n}|^r |\sqrt{\theta^n \varphi^n}|^r \dd \bq \dd \bx \dd t \\
& \leq 2^r \left(\int_{Q \times \wtD} \theta^n |\nabxq \sqrt{\varphi^n}|^2 \dd \bq \dd \bx \dd t  \right)^{\frac{r}{2}}
\left(\int_{Q \times \wtD}|\theta^n \varphi^n|^{\frac{r}{2-r}} \dd \bq \dd \bx \dd t\right)^{\frac{2-r}{r}} \\
& \leq C\left(\int_Q |\theta^n|^{\frac{r}{2-r}} \left(\int_{\wtD} |\varphi^n|^{\frac{r}{2-r}} \dd \bq\right) \dd \bx \dd t\right)^{\frac{2-r}{2}}.
\end{align*}
By interpolation, 
\[ \|\varphi^n\|^{\frac{r}{2-r}}_{L^{\frac{r}{2-r}}(\wtD)} \leq \|\varphi^n\|^{\frac{(1-\lambda)r}{2-r}}_{L^1(\wtD)} \|\varphi^n\|^{\frac{\lambda r}{2-r}}_{L^{1+\delta}(\wtD)},\]
where
\[ \frac{2-r}{r} = \frac{1-\lambda}{1} + \frac{\lambda}{1+\delta} \quad \Leftrightarrow \quad \frac{2(1-r)}{r} = -\lambda \frac{\delta}{1+ \delta} \quad \Leftrightarrow \quad \lambda = \frac{2(1+\delta) (r-1)}{\delta r} .\]
Clearly, $\lambda \in (0,1)$. Hence, 
\begin{align*}
\int_{Q \times \wtD} |\theta^n \nabxq \varphi^n|^r \dd \bq \dd \bx \dd t &\leq C\left(\int_Q |\theta^n|^{\frac{r}{2-r}} \|\varphi^n\|^{\frac{\lambda r}{2-r}}_{L^{1+\delta}(\wtD)}  \dd \bx \dd t\right)^{\frac{2-r}{2}}\\
& = C \left(\int_Q |\theta^n|^{\frac{r}{2-r}} \|\varphi^n\|^{\frac{2(1+\delta)(r-1)}{\delta(2-r)}}_{L^{1+\delta}(\wtD)}  \dd \bx \dd t\right)^{\frac{2-r}{2}}.
\end{align*}
Now, let $\gamma=\frac{\delta(2-r)}{2(r-1)}$ and let $\gamma'= \frac{\gamma}{\gamma-1}$; evidently, $\gamma, \gamma' >1$. Hence, by H\"older's inequality, 
\begin{align*}
 \int_{Q \times \wtD} |\theta^n \nabxq \varphi^n|^r \dd \bq \dd \bx \dd t &\leq C  \left(\int_Q |\theta^n| 
 ^{\frac{r\gamma'}{2-r}} \dd \bx \dd t\right)^{\frac{2-r}{2\gamma'}}
 \left(\int_Q
 \|\varphi^n\|^{\frac{2(1+\delta)(r-1)}{\delta(2-r)}\gamma}_{L^{1+\delta}(\wtD)}  \dd \bx \dd t\right)^{\frac{2-r}{2\gamma}} \\
 &= C  \left(\int_Q |\theta^n| 
 ^{\frac{r\gamma'}{2-r}} \dd \bx \dd t\right)^{\frac{2-r}{2\gamma'}}
 \left(\int_Q
 \|\varphi^n\|^{1+\delta}_{L^{1+\delta}(\wtD)}  \dd \bx \dd t\right)^{\frac{2-r}{2\gamma}} \\
 &= C  \left(\int_Q |\theta^n| 
 ^{\frac{r\gamma'}{2-r}} \dd \bx \dd t\right)^{\frac{2-r}{2\gamma'}}
 \|\varphi^n\|^{\frac{(1+\delta)(2-r)}{2\gamma}}_{L^{1+\delta}(Q \times \wtD)}\\
 &\leq C_{\wtD} \left(\int_Q |\theta^n| 
 ^{\frac{r\gamma'}{2-r}} \dd \bx \dd t\right)^{\frac{2-r}{2\gamma'}},
\end{align*}
where the last inequality follows from \eqref{eq:29}. According to \eqref{eq:22} the right-hand side of the last inequality will be bounded by a constant, independent of $n$, provided that 
\begin{align*} 
\frac{r\gamma'}{2-r} \leq \frac{2}{3} + \beta \quad \Leftrightarrow \quad \frac{r}{2-r} \,\frac{\delta (2-r)}{\delta(2-r) - 2(r-1)} \leq \frac{2}{3} + \beta\quad &\Leftrightarrow \quad \frac{\delta r}{\delta(2-r) - 2(r-1)} \leq \frac{2}{3} + \beta\\
\Leftrightarrow \quad 3\delta r\leq (2+3\beta)(2\delta + 2) - (2+3\beta)(\delta + 2)r\quad &\Leftrightarrow \quad 
r \leq \frac{2(2+3\beta)(\delta + 1)}{3\delta + (2+3\beta)(\delta + 2)};
\end{align*}
that is, if
\[ 1 < r \leq \frac{2(2+3\beta)(\delta + 1)}{3\delta + (2+3\beta)(\delta + 2)}.\]
We note that the fraction on the right-hand side of this inequality is $>1$ thanks to our assumption that $\beta>\frac{5}{6} ~(> \frac{1}{3})$ and is $<\frac{2+2\delta}{2+\delta}<2$, so if this inequality is satisfied, then automatically $r \in \left(1,\frac{2+2\delta}{2+\delta}\right)$. In particular, we can set 
\[ r = 
\frac{2(2+3\beta)(\delta + 1)}{3\delta + (2+3\beta)(\delta + 2)}.\]
With this choice of $r$ we then have that
\begin{align}\label{eq:30a}
 \int_{Q \times \wtD} |\theta^n \nabxq \varphi^n|^r \dd \bq \dd \bx \dd t &\leq C_{\wtD}.
 \end{align}
It remains to derive a bound on the final component of $\boldsymbol{H}^n$, which is $(\nabx \bv^n) \bq \varphi^n + \bj_{\varphi,\bq}^n$.  For any $r \in (1,p)$, where $p=2 - \frac{6}{5+3\beta}$, by H\"older's inequality and interpolation with $\lambda:=\frac{3-r}{2r} \in (0,1)$, recalling that $|\bq|^2 \leq b$ for all $\bq \in D$, we have
\begin{align}\label{eq:30c}
\begin{aligned}
\int_{Q \times \wtD}|(\nabx \bv^n)\bq \varphi^n|^r \dd \bq \dd \bx \dd t   &\leq b^{\frac{r}{2}} \int_Q |\nabx \bv^n|^r \|\varphi^n\|_{L^r(\wtD)}^r \dd \bx \dd t \\
&\leq b^{\frac{r}{2}} \left(\int_Q |\nabx \bv^n|^p \dd \bx \dd t\right)^{\frac{r}{p}} \left(\int_Q \|\varphi^n\|^{\frac{pr}{p-r}}_{L^r(\wtD)} \dd \bx \dd t\right)^{\frac{p-r}{p}}\\
&\leq b^{\frac{r}{2}}  \left(\int_Q |\nabx \bv^n|^p \dd \bx \dd t\right)^{\frac{r}{p}} \left(\int_Q \|\varphi^n\|^{\frac{\lambda pr}{p-r}}_{L^1(\wtD)} \|\varphi^n\|^{\frac{(1-\lambda) pr}{p-r}}_{L^3(\wtD)} \dd \bx \dd t\right)^{\frac{p-r}{p}}\\
&\leq b^{\frac{r}{2}} \|\nabx \bv^n\|_{L^p(Q)}^r \|\varphi^n\|_{L^\infty(Q;L^1(D))}^{\lambda r}
\left(\int_D \|\varphi^n\|_{L^3(\wtD)}^{\frac{3p(r-1)}{2(p-r)}}\dd \bx \dd t\right)^{\frac{p-r}{p}} \leq C_{\wtD}
\end{aligned}
\end{align}
for all $n \geq 1$, provided that $r \in (1,p)$ is such that 
\[ 0 < \frac{3p(r-1)}{2(p-r)} \leq 1 \quad \Leftrightarrow \quad 1 < r \leq \frac{5p}{3p+2} (<p) .\]
For example, we can take (again, because by assumption $\beta > \frac{5}{6} ~(>\frac{1}{3})$)
\[ r = \frac{5p}{3p+2} = \frac{5(2+3\beta)}{11+12\beta} \in (1,2).\]
We move on to bound $\bj_{\varphi,\bq}^n := - \boldsymbol{F}^n\varphi^n - \theta^n \nabq \varphi^n$ and recall that $\boldsymbol{F}^n= \nabq\, \cU_e + \theta^n \nabq \,\cU_\eta$; 
hence we need to bound $\bj_{\varphi,\bq}^n := - (\nabq\, \cU_e + \theta^n \nabq \,\cU_\eta)\varphi^n - \theta^n \nabq \varphi^n$. Thanks to \eqref{eq:30a} we already have a bound on the final summand in $\bj_{\varphi,\bq}^n$; therefore, it remains to deal with the first two summands of $\bj_{\varphi,\bq}^n$. Since $\nabq\,\cU_e$ and $\nabq\,\cU_\eta$ are uniformly continuous and therefore bounded functions on $\overline{\wtD}$,  we only need to derive bounds on $\varphi^n$ and $\theta^n \varphi^n$ in the norm of $L^r(Q \times \widetilde{D})$ for some $r \in (1,2)$.

For any $r \in (1,\frac{5}{3}]$ and $\lambda:= \frac{3-r}{2r} \in (0,1)$, we have
\begin{align*}
   \int_{Q \times \wtD} |\varphi^n|^r \dd \bq \dd\bx \dd t &= \int_Q \|\varphi^n\|^r_{L^r(\wtD)} \dd \bx \dd t \leq \int_Q \|\varphi^n\|_{L^1(\wtD)}^{\lambda r} \|\varphi^n\|_{L^3(\wtD)}^{(1-\lambda)r} \dd \bx \dd t\\
   &\leq \|\varphi^n\|_{L^\infty(Q;L^1(D))}^{\lambda r}\int_D \|\varphi^n\|_{L^3(\wtD)}^{(1-\lambda)r}\dd \bx \dd t \leq C_{\wtD},
\end{align*}
because
\[ (1-\lambda) r = \frac{3(r-1)}{2}  \leq 1 \quad \mbox{for $r \in (1,\frac{5}{3}]$} .\]

Similarly, for any $r \in (1,\frac{3}{2}]$ and $\lambda:=\frac{3-r}{2r}$, we have $r < \beta + \frac{2}{3}$ for $\beta>\frac{5}{6}$; therefore, $\gamma:=\frac{2+3\beta }{3r}>1$. Since $\frac{2}{3r-1}<1$, it follows that 
\[ 0< (1-\lambda)r \frac{\gamma}{\gamma-1} < 1.\]
Hence, 
\begin{align*}
   \int_{Q \times \wtD} |\theta^n \varphi^n|^r \dd \bq \dd \bx \dd t &= \int_Q |\theta^n|^r \|\varphi\|^r_{L^r(\wtD)}   \dd \bx \dd t 
   \leq \int_Q |\theta^n|^r \|\varphi^n\|_{L^1(\wtD)}^{\lambda r} \|\varphi^n\|_{L^3(\wtD)}^{(1-\lambda)r}  \dd \bx \dd t \\
   &\leq \|\varphi^n\|^{\lambda r}_{L^\infty(Q;L^1(\wtD))} \left(\int_Q |\theta^n|^{\gamma r}  \dd \bx \dd t\right)^{\frac{1}{\gamma}}  
   \left( \int_Q \|\varphi^n\|_{L^3(\wtD)}^{(1-\lambda)r \frac{\gamma}{\gamma -1}}  \dd \bx \dd t \right)^{\frac{\gamma-1}{\gamma}}\\
   & = \|\varphi^n\|^{\lambda r}_{L^\infty(Q;L^1(\wtD))} 
   \|\theta^n\|^r_{L^{\frac{2}{3}+\beta}(Q)} 
   \left( \int_Q \|\varphi^n\|_{L^3(\wtD)}^{(1-\lambda)r \frac{\gamma}{\gamma -1}}  \dd \bx \dd t \right)^{\frac{\gamma-1}{\gamma}} \leq C_{\wtD}.
\end{align*}
These bounds on $(\nabx\bv^n)\bq\varphi^n$, $\varphi^n$ and $\theta^n\varphi^n$, together with \eqref{eq:30a}, imply that, for some $r \in (1,2)$, independent of $\wtD$ and $n$,
\begin{align}\label{eq:30d}  \int_{Q \times \wtD} |(\nabx \bv^n)\bq \varphi^n + \bj_{\varphi,\bq}^n|^r \dd \bq \dd \bx \dd t \leq C_{\wtD}\quad \mbox{for all $n \geq 1$.}
\end{align}

Combining the bounds \eqref{eq:29c}, \eqref{eq:30b}, \eqref{eq:30a}, \eqref{eq:30c} and \eqref{eq:30d} it then follows that, for some $\delta \in (0,1)$, 
\[\|\boldsymbol{H}^n\|_{L^{1+\delta}((0,T) \times \Omega \times \wtD)} \leq C_{\wtD}\quad \mbox{for all $n \geq 1$}.\]
Recall further that $\mbox{div}_{t,\bx,\bq}(\boldsymbol{H}^n)=0$. 

\medskip

In preparation for the application of the div-curl lemma, we consider another 7-component vector-function: 
\[ \boldsymbol{Q}^n:=(T_{k,\varepsilon}(\varphi^n),\boldsymbol{0}, \boldsymbol{0}),\]
where $\boldsymbol{0} \in \mathbb{R}^3$, 
and we note that $\|\boldsymbol{Q}^n \|_{L^\infty(Q \times D)}\leq k$ for each $k\geq 1$ and every $\varepsilon \in (0,1)$. Next, we bound $\mbox{curl}_{t,\bx,\bq}(\boldsymbol{Q}^n):= \nabla_{t,\bx,\bq}\boldsymbol{Q}^n- (\nabla_{t,\bx,\bq} \boldsymbol{Q}^n)^{\mathrm{T}}$. Clearly, 
\begin{align*}
\| \mbox{curl}_{t,\bx,\bq}(\boldsymbol{Q}^n)\|_{L^2(Q \times \wtD)} &\leq   2 \|\nabla_{\bx,\bq} T_{k,\varepsilon}(\varphi^n)\|_{L^2(Q \times \wtD)} = 2 \| T'_{k,\varepsilon}(\varphi^n)\nabla_{\bx,\bq}\varphi^n \|_{L^2(Q \times \wtD)}\\
& = 4 \| T'_{k,\varepsilon}(\varphi^n) \sqrt{\varphi^n} \nabla_{\bx,\bq}\sqrt{\varphi^n} \|_{L^2(Q \times \wtD)}\\
& \leq C \sqrt{k+\varepsilon} \| \nabla_{\bx,\bq}\sqrt{\varphi^n} \|_{L^2(Q \times \wtD)} \leq 2 \sqrt{k+\varepsilon}\, C_{\wtD},
\end{align*}
where the last inequality follows from \eqref{eq:28a}. 

In summary then, $(\boldsymbol{H}^n)_{n=1}^\infty$ is a bounded sequence in $L^{1+\delta}((0,T) \times \Omega \times \wtD;\mathbb{R}^7)$ for some $\delta \in (0,1)$ and $\mbox{div}_{t,\bx,\bq}(\boldsymbol{H}^n)=0$, whereby $\mbox{div}_{t,\bx,\bq}(\boldsymbol{H}^n)$ is, trivially, precompact in $W^{-1,2}((0,T) \times \Omega \times \wtD):= [W^{1,2}_0((0,T) \times \Omega \times \wtD)]^\ast$. Furthermore, $(\boldsymbol{Q}^n)_{n=1}^\infty$ is a bounded sequence in $L^\infty((0,T) \times \Omega \times \wtD;\mathbb{R}^7)$  and $\mbox{curl}_{t,\bx,\bq}(\boldsymbol{H}^n)$, is precompact in $W^{-1,2}((0,T) \times \Omega \times \wtD;\mathbb{R}^{7 \times 7}):= [W^{1,2}_0((0,T) \times \Omega \times \wtD;\mathbb{R}^{7 \times 7})]^\ast$ thanks to the compact embedding of $L^2((0,T) \times \Omega \times \wtD;\mathbb{R}^{7 \times 7})$ into $W^{-1,2}((0,T) \times \Omega \times \wtD;\mathbb{R}^{7 \times 7})$.

Thus, by the div-curl lemma, for each $k \geq 1$ and $\varepsilon \in (0,1)$, as $n \rightarrow \infty$, 
\begin{align}\label{eq:31}
 \boldsymbol{H}^n \cdot \boldsymbol{Q}^n = \varphi^n T_{k,\varepsilon}(\varphi^n) \rightharpoonup \varphi \,\overline{T_{k,\varepsilon}(\varphi)}\quad \mbox{weakly in $L^{1+\delta'}((0,T) \times \Omega \times \wtD;\mathbb{R}^7)$},
\end{align}
where $0<\delta'<\delta$ and $\overline{T_{k,\varepsilon}(\varphi)} \in L^\infty((0,T) \times \Omega \times \wtD)= L^\infty(Q \times \wtD)$ is the weak-$\ast$ limit of the bounded sequence $(T_{k,\varepsilon}(\varphi^n))_{n=1}^\infty \subset L^\infty(Q \times \wtD)$, to be identified. In order to identify $\overline{T_{k,\varepsilon}(\varphi)}$ we note that, thanks to the monotonicity of the function $T_{k,\varepsilon}$, one has that 
\[ (T_{k,\varepsilon}(a) 
- T_{k,\varepsilon}(b))^2 
\leq (T_{k,\varepsilon}(a) - T_{k,\varepsilon}(b))
(a-b) \quad \forall\, a,b \in \mathbb{R}.
\]
Hence, 
\begin{align*}
\int_{Q \times \wtD}(T_{k,\varepsilon}(\varphi^n) - T_{k,\varepsilon}(\varphi^m))^2 \dd \bq \dd \bx \dd t \leq \int_{Q \times \wtD} (T_{k,\varepsilon}(\varphi^n) - T_{k,\varepsilon}(\varphi^m))(\varphi^n - \varphi^m)\dd \bq \dd \bx \dd t.    
\end{align*}
Passing to the limit $m \rightarrow \infty$, with  $n\geq 1$, $k\geq 1$ and $\varepsilon \in (0,1)$ fixed, using the weak lower-semicontinuity of the $L^2(Q \times \wtD)$ norm, it then follows that 
\begin{align*}
\int_{Q \times \wtD}(T_{k,\varepsilon}(\varphi^n) - \overline{T_{k,\varepsilon}(\varphi)})^2 \dd \bq \dd \bx \dd t &\leq \mathrm{lim.inf}_{m \to \infty}\int_{Q \times \wtD} (T_{k,\varepsilon}(\varphi^n) - T_{k,\varepsilon}(\varphi^m))(\varphi^n - \varphi^m)\dd \bq \dd \bx \dd t\\
& =\int_{Q \times \wtD} (T_{k,\varepsilon}(\varphi^n) - \overline{T_{k,\varepsilon}(\varphi)})(\varphi^n - \varphi)\dd \bq \dd \bx \dd t\\
& =\int_{Q \times \wtD} (T_{k,\varepsilon}(\varphi^n) \varphi^n
- \overline{T_{k,\varepsilon}(\varphi)}(\varphi^n - \varphi)
- T_{k,\varepsilon}(\varphi^n)\varphi)\dd \bq \dd \bx \dd t, 
\end{align*}
where in the transition to the second line we have used \eqref{eq:30} and \eqref{eq:31}. Now we pass to the limit $n \rightarrow \infty$ using \eqref{eq:31} in the first term on the right-hand side, \eqref{eq:30} in the second term on the right-hand side and the definition of $\overline{T_{k,\varepsilon}(\varphi)}$ in the final term on the right-hand side. Hence, 
\begin{align*}
\lim_{n \to \infty}\int_{Q \times \wtD}(T_{k,\varepsilon}(\varphi^n) - \overline{T_{k,\varepsilon}(\varphi)})^2 \dd \bq \dd \bx \dd t &\leq \int_{Q \times \wtD} (\overline{T_{k,\varepsilon}(\varphi)} \varphi
- 0
- \overline{T_{k,\varepsilon}(\varphi)}\varphi)\dd \bq \dd \bx \dd t =0.
\end{align*}
Thus we have shown that, as $n \rightarrow \infty$,  
\begin{align*}
   T_{k,\varepsilon}(\varphi^n) \to \overline{T_{k,\varepsilon}(\varphi)}\quad \mbox{strongly in $L^2(Q \times \wtD)$} 
\end{align*}
for each $k\geq 1$ and $\varepsilon \in (0,1)$, and therefore also 
\begin{align}\label{eq:32}
   T_{k,\varepsilon}(\varphi^n) \to \overline{T_{k,\varepsilon}(\varphi)}\quad \mbox{strongly in $L^1(Q \times \wtD)$}.
\end{align}

Now, by the triangle inequality,
\[\int_{\Omega\times \widetilde{D}} |\varphi^n - \varphi^m| \dd \bq \dd \bx \dd t \leq   \int_{Q \times \wtD} |T_{k,\varepsilon} (\varphi^n) - {T_{k,\varepsilon} (\varphi^m)}| + |\varphi^n - T_{k,\varepsilon} (\varphi^n)| + |\varphi^m - T_{k,\varepsilon} (\varphi^m)|\dd \bq \dd \bx \dd t. \]
Hence, using \eqref{eq:29b} and weak lower-semicontinuity in the $L^1(Q \times D)$ norm in conjunction with \eqref{eq:32}, it follows by passing to the $\mathrm{lim.inf}_{m \to \infty}$ that
\begin{align*}
   \int_{Q \times \wtD} |\varphi^n -  \varphi| \dd \bq \dd \bx \dd t &\leq
\int_{Q \times \wtD} |T_{k,\varepsilon} (\varphi^n) - \overline{T_{k,\varepsilon} (\varphi)}|\dd \bq \dd \bx \dd t  + \int_{Q \times \wtD}|\varphi^n|\,1\!\!1_{\{|\varphi^n| > k\}} \dd \bq \dd \bx \dd t \\
&\quad + \mathrm{lim.inf}_{m\to \infty}\int_{Q \times \wtD} |\varphi^m - T_{k,\varepsilon} (\varphi^m)|\dd \bq \dd \bx \dd t\\
&\leq
\int_{Q \times \wtD} |T_{k,\varepsilon} (\varphi^n) - \overline{T_{k,\varepsilon} (\varphi)}|\dd \bq \dd \bx \dd t  + \int_{Q \times \wtD}|\varphi^n \log \varphi^n|\frac{\,1\!\!1_{\{|\varphi^n| > k\}}}{|\log \varphi^n|} \dd \bq \dd \bx \dd t \\
&\quad + \mathrm{lim.inf}_{m\to \infty}\int_{Q \times \wtD} |\varphi^m - T_{k,\varepsilon} (\varphi^m)|\dd \bq \dd \bx \dd t\\
 &\leq
\int_{Q \times \wtD} |T_{k,\varepsilon} (\varphi^n) - \overline{T_{k,\varepsilon} (\varphi)}|\dd \bq \dd \bx \dd t  + \frac{C}{\log k} \\
&\quad + \mathrm{lim.inf}_{m\to \infty}\int_{Q \times \wtD} |\varphi^m - T_{k,\varepsilon} (\varphi^m)|\dd \bq \dd \bx \dd t\\
&\leq \int_{Q \times \wtD} |T_{k,\varepsilon} (\varphi^n) - \overline{T_{k,\varepsilon} (\varphi)}|\dd \bq \dd \bx \dd t  + \frac{C}{\log k} \\
&\quad + \mathrm{lim.sup}_{m\to \infty} \int_{Q \times \wtD}|\varphi^m \log \varphi^m|\frac{\,1\!\!1_{\{|\varphi^m| > k\}}}{|\log \varphi^m|} \dd \bq \dd \bx \dd t \\
&\leq \int_{Q \times \wtD} |T_{k,\varepsilon} (\varphi^n) - \overline{T_{k,\varepsilon} (\varphi)}|\dd \bq \dd \bx \dd t  + \frac{2C}{\log k},
\end{align*}
for all $k>1$ and all $\varepsilon \in (0,1)$, where we have used \eqref{eq:29a} in transitioning from the second inequality to the third inequality and from the fourth inequality to the fifth inequality. Hence, by first passing to the limit $n \to \infty$ using \eqref{eq:32} with $k>1$ held fixed, and then by passing to the limit $k \to \infty$, we deduce that 
\begin{align}\label{eq:32a}
\varphi^n \to \varphi \quad \mbox{strongly in $L^1(Q \times \wtD)$}.
\end{align}
This then implies that, for a subsequence, not indicated, 
\begin{align*}
\varphi^n \to \varphi \quad \mbox{a.e.~in $Q \times \wtD$}.
\end{align*}
By selecting a sequence of nested sets $(\wtD_\ell)_{\ell=1}^\infty$ such that $\wtD_\ell\subset\joinrel\subset D$ and $\bigcup_{\ell=1}^\infty \wtD_\ell = D$, and using that $\varphi^n \to \varphi$ a.e.~in $\Omega \times \wtD_\ell$ for each $\ell \geq 1$, the application of a diagonal procedure then produces a subsequence (that is, once again, not relabelled) such that, as $n \to \infty$, 
\begin{align}\label{eq:33}
\varphi^n \to \varphi \quad \mbox{a.e.~in $Q \times D$}.
\end{align}
Thanks to Vitali's convergence theorem, \eqref{eq:29b} and \eqref{eq:33} finally imply that 
\begin{align}\label{eq:33a}
\varphi^n \to \varphi \quad \mbox{strongly in $L^1(Q \times D)$}.
\end{align}

\subsection{Limit passage in the sequence of stress tensors}

The pointwise convergence results \eqref{eq:thp} and \eqref{eq:33} now enable us to identify the still unidentified weak limit $\overline\bSS$ of the sequence $(\bSS^n)_{n=1}^\infty \subset L^p(Q;\mathbb{R}^{d \times d})$, with $p=2-\frac{6}{5+3\beta}$ and $\beta>\frac{5}{6}$, which appeared above:
\[ \bSS^n = 2\nu(\theta^n) \bDD(\bv^n) + \int_D (\boldsymbol{F}^n \otimes \bq) \varphi^n  \dd \bq,  \]

We begin by considering the first summand in $\bSS^n$. As $\nu \in BC([0,\infty);\mathbb{R}_{\geq 0})$ and $\theta^n \to \theta$ a.e.~on $Q$, the sequence of measurable functions $(\nu(\theta^n))_{n=1}^\infty$ converges to $\nu(\theta)$ a.e.~on $Q$ and $\|\nu(\theta^n)\|_{L^\infty(Q)} \leq C$ for all $n \geq 1$. Furthermore, because of \eqref{eq:v}, we have that, as $n \to \infty$,
\[ \bDD(\bv^n)\rightharpoonup \bDD(\bv)\quad \mbox{weakly in $L^{\frac{4+6\beta}{5+3\beta}}(Q;\mathbb{R}^{d\times d})$},\quad \beta>\frac{5}{6},
\]
and therefore also $\bDD(\bv^n)\rightharpoonup \bDD(\bv)$ weakly in $L^1(Q;\mathbb{R}^{d \times d})$, as $n\to \infty$. Hence, by Proposition 2.61 on p.183 in \cite{Fonseca}, $\nu(\theta^n) \bDD(\bv^n) \rightharpoonup \nu(\theta) \bDD(\bv)$ in $L^1(Q;\mathbb{R}^{d \times d})$ as $n\to \infty$. Since $(\nu(\theta^n) \bDD(\bv^n))_{n=1}^{\infty}$ is a bounded sequence in $L^{\frac{4+6\beta}{5+3\beta}}(Q;\mathbb{R}^{d\times d})$ with $\beta>\frac{5}{6}$, and $\frac{4+6\beta}{5+3\beta}>1$, there exists a weakly convergent subsequence (not indicated) such that $\nu(\theta^n) \bDD(\bv^n) \rightharpoonup \overline{\nu(\theta) \bDD(\bv)}$ in $L^{\frac{4+6\beta}{5+3\beta}}(Q;\mathbb{R}^{d\times d}) \subset L^1(Q;\mathbb{R}^{d \times d})$ as $n \to \infty$. By the uniqueness of the weak limit, $\overline{\nu(\theta) \bDD(\bv)} = \nu(\theta) \bDD(\bv)$, and therefore as $n \to \infty$,
\begin{align}\label{eq:34}
 \nu(\theta^n) \bDD(\bv^n) \rightharpoonup \nu(\theta) \bDD(\bv)\quad \mbox{weakly in $L^{\frac{4+6\beta}{5+3\beta}}(Q;\mathbb{R}^{d \times d})$ with $\beta>\frac{5}{6}$}.
\end{align}

We continue by considering the weak limit of the second summand of $\bSS^n$. Let us suppose to this end that $\bw \in C([0,T];C^1_0(\overline{\Omega};\mathbb{R}^d))$ and that $\divx\bw=0$ on $Q$. Then 
\begin{align}\label{eq:34a}
&\int_Q \left(\int_D (\bF^n \otimes \bq)\varphi^n \dd \bq \right)  : \nabx \bw \dd \bx \dd t  = \int_Q \left( \int_D ((\nabq\, \cU_e + \theta^n \nabq \,\cU_\eta) \otimes \bq) \varphi^n  \dd \bq\right) : \nabx \bw \dd \bx \dd t
 \nonumber\\
&\qquad \qquad = \int_Q \left(  \int_D \sqrt{\theta^n}\, \sqrt{\varphi^n} \left(\bq \otimes  
 \sqrt{\theta^n}\, \mathrm{e}^{-\frac{\cU_e +\theta^n \,\cU_\eta}{2\theta^n}}\, \frac{\nabq \widehat{\varphi^n}}{\sqrt{\widehat{\varphi^n}}}\right) \dd \bq\right) : \nabx \bw \dd \bx \dd t \nonumber\\
 &\qquad \qquad = 2\int_Q \left(  \int_D \sqrt{\theta^n}\, \sqrt{\varphi^n} \left(\bq \otimes  
 \sqrt{\theta^n \, \mathrm{e}^{-\frac{\cU_e +\theta^n \,\cU_\eta}{\theta^n}}}\, \nabq \sqrt{{\varphi^n \, \mathrm{e}^{\frac{\cU_e +\theta^n \,\cU_\eta}{\theta^n}}}}\right) \dd \bq\right) : \nabx \bw \dd \bx \dd t, 
 \end{align}
where in the transition to the second line we used \eqref{eq:19a} together with the fact that $\divx\bw=0$ in order to eliminate the first term on the right-hand side of \eqref{eq:19a}.

Motivated by the form of the expression in the last line of \eqref{eq:34a}, we begin by showing that
\begin{align}\label{eq:35} 
\sqrt{\theta^n}\, \sqrt{\varphi^n}  \rightharpoonup \sqrt{\theta}\, \sqrt{\varphi}\quad \mbox{strongly in $L^2(Q \times D)$}
\end{align}
and
\begin{align}\label{eq:36}
\begin{aligned}\bq \otimes  
 \sqrt{\theta^n \, \mathrm{e}^{-\frac{\cU_e +\theta^n \,\cU_\eta}{\theta^n}}}\, \nabq \sqrt{{\varphi^n \, \mathrm{e}^{\frac{\cU_e +\theta^n \,\cU_\eta}{\theta^n}}}} \rightharpoonup\bq \otimes  
 \sqrt{\theta \, \mathrm{e}^{-\frac{\cU_e +\theta \,\cU_\eta}{\theta}}}\, \nabq \sqrt{{\varphi \, \mathrm{e}^{\frac{\cU_e +\theta \,\cU_\eta}{\theta}}}}\\\
 \mbox{weakly in $L^2(Q \times D;\mathbb{R}^{d \times d})$}.   
\end{aligned}
\end{align}
as $n \to \infty$.

Because $|\sqrt{a} - \sqrt{b}| \leq \sqrt{|a-b|}$ for all $a,b \in \mathbb{R}_{\geq 0}$, it follows that 
\[ \|\sqrt{\theta^n}\, \sqrt{\varphi^n}  - \sqrt{\theta}\, \sqrt{\varphi}\|^2_{L^2(Q \times D)} \leq  \|\theta^n \varphi^n  - \theta \varphi\|_{L^1(Q \times D)}.\]
We will show that the right-hand side of this inequality tends to 0 as $n \to \infty$. Because $\theta^n \varphi^n \to \theta \varphi$ a.e.~on $Q \times D$, thanks to Vitali's theorem, it suffices to show that the sequence $(\theta^n\varphi^n)_{n=1}^\infty$ is equi-integrable. We begin by noting that, for $\alpha \in (0,1)$, H\"older's inequality implies that 
\begin{align*}
\int_D \varphi^n |\log \varphi^n|^\alpha \dd \bq &= \int_D (\varphi^n)^{1-\alpha} (\varphi^n  |\log \varphi^n|) ^\alpha \dd \bq \leq \|\varphi^n\|_{L^1(D)}^{1-\alpha}\left(\int_D \varphi^n |\log \varphi^n|\dd \bq \right)^\alpha.
\end{align*}
Hence, by Young's inequality, 
\begin{align*}
&\int_{Q \times D}|\theta^n|^{1+\alpha} \varphi^n|\log \varphi^n|^\alpha \dd \bq \dd \bx \dd t \leq  \|\varphi^n\|_{L^\infty(Q;L^1(D))}^{1-\alpha} \int_Q |\theta^n|^{1+\alpha}   \left(\int_D \varphi^n |\log \varphi^n|\dd \bq \right)^\alpha \dd \bx \dd t\\
& \leq \|\varphi^n\|_{L^\infty(Q;L^1(D))}^{1-\alpha} \int_Q \left( (1-\alpha)|\theta^n|^{\frac{1+\alpha}{1-\alpha}} + \alpha \int_D \varphi^n |\log \varphi^n|\dd \bq \right)\dd \bx \dd t\\
&  \leq \|\varphi^n\|_{L^\infty(Q;L^1(D))}^{1-\alpha} \left(\|\theta^n\|_{L^{\frac{1+\alpha}{1-\alpha}}(Q)}^{\frac{1+\alpha}{1-\alpha}}  + \|\varphi^n \log \varphi^n\|_{L^1(Q \times D)} \right).
\end{align*}
We set $\frac{1+\alpha}{1-\alpha} = \frac{2}{3} + \beta$ with $\beta > \frac{5}{6}$, whereby $\alpha= \frac{3\beta -1}{3\beta + 5} \in (0,1)$. Then it follows that, for all $n \geq 1$,
\[ \int_{Q \times D}|\theta^n|^{1+\alpha} \varphi^n|\log \varphi^n|^\alpha \dd \bq \dd \bx \dd t \leq  C,\]
where $C$ is a positive constant independent of $n$, which then implies that the sequence $(\theta^n\varphi^n)_{n=1}^\infty$ is equi-integrable. That completes the proof of \eqref{eq:35}. 

Concerning \eqref{eq:36}, it directly follows from \eqref{eq:wL2} that the weak limit exists, and therefore
\begin{align}\label{eq:37}
\begin{aligned}
\bq \otimes  
 \sqrt{\theta^n \, \mathrm{e}^{-\frac{\cU_e +\theta^n \,\cU_\eta}{\theta^n}}}\, \nabq \sqrt{{\varphi^n \, \mathrm{e}^{\frac{\cU_e +\theta^n \,\cU_\eta}{\theta^n}}}} \rightharpoonup \bq \otimes  
 \overline{\sqrt{\theta \, \mathrm{e}^{-\frac{\cU_e +\theta \,\cU_\eta}{\theta}}}\, \nabq \sqrt{{\varphi \, \mathrm{e}^{\frac{\cU_e +\theta \,\cU_\eta}{\theta}}}}}\\\
 \mbox{weakly in $L^2(Q \times D;\mathbb{R}^{d \times d})$}
 \end{aligned}
 \end{align}
as $n \to \infty$. Therefore, also 
\begin{align*}
\bq \otimes  
 \sqrt{\theta^n \, \mathrm{e}^{-\frac{\cU_e +\theta^n \,\cU_\eta}{\theta^n}}}\, \nabq \sqrt{{\varphi^n \, \mathrm{e}^{\frac{\cU_e +\theta^n \,\cU_\eta}{\theta^n}}}} \rightharpoonup \bq \otimes  
 \overline{\sqrt{\theta \, \mathrm{e}^{-\frac{\cU_e +\theta \,\cU_\eta}{\theta}}}\, \nabq \sqrt{{\varphi \, \mathrm{e}^{\frac{\cU_e +\theta \,\cU_\eta}{\theta}}}}} 
\end{align*}
weakly in $L^2(Q \times \wtD ;\mathbb{R}^{d \times d})$ as $n \to \infty$ for any $\wtD \subset\joinrel\subset D$, and the weak limit is independent of the choice of $\wtD$. We shall identify the weak limit by proving that, as $n \to \infty$,
\begin{align*}
\bq \otimes  
 \sqrt{\theta^n \, \mathrm{e}^{-\frac{\cU_e +\theta^n \,\cU_\eta}{\theta^n}}}\, \nabq \sqrt{{\varphi^n \, \mathrm{e}^{\frac{\cU_e +\theta^n \,\cU_\eta}{\theta^n}}}} \rightharpoonup \bq \otimes  
 \sqrt{\theta \, \mathrm{e}^{-\frac{\cU_e +\theta \,\cU_\eta}{\theta}}}\, \nabq \sqrt{{\varphi \, \mathrm{e}^{\frac{\cU_e +\theta \,\cU_\eta}{\theta}}}}\\\
 \mbox{weakly in $L^1(Q \times \wtD ;\mathbb{R}^{d \times d})$}.     
\end{align*}
By the uniqueness of the weak limit it will then follow that, as $n \to \infty$,
\[ \overline{\sqrt{\theta \, \mathrm{e}^{-\frac{\cU_e +\theta \,\cU_\eta}{\theta}}}\, \nabq \sqrt{{\varphi \, \mathrm{e}^{\frac{\cU_e +\theta \,\cU_\eta}{\theta}}}}} = \sqrt{\theta \, \mathrm{e}^{-\frac{\cU_e +\theta \,\cU_\eta}{\theta}}}\, \nabq \sqrt{{\varphi \, \mathrm{e}^{\frac{\cU_e +\theta \,\cU_\eta}{\theta}}}}, \]
as required. 

To this end, we shall prove that 
\begin{align}\label{eq:38}
 \sqrt{\theta^n \, \mathrm{e}^{-\frac{\cU_e +\theta^n \,\cU_\eta}{\theta^n}}} \to \sqrt{\theta \, \mathrm{e}^{-\frac{\cU_e +\theta \,\cU_\eta}{\theta}}} \quad \mbox{strongly in $L^2(Q \times \wtD)$}
\end{align}
and that 
\begin{align}\label{eq:39}
    \nabq \sqrt{{\varphi^n \, \mathrm{e}^{\frac{\cU_e +\theta^n \,\cU_\eta}{\theta^n}}}} \rightharpoonup \nabq \sqrt{{\varphi \, \mathrm{e}^{\frac{\cU_e +\theta \,\cU_\eta}{\theta}}}}\quad
 \mbox{weakly in $L^2(Q \times \wtD ;\mathbb{R}^{d})$}
\end{align}
as $n \to \infty$.

To prove \eqref{eq:38}, we begin by noting that, by \eqref{eq:str}, $\theta^n \to \theta$ strongly in $L^1(Q)$; also, by \eqref{eq:**}$_1$, $\theta^n \rightharpoonup \theta$  weakly in $L^{\frac{2}{3}+\beta}(Q)$ with $\beta>\frac{5}{6}$. Therefore, $\theta^n \to \theta$ strongly in $L^r(Q)$ for all $r \in [1,\frac{2}{3}+\beta)$. For a subsequence (not indicated) $\theta^n \to \theta$ a.e.~on $Q$. Because $\theta^n \geq \theta_{\min}>0$ on $\overline{Q}$, we have  
\[ \frac{\cU_e + \theta^n\cU_\eta}{\theta^n} = \frac{\cU_e}{\theta^n} + \cU_\eta \in C(\overline{Q \times \wtD}), \]
\[\mathrm{e}^{-\frac{\cU_e + \theta^n\cU_\eta}{\theta^n}} \to \mathrm{e}^{-\frac{\cU_e + \theta\cU_\eta}{\theta}} \quad \mbox{a.e.~on $Q \times \wtD$} \qquad \mbox{and}\qquad 0<\mathrm{e}^{-\frac{\cU_e + \theta^n\cU_\eta}{\theta^n}} \leq 1 \quad \mbox{a.e.~on $Q \times \wtD$}. \]
Thus, by the dominated convergence theorem, 
\[ \mathrm{e}^{-\frac{\cU_e + \theta^n\cU_\eta}{\theta^n}} \to \mathrm{e}^{-\frac{\cU_e + \theta\cU_\eta}{\theta}} \quad \mbox{strongly in $L^1(Q \times \wtD)$}.\]
However, since $\mathrm{e}^{-\frac{\cU_e + \theta^n\cU_\eta}{\theta^n}} \in L^\infty(Q \times \wtD)$ and $\mathrm{e}^{-\frac{\cU_e + \theta\cU_\eta}{\theta}}\in L^\infty(Q \times \wtD)$, this implies that 
\[ 
 \mathrm{e}^{-\frac{\cU_e + \theta^n\cU_\eta}{\theta^n}} \to \mathrm{e}^{-\frac{\cU_e + \theta\cU_\eta}{\theta}} \quad \mbox{strongly in $L^s(Q \times \wtD)$}\quad \mbox{for all $s \in [1,\infty)$},
\]
and, in particular, for $s\in (\frac{3\beta+2}{3\beta-1},\infty)$ (where $\frac{3\beta+2}{3\beta-1}$ is the H\"older conjugate of $\frac{2}{3}+\beta$, with $\beta>\frac{5}{6}$). It therefore follows that the product
\[
 \theta^n \, \mathrm{e}^{-\frac{\cU_e +\theta^n \,\cU_\eta}{\theta^n}} \to \theta \, \mathrm{e}^{-\frac{\cU_e +\theta \,\cU_\eta}{\theta}}\quad \mbox{strongly in $L^1(Q \times \wtD)$}.
\]
This directly implies \eqref{eq:38} using the elementary inequality $|\sqrt{a} - \sqrt{b}| \leq \sqrt{|a-b|}$ for $a,b \in \mathbb{R}_{\geq 0}$.

The proof of \eqref{eq:39} is similar. According to \eqref{eq:32a}, $\varphi^n \to \varphi$ strongly in $L^1(Q \times \wtD)$ and, by \eqref{eq:30} (cf.~also \eqref{eq:29}), $\varphi^n \to \varphi$ weakly in $L^{1+\delta}(Q \times \wtD)$; it follows that $\varphi^n \to \varphi$ strongly in $L^r(Q \times \wtD)$ for all $r \in [1,1+\delta)$ as $n \to \infty$. Also, by the same reasoning as above,
\[ \mathrm{e}^{\frac{\cU_e + \theta^n\cU_\eta}{\theta^n}} \to \mathrm{e}^{\frac{\cU_e + \theta\cU_\eta}{\theta}} \quad \mbox{strongly in $L^s(Q \times \wtD)$}\quad \mbox{for all $s \in [1,\infty)$},\]
and in particular for $s \in (1+\frac{1}{\delta},\infty)$ (where $1+\frac{1}{\delta}$ is the H\"older conjugate of $1+\delta$). Therefore, 
\[
 \varphi^n \, \mathrm{e}^{\frac{\cU_e +\theta^n \,\cU_\eta}{\theta^n}} \to \varphi \, \mathrm{e}^{\frac{\cU_e +\theta \,\cU_\eta}{\theta}}\quad \mbox{strongly in $L^1(Q \times \wtD)$,}
\]
whereby 
\begin{align}\label{eq:40}
 \sqrt{\varphi^n \, \mathrm{e}^{\frac{\cU_e +\theta^n \,\cU_\eta}{\theta^n}}} \to \sqrt{\varphi \, \mathrm{e}^{\frac{\cU_e +\theta \,\cU_\eta}{\theta}}}\quad \mbox{strongly in $L^2(Q \times \wtD)$}
\end{align}
as $n \to \infty$. It follows from \eqref{eq:hat} and the lower bound
\[ 
\theta^n \, \mathrm{e}^{-\frac{\cU_e + \theta^n\, \cU_\eta}{\theta^n}} \geq \theta_{\min} 
\mathrm{e}^{-\frac{\cU_e}{\theta_{\min}} - \cU_\eta}  \geq C_{\wtD}>0 \quad \mbox{on $\overline{Q \times \wtD}$},
\]
guaranteed by the uniform continuity of $\cU_e$ and $\cU_\eta$ on $\overline{\wtD}$, that 
\begin{align*}
\int_{Q \times \wtD}
|\nabq \sqrt{\varphi^{n}  \, \mathrm{e}^{\frac{\cU_e + \theta^n\, \cU_\eta}{\theta^n}}  }|^2\dd \bq \dd \bx \dd t  \leq C_{\wtD}
\end{align*}
for all $n \geq 1$. Hence, for a subsequence (not indicated), 
\[ \nabq \sqrt{\varphi^{n}  \, \mathrm{e}^{\frac{\cU_e + \theta^n\, \cU_\eta}{\theta^n}}} \to \overline{ \nabq \sqrt{\varphi^{n}  \, \mathrm{e}^{\frac{\cU_e + \theta^n\, \cU_\eta}{\theta^n}}}}  \quad \mbox{weakly in $L^2(Q \times \wtD;\mathbb{R}^d)$}
\] 
as $n \to \infty$. The uniqueness of the weak limit and \eqref{eq:40} then imply that \eqref{eq:39} holds.

Having shown both \eqref{eq:38} and \eqref{eq:39}, the assertion \eqref{eq:36} follows thanks to the argument described in the paragraph between equations \eqref{eq:37} and \eqref{eq:38}.

This then completes the proofs of \eqref{eq:35} and \eqref{eq:36}. Using \eqref{eq:35} and \eqref{eq:36}, we can now pass to the limit in \eqref{eq:34a}; the argument proceeds as follows: 
\begin{align*}
&\lim_{n \to \infty} \int_Q \left(\int_D (\bF^n \otimes \bq)\varphi^n \dd \bq \right)  : \nabx \bw \dd \bx \dd t \\ 
&= 2\int_Q \left( \int_D \sqrt{\theta}\, \sqrt{\varphi} \left(\bq \otimes  
 \sqrt{\theta \, \mathrm{e}^{-\frac{\cU_e +\theta \,\cU_\eta}{\theta}}}\, \nabq \sqrt{{\varphi \, \mathrm{e}^{\frac{\cU_e +\theta \,\cU_\eta}{\theta}}}}\right) \dd \bq\right) : \nabx \bw \dd \bx \dd t\\
&= 2\int_Q \left( \int_D \sqrt{\theta}\, \sqrt{\varphi} \left(  
 \sqrt{\theta \, \mathrm{e}^{-\frac{\cU_e +\theta \,\cU_\eta}{\theta}}}\, \nabq \sqrt{{\varphi \, \mathrm{e}^{\frac{\cU_e +\theta \,\cU_\eta}{\theta}}}} \otimes \bq\right) \dd \bq\right) : \nabx \bw \dd \bx \dd t\\
%
&= 2\int_Q \left( \int_D \sqrt{\frac{\theta}{\varphi}}  \left(  
 \sqrt{\theta \, \mathrm{e}^{-\frac{\cU_e +\theta \,\cU_\eta}{\theta}}}\, \nabq \sqrt{{\varphi \, \mathrm{e}^{\frac{\cU_e +\theta \,\cU_\eta}{\theta}}}}\otimes \bq\right)\varphi \dd \bq\right) : \nabx \bw \dd \bx \dd t\\
 &= 2\int_Q \left( \int_D \theta  \left( \frac{\nabq \sqrt{\widehat{\varphi}}}{\sqrt{\widehat{\varphi}}}\otimes \bq\right)\varphi \dd \bq\right) : \nabx \bw \dd \bx \dd t\\
 &= \int_Q\left(  \theta \int_D \left( \frac{\nabq \widehat{\varphi}}{\widehat{\varphi}}\otimes \bq\right)\varphi \dd \bq\right) : \nabx \bw \dd \bx \dd t = \int_Q \left( \theta \int_D \left(\nabq \log \widehat{{\varphi}} \otimes \bq\right)\varphi \dd \bq \right):\nabx \bw \dd \bx \dd t.
\end{align*}
Comparing the expression on the right-hand side of this equality with the right-hand side of \eqref{eq:19b} and recalling that, by hypothesis, $\bII : \nabx \bw = \divx \bw = 0$, it follows that
\begin{align*}
    \lim_{n \to \infty} \int_Q \left(\int_D (\bF^n \otimes \bq)\varphi^n \dd \bq \right)  : \nabx \bw \dd \bx \dd t  &= 
    \int_Q \left( \theta \int_D \left(\nabq \log \widehat{{\varphi}} \otimes \bq\right)\varphi \dd \bq \right):\nabx \bw \dd \bx \dd t\\
    &=
    \int_Q \left(\int_D (\bF \otimes \bq)\varphi \dd \bq \right)  : \nabx \bw \dd \bx \dd t
\end{align*}
for all $\bw \in C([0,T];C^1_0(\overline{\Omega};\mathbb{R}^d))$ such that $\divx\bw=0$ on $Q$. Thus, we have shown that 
\[ \overline{\bSS} =2\nu(\theta) \bDD(\bv) + \int_D (\boldsymbol{F} \otimes \bq) \varphi  \dd \bq = \bSS, \]
whereby, thanks to the density of the set of all divergence-free $C^1_0(\overline\Omega;\mathbb{R}^d)$ functions in $W^{1,s}_{0,\mathrm{div}}(\Omega;\mathbb{R}^d)$, $s>1$, we have from \eqref{eq:E1F} that, for a.e.~$t \in (0,T)$, 
\begin{align*}
- \int_Q \bv\cdot \partial_t \bw \dd \bx \dd t- \int_Q (\bv \otimes \bv) : \nabx \bw \dd \bx \dd t + \int_Q \bSS:\mathbb{D}(\bw)  \dd \bx \dd t = \int_\Omega \bv_0 \cdot \bw(0) \dd \bx + \int_Q \bff \cdot \bw \dd \bx \dd t
\end{align*}
for all test functions $\bw \in C^1_0([0,T); W^{1,s}_{0,\mathrm{div}}(\Omega;\mathbb{R}^d))$, 
for some $s>1$.

\subsection{Limit passage in the Fokker--Planck equation}

Next we consider passing to the limit in the Fokker--Planck equation \eqref{eq:m3} satisfied by $\varphi^n$ (and with $\bv$, $\theta$ and $\bj_{\varphi,\bq}$ there replaced by $\bv^n$, $\theta^n$ and $\bj_{\varphi,\bq}^n$, respectively). To this end, we need to pass to the respective limits in the sequences
\[ \bv^n\varphi^n, \quad \theta^n \nabx \varphi^n, \quad (\nabx \bv^n)\bq \varphi^n, \quad \mbox{and}\quad \bj_{\varphi}^n.\]

We begin by considering the term in the Fokker--Planck equation involving $\bv^n\varphi^n$. Let $\zeta \in C([0,T];C^1(\overline{\Omega} \times \overline{D}))=C([0,T];C^1(\overline{\Omega \times D}))$. Thanks to \eqref{eq:33a}, we have that
\[ \bz^n:= \int_D\varphi^n \nabx \zeta \dd \bq \to \int_D\varphi \nabx \zeta \dd \bq=: \bz \quad \mbox{strongly in $L^1(Q;\mathbb{R}^d)$}\]
as $n \rightarrow \infty$. Because 
\[ \|\bz^n\|_{L^\infty(Q)}\leq \|\nabx \zeta\|_{L^\infty(Q \times D)} \|\varphi^n\|_{L^\infty(Q;L^1(D))} \leq C\]
and
\[\|\bz\|_{L^\infty(Q)}\leq \|\nabx \zeta\|_{L^\infty(Q \times D)} \|\varphi\|_{L^\infty(Q;L^1(D))} \leq C, \]
~\\
it follows by H\"older's inequality that $\bz^n \to \bz$ strongly in $L^s(Q;\mathbb{R}^d)$ for all $s \in [1,\infty)$, and in particular for $s=\frac{r}{r-1}$ where $r=\frac{10}{3} \frac{2+3\beta}{5+3\beta}$ and $\beta>\frac{5}{6}$. Then it follows further from \eqref{eq:va} that $\bv^n \cdot \bz^n \rightarrow \bv \cdot \bz$ strongly in $L^1(Q)$. Therefore, 
\[ \lim_{n\to \infty}\int_{Q \times D} \bv^n \varphi^n \cdot \nabx\zeta \dd \bq \dd \bx \dd t = \int_{Q \times D} \bv\varphi \cdot \nabx\zeta \dd \bq \dd \bx \dd t\]
for all $\zeta \in C([0,T];C^1(\overline{\Omega \times D}))$.

\medskip

To pass to the limit in the term that involves $\theta^n \nabx \varphi^n$, we note that 
\[ \theta^n \nabx \varphi^n = \big(\sqrt{\theta^n} \sqrt{\varphi^n} \big) \left(2\sqrt{\theta^n} \nabx \sqrt{\varphi^n}\right).
\]
According to \eqref{eq:35}, the first factor strongly converges to $\sqrt{\theta} \sqrt{\varphi}$ in $L^2(Q \times D)$. By virtue of \eqref{eq:26}, 
\[ \|2\sqrt{\theta^n} \nabx \sqrt{\varphi^n}\|_{L^2(Q \times D)} =  \left\|\sqrt{\frac{\theta^n}{\varphi^n}} \nabx \varphi^n\right\|_{L^2(Q \times D)} \leq C\]
for all $n \geq 1$, with $C$ independent of $n$, and therefore, for a subsequence (not indicated),
\[ 2\sqrt{\theta^n} \nabx \sqrt{\varphi^n} \rightharpoonup 2\overline{\sqrt{\theta} \nabx \sqrt{\varphi}}\quad \mbox{weakly in $L^2(Q \times D;\mathbb{R}^d)$} \]
as $n \to \infty$, where $\overline{\sqrt{\theta} \nabx \sqrt{\varphi}} \in L^2(Q \times D)$ is a weak limit to be identified. Hence, 
also
\[ 2\sqrt{\theta^n} \nabx \sqrt{\varphi^n}  \rightharpoonup 2\overline{\sqrt{\theta} \nabx \sqrt{\varphi}}\quad \mbox{weakly in $L^1(Q \times \wtD\mathbb{R}^d))$} \]
as $n \to \infty$ for any $\wtD \subset\joinrel\subset D$, and the weak limit on the right-hand side is independent of the choice of $\wtD$. By \eqref{eq:str}, $\theta^n \to \theta$ strongly in $L^1(Q)$ and therefore $\sqrt{\theta^n} \to \sqrt{\theta}$ strongly in $L^2(Q)$ as $n \to \infty$. Also, thanks to \eqref{eq:28a} applied to $\varphi^n$, 
\[ \|\nabx \sqrt{\varphi^n}\|_{L^2(Q \times \wtD)} \leq C_{\wtD}\]
for all $n \geq 1$ and for any $\wtD \subset\joinrel\subset D$. Therefore, there exists a subsequence (not indicated) such that $\nabx \sqrt{\varphi^n}$ is weakly convergent in $L^2(Q \times \wtD;\mathbb{R}^d)$. Because by \eqref{eq:33a} $\varphi^n \to \varphi$ strongly in $L^1(Q \times D)$ and therefore $\sqrt{\varphi^n} \to \sqrt{\varphi}$ strongly in $L^2(Q \times D)$, also $\sqrt{\varphi^n} \to \sqrt{\varphi}$ strongly in $L^2(Q \times \wtD)$. The uniqueness of the weak limit then implies that the weak limit of $\nabx \sqrt{\varphi^n}$ must be $\nabx \sqrt{\varphi}$; that is, 
\[ \nabx \sqrt{\varphi^n} \rightharpoonup \nabx \sqrt{\varphi}\quad \mbox{weakly in $L^2(Q \times \wtD;\mathbb{R}^d)$}\]
as $n \to \infty$. Thus, we deduce that 
\[ 2\sqrt{\theta^n} \nabx \sqrt{\varphi^n}  \rightharpoonup 2\sqrt{\theta} \nabx \sqrt{\varphi}\quad \mbox{weakly in $L^1(Q \times \wtD;\mathbb{R}^d)$} \]
as $n \to \infty$ for any $\wtD \subset\joinrel\subset D$. Therefore, because of the independence of the weak limit of the choice of $\wtD$, 
\begin{align}\label{eq:tf} 
2\sqrt{\theta^n} \nabx \sqrt{\varphi^n} \rightharpoonup 2 \sqrt{\theta} \nabx \sqrt{\varphi} \quad \mbox{weakly in $L^2(Q \times D;\mathbb{R}^d)$}.
\end{align}
Since $\sqrt{\theta^n} \sqrt{\varphi^n} \to \sqrt{\theta} \sqrt{\varphi}$ strongly in $L^2(Q \times D)$, it follows that 
\[ 
\theta^n \nabx \varphi^n= \big(\sqrt{\theta^n} \sqrt{\varphi^n} \big) \left(2\sqrt{\theta^n} \nabx \sqrt{\varphi^n}\right)
\rightharpoonup \big(\sqrt{\theta} \sqrt{\varphi} \big) \left(2\sqrt{\theta} \nabx \sqrt{\varphi}\right) =  \theta \nabx \varphi
\]
weakly in $L^1(Q \times D;\mathbb{R}^d)$ as $n \to \infty$. Hence, 
\[
\lim_{n \to \infty}  \int_{Q \times D} \theta^n \nabx \varphi^n \cdot \nabx \zeta \dd \bq \dd \bx \dd t = \int_{Q \times D} \theta \nabx \varphi \cdot \nabx \zeta \dd \bq \dd \bx \dd t 
\]
for all $\zeta \in C([0,T];C^1(\overline{\Omega \times D}))$.

\medskip

Next, we consider the term in the Fokker--Planck equation containing $(\nabx \bv^n)\bq \varphi^n$. Using an identical argument as in the case of the term $\bv^n\varphi^n$ discussed above, thanks to \eqref{eq:33} we have the following:
\begin{align}\label{eq:tg} 
\mathbb{Z}^n:= \int_D(\nabx \zeta \otimes \bq)\varphi^n \dd \bq \to \int_D(\nabx \zeta \otimes \bq)\varphi \dd \bq=: \mathbb{Z} \quad \mbox{strongly in $L^s(Q;\mathbb{R}^{d \times d})$}
\end{align}
as $n \rightarrow \infty$ for all $\zeta \in C([0,T];C^1(\overline{\Omega \times D}))$, for all $s \in [1,\infty)$,  and in particular for $s=\frac{r}{r-1}$, where $r = \frac{4+6\beta}{5+3\beta}$, $\beta>\frac{5}{6}$. Hence, thanks to \eqref{eq:v}, 
\[ \lim_{n \to \infty}\int_Q \nabx \bv^n : \mathbb{Z}^n \dd \bx \dd t = \int_Q \nabx \bv : \mathbb{Z} \dd \bx \dd t\]
as $n \to \infty$. Consequently,
\[\lim_{n \to \infty} \int_{Q \times D}(\nabx \bv^n)\bq\varphi^n \cdot \nabx \zeta \dd \bq \dd \bx \dd t =  \int_{Q \times D}(\nabx \bv)\bq\varphi \cdot \nabx \zeta \dd \bq \dd \bx \dd t\]
as $n \to \infty$ for all $\zeta \in C([0,T];C^1(\overline{\Omega \times D}))$.

\medskip

It remains to pass to the limit in the term containing $\bj_{\varphi}^n$ in the Fokker--Planck equation. We write
\begin{align}\label{eq:jna}
\bj_{\varphi}^n = \left(\sqrt{\theta^n} \sqrt{\varphi^n} \right) \left( \frac{1}{\sqrt{\theta^n \varphi^n}} \bj_{\varphi}^n \right).
\end{align}
According to \eqref{eq:35}, $\sqrt{\theta^n} \sqrt{\varphi^n} \to \sqrt{\theta} \sqrt{\varphi}$ strongly in $L^2(Q \times D)$, and, thanks to the bound \eqref{eq:jbd}, also,
\begin{align}\label{eq:jn0}
\frac{1}{\sqrt{\theta^n \varphi^n}}\bj^n_\varphi \rightharpoonup \overline{\frac{1}{\sqrt{\theta \varphi}}\bj_{\varphi,\bq}}\quad \mbox{weakly in $L^2(Q \times D;\mathbb{R}^d)$},
\end{align}
and therefore also weakly in $L^2(Q \times \wtD;\mathbb{R}^d)$ for any $\wtD \subset\joinrel\subset D$, and the weak limit, to be identified,  is independent of the choice of $\wtD$; we shall use this observation to identify the weak limit. We begin by noting that, thanks to the definitions of $\bj^n_\varphi$ and $\boldsymbol{F}^n$, we have 
\begin{align}\label{eq:jn}
\frac{1}{\sqrt{\theta^n \varphi^n}}\bj^n_\varphi = -\left[\sqrt{\varphi^n}\frac{1}{\sqrt{\theta^n}}\,(\nabq\, \cU_e) + \sqrt{\theta^n}\sqrt{\varphi^n} \,(\nabq\,\cU_\eta) + 2 \sqrt{\theta^n}\nabq \sqrt{\varphi^n} \right].
\end{align}
By \eqref{eq:tf}, the final term on the right-hand side of \eqref{eq:jn} satisfies
\[
2\sqrt{\theta^n} \nabx \sqrt{\varphi^n} \rightharpoonup 2 \sqrt{\theta} \nabx \sqrt{\varphi} \quad \mbox{weakly in $L^2(Q \times D;\mathbb{R}^d)$}
\]
as $n \to \infty$. Since $\cU_\eta \in C^1(\overline{D})$, it follows from \eqref{eq:36} that for the middle term on the right-hand side of \eqref{eq:jn} we have
\[  \sqrt{\theta^n}\, \sqrt{\varphi^n}\, (\nabq\,\cU_\eta) \to \sqrt{\theta}\, \sqrt{\varphi}\, (\nabq\,\cU_\eta)  \quad \mbox{strongly in $L^2(Q \times D;\mathbb{R}^d)$}.\]
Therefore, it remains to consider the first term on the right-hand side of \eqref{eq:jn}. Thanks to \eqref{eq:33a}, $\sqrt{\varphi^n} \to \sqrt{\varphi}$ strongly in $L^2(Q \times D)$. Furthermore, since $\theta^n \geq \theta_{\min}\, \mathrm{e}^{-\alpha T}>0$ with $\alpha>0$ for all $n \geq 1$ (cf. \eqref{eq:20}), whereby $0<1/\theta^n \leq \mathrm{e}^{\alpha T}/\theta_{\min}$ a.e.~on $Q$, the a.e.~convergence of $\theta^n$ to $\theta$ on $Q$ asserted in \eqref{eq:thp} implies by the dominated convergence theorem that $1/\theta^n \to 1/\theta$  strongly in $L^1(Q \times D)$. Therefore,  $1/\sqrt{\theta^n} \to 1/\sqrt{\theta}$  strongly in $L^2(Q \times D)$. Hence, the product $\sqrt{\varphi^n}/\sqrt{\theta^n}$ converges strongly to $\sqrt{\varphi}/\sqrt{\theta}$ in $L^1(Q \times D)$. Furthermore, because $(1/\sqrt{\theta^n})_{n=1}^\infty$ is a bounded sequence in $L^\infty(Q)$, it has a weakly-$\ast$ convergent subsequence (not indicated), which by the uniqueness of the weak limit must converge to $1/\sqrt{\theta}$. Therefore, 
\[\sqrt{\varphi^n}\frac{1}{\sqrt{\theta^n}} \rightharpoonup \sqrt{\varphi}\frac{1}{\sqrt{\theta}} \quad \mbox{weakly in $L^2(Q \times D)$}.\]
Because $\cU_e \in C^1(\wtD)$ for each $\wtD \subset\joinrel\subset D$, we therefore deduce that 
\[\sqrt{\varphi^n}\frac{1}{\sqrt{\theta^n}}\,(\nabq\, \cU_e) \rightharpoonup \sqrt{\varphi}\frac{1}{\sqrt{\theta}}\,(\nabq\, \cU_e)\quad \mbox{weakly in $L^2(Q \times \wtD;\mathbb{R}^d)$}\]
as $n \to \infty$ for each $\wtD \subset\joinrel\subset D$. We are now ready to pass to the limit in \eqref{eq:jn}:
\begin{align*}
\frac{1}{\sqrt{\theta^n \varphi^n}}\bj^n_\varphi &= -\left[\sqrt{\varphi^n}\frac{1}{\sqrt{\theta^n}}\,(\nabq\, \cU_e) + \sqrt{\theta^n}\sqrt{\varphi^n} \,(\nabq\,\cU_\eta) + 2 \sqrt{\theta^n}\nabq \sqrt{\varphi^n} \right]\\
&\rightharpoonup  -\left[\sqrt{\varphi}\frac{1}{\sqrt{\theta}}\,(\nabq\, \cU_e) + \sqrt{\theta}\sqrt{\varphi} \,(\nabq\,\cU_\eta) + 2 \sqrt{\theta}\nabq \sqrt{\varphi} \right] = \frac{1}{\sqrt{\theta \varphi}}\bj_{\varphi,\bq}
\end{align*}
weakly in $L^2(Q \times \wtD;\mathbb{R}^d)$ as $n \to \infty$ for each $\wtD \subset\joinrel\subset D$. Hence, because of the uniqueness of the weak limit, \eqref{eq:jn0} implies that 
\begin{align}\label{eq:tfj}
\frac{1}{\sqrt{\theta^n \varphi^n}}\bj^n_\varphi 
&\rightharpoonup \frac{1}{\sqrt{\theta \varphi}}\bj_{\varphi,\bq}\quad \mbox{weakly in $L^2(Q \times D;\mathbb{R}^d)$}
\end{align}
as $n \to \infty$. Returning with this to \eqref{eq:jna} and recalling that $\sqrt{\theta^n} \sqrt{\varphi^n} \to \sqrt{\theta} \sqrt{\varphi}$ strongly in $L^2(Q \times D)$, we finally find that 
\[ \bj_{\varphi}^n \rightharpoonup \bj_{\varphi} \quad \mbox{weakly in $L^1(Q \times D;\mathbb{R}^d)$}\]
as $n \to \infty$. Therefore, 
\[\lim_{n \to \infty} \int_{Q \times D} \bj_{\varphi,\bq}^n \cdot \nabx \zeta \dd \bq \dd \bx \dd t =  \int_{Q \times D} \bj_{\varphi,\bq}\cdot \nabx \zeta \dd \bq \dd \bx \dd t\]
for all $\zeta \in C([0,T];C^1(\overline{\Omega \times D}))$. Having passed to the limits in the sequences 
\[ \bv^n\varphi^n, \quad \theta^n \nabx \varphi^n, \quad 
(\nabx \bv^n)\bq \varphi^n \quad \mbox{and}\quad \bj_{\varphi}^n,\]
one can now pass to the limit $n \to \infty$ in the weak formulation of the Fokker--Planck equation satisfied by the function $\varphi^n$ to obtain:
\begin{align*}
& - \int_{Q \times D} \varphi\, \partial_t\zeta \dd \bq \dd \bx \dd t - \int_{Q \times D} \bv \varphi \cdot \nabx \zeta \dd \bq \dd \bx \dd t  +  \int_{Q \times D}\theta \nabx \varphi \cdot \nabx \zeta \dd \bq \dd \bx \dd t \\
&\quad - \int_{Q \times D} ((\nabx \bv) \bq \varphi + \bj_{\varphi,\bq})\cdot \nabq \zeta \dd \bq \dd \bx \dd t = \int_{\Omega \times D} \varphi_0 \,\zeta(0) \dd \bq \dd \bx
\qquad \forall\, \zeta \in C^1_0([0,T);C^1(\overline{\Omega \times D})).
\end{align*}

\subsection{Limit passage in the renormalized temperature equation}

Finally, we pass to the limit $n \to \infty$ in the renormalized temperature equation \eqref{eq:ren}. As noted in the paragraph following \eqref{eq:*r}, we already have in place the required convergence results to pass to the limit $n \to \infty$ in the weak form of the left-hand side 
\[ \partial_t T_{k,\varepsilon}(\theta^n) + \divx(\bv^n\, T_{k,\varepsilon}(\theta^n)) - \divx(\kappa(\theta^n) \nabx T_{k,\varepsilon}(\theta^n))\]
of \eqref{eq:ren}. Therefore, it remains to deal with the various terms on the right-hand side of \eqref{eq:ren}. Motivated by the analysis in  \cite{BFM}, we shall pass to the limit in those terms from one side; this will result in a renormalized variational \textit{inequality} satisfied by the absolute temperature. To this end, we take $\psi \in C([0,T]; W^{1,z}(\Omega))$ as a test function for \eqref{eq:ren} with $z>d$ (whereby, thanks to the Sobolev embedding $W^{1,z}(\Omega) \hookrightarrow C(\overline{\Omega})$, $\psi \in C(\overline{Q})$), and assume that $\psi \geq 0$ on $\overline{Q}$. We define
\begin{align*}
\mathcal{T}_1 &:=\int_{Q \times D} \left(2\nu(\theta^n) T_{k,\varepsilon}'(\theta^n)|\bDD(\bv^n)|^2 - T_{k,\varepsilon}''(\theta^n)\kappa(\theta^n)|\nabx\theta^n|^2 \right) \psi \dd \bq \dd \bx \dd t,  \\
\mathcal{T}_2 &:= \int_{Q \times D} \left(T_{k,\varepsilon}'(\theta^n) \bDD(\bv^n):\theta^n (\nabq\, \cU_\eta \otimes \bq) \varphi^n \right) \psi \dd \bq \dd \bx \dd t, \\
\mathcal{T}_3 &:= \int_{Q \times D} \left( T_{k,\varepsilon}'(\theta^n) \theta^n  \varphi^n \nabq\, \cU_\eta \cdot \nabq\, \cU_e \right) \psi \dd \bq \dd \bx \dd t,  \\
\mathcal{T}_4 &:= \int_{Q \times D}  T_{k,\varepsilon}'(\theta^n)\left( \theta^n \nabq \varphi^n \cdot \nabq\, \cU_e + |\nabq\, \cU_e|^2 \varphi^n\right)  \psi \dd \bq \dd \bx \dd t.
\end{align*}
We note that these integrals are well-defined. Specifically, $\mathcal{T}_3$ and $\mathcal{T}_4$, which involve the singular potential $\cU_e$, are finite thanks to \eqref{eq:up}, with $\theta$ and $\varphi$ replaced there by $\theta^n$ and $\varphi^n$, respectively. 

We start by considering the term $\mathcal{T}_1$. As $T'_{\kappa,\varepsilon}(s) \geq 0$ and $T''_{\kappa,\varepsilon}(s) \leq 0$ for all $s \geq 0$, it follows that the integrand of $\mathcal{T}_1$ is nonnegative.  To pass to the limit in $\mathcal{T}_1$ as $n \to \infty$, note that 
\begin{align*}
2\nu(\theta^n) T_{k,\varepsilon}'(\theta^n)\left|\bDD(\bv^n)\right|^2  \psi &= 2\nu(\theta^n) T_{k,\varepsilon}'(\theta^n)\,\theta^n \left|\frac{\bDD(\bv^n)}{\sqrt{\theta^n}}- \frac{\bDD(\bv)}{\sqrt{\theta}} + \frac{\bDD(\bv)}{\sqrt{\theta}}\right|^2  \psi \\
&= 2\nu(\theta^n) T_{k,\varepsilon}'(\theta^n)\, \theta^n 
\left|\frac{\bDD(\bv^n)}{\sqrt{\theta^n}}- \frac{\bDD(\bv)}{\sqrt{\theta}}\right|^2\, \psi + 2\nu(\theta^n) T_{k,\varepsilon}'(\theta^n)\, \theta^n \left| \frac{\bDD(\bv)}{\sqrt{\theta}}\right|^2 \psi \\
&\quad + 4\nu(\theta^n) T_{k,\varepsilon}'(\theta^n) \,\theta^n
\left(\frac{\bDD(\bv^n)}{\sqrt{\theta^n}}- \frac{\bDD(\bv)}{\sqrt{\theta}}\right) :\frac{\bDD(\bv)}{\sqrt{\theta}}\, \psi  \\
& \hspace{-14mm} \geq 2\nu(\theta^n) T_{k,\varepsilon}'(\theta^n)\, \theta^n \left| \frac{\bDD(\bv)}{\sqrt{\theta}}\right|^2 \psi + 4\nu(\theta^n) T_{k,\varepsilon}'(\theta^n)\, \theta^n
\left(\frac{\bDD(\bv^n)}{\sqrt{\theta^n}}- \frac{\bDD(\bv)}{\sqrt{\theta}}\right) :\frac{\bDD(\bv)}{\sqrt{\theta}}\, \psi.
\end{align*}
As $n \to \infty$, the first term on the right-hand side converges to 
$2\nu(\theta) T_{k,\varepsilon}'(\theta)|\bDD(\bv)|^2\psi$
strongly in $L^1(Q \times D)$ using \eqref{eq:thp} by the dominated convergence theorem, while the second term on the right-hand side converges weakly in $L^1(Q \times D)$ to $0$ because of \eqref{eq:thp}, the boundedness of the sequence $(4\nu(\theta^n) T_{k,\varepsilon}'(\theta^n)\, \theta^n \psi)_{n=1}^\infty$ in $L^\infty(Q \times D)$ and the weak convergence
\begin{align} \label{eq:L2w}
\left(\frac{\bDD(\bv^n)}{\sqrt{\theta^n}}- \frac{\bDD(\bv)}{\sqrt{\theta}}\right) :\frac{\bDD(\bv)}{\sqrt{\theta}}\rightharpoonup 0 \quad \mbox{in $L^1(Q \times D)$},
\end{align}
using Proposition 2.61 on p.183 in \cite{Fonseca}.
The weak convergence result \eqref{eq:L2w} follows by recalling \eqref{eq:thp} and the fact that the sequence $(\theta^n)_{n=1}^\infty$ is equi-bounded from below by a positive constant, which, by the dominated convergence theorem, together guarantee that the sequence $(1/\sqrt{\theta^n})_{n=1}^\infty$ converges strongly to $(1/\sqrt{\theta})$ in $L^s(Q)$ for all $s \in [1,\infty)$, and in particular for $s = \frac{4+6\beta}{3\beta-1}$, where $\beta>1/3$. On the other hand, by \eqref{eq:v}$_3$, $\bDD(\bv^n) \rightharpoonup \bDD(\bv)$ weakly in $L^{\frac{4+6\beta}{5+3\beta}}(Q;\mathbb{R}^{d \times d})$, and therefore $\bDD(\bv^n)/\sqrt{\theta^n} \rightharpoonup \bDD(\bv)/\sqrt{\theta}$ weakly in $L^1(Q;\mathbb{R}^{d \times d})$. However, $(\bDD(\bv^n)/\sqrt{\theta^n})_{n=1}^\infty$ is a bounded sequence in $L^2(Q;\mathbb{R}^{d \times d})$ and therefore has a weakly convergent subsequence in $L^2(Q;\mathbb{R}^{d \times d})$; by the uniqueness of the weak limit (for a subsequence, not indicated), $\bDD(\bv^n)/\sqrt{\theta^n} \rightharpoonup \bDD(\bv)/\sqrt{\theta}$ weakly in $L^2(Q;\mathbb{R}^{d \times d})$. This then implies the weak convergence result \eqref{eq:L2w} asserted above.

Therefore, by the weak lower-semicontinuity of the $L^1(Q \times D)$ norm, 
\[ \mbox{lim.inf}_{n \to \infty} \int_{Q \times D}
2\nu(\theta^n) T_{k,\varepsilon}'(\theta^n)\left|\bDD(\bv^n)\right|^2  \psi \dd \bq \dd \bx \dd t \geq 
\int_{Q \times D} 2\nu(\theta) T_{k,\varepsilon}'(\theta)\left|\bDD(\bv)\right|^2  \psi\dd \bq \dd \bx \dd t
\]
which completes the passage to the limit in the first term of $\mathcal{T}_1$. 

Passage to the limit in the second term of $\mathcal{T}_1$ proceeds analogously, by writing
\begin{align*}
- T_{k,\varepsilon}''(\theta^n)\kappa(\theta^n)|\nabx\theta^n|^2 \psi 
&=  - T_{k,\varepsilon}''(\theta^n)\kappa(\theta^n)\, \theta^n 
\left|\frac{\nabx\theta^n}{\theta^n}-\frac{\nabx\theta}{\theta} + \frac{\nabx\theta}{\theta}\right|^2 \psi\\
&\hspace{-10mm} = - T_{k,\varepsilon}''(\theta^n)\kappa(\theta^n)\, \theta^n 
\left|\nabx \log \theta^n -\nabx\log \theta + \nabx \log \theta \right|^2 \psi\\
&\hspace{-10mm} = - T_{k,\varepsilon}''(\theta^n)\kappa(\theta^n)\, \theta^n 
|\nabx \log \theta^n -\nabx\log \theta|^2 \psi  
- T_{k,\varepsilon}''(\theta^n)\kappa(\theta^n)\, \theta^n 
\left|\nabx \log \theta \right|^2 \psi\\
&\hspace{-10mm}\quad 
- 2T_{k,\varepsilon}''(\theta^n)\kappa(\theta^n)\, \theta^n\, \psi
\left(\nabx \log \theta^n -\nabx\log \theta\right) \cdot \nabx \log \theta
\\
&\hspace{-10mm} 
\geq - T_{k,\varepsilon}''(\theta^n)\kappa(\theta^n)\, \theta^n 
\left|\nabx \log \theta \right|^2 \psi\\
&\hspace{-10mm}\quad 
- 2T_{k,\varepsilon}''(\theta^n)\kappa(\theta^n)\, \theta^n\,  \psi 
\left(\nabx \log \theta^n -\nabx\log \theta\right) \cdot \nabx \log \theta.
\end{align*}
As $n \to \infty$, the first term on the right-hand side converges to 
$- T_{k,\varepsilon}''(\theta)\kappa(\theta)\, \theta
\left|\nabx \log \theta \right|^2 \psi$  strongly in $L^1(Q\times D)$ by 
\eqref{eq:thp} and the dominated convergence theorem, while the second term converges to 0 weakly in $L^1(Q \times D)$ by \eqref{eq:logfi}.
Therefore, by the weak lower-semicontinuity of the $L^1(Q \times D)$ norm, 
\[ 
\mbox{lim.inf}_{n \to \infty} \int_{Q \times D}
- T_{k,\varepsilon}''(\theta^n)\kappa(\theta^n)|\nabx\theta^n|^2 \psi
\dd \bq \dd \bx \dd t \geq  \int_{Q \times D}
- T_{k,\varepsilon}''(\theta)\kappa(\theta)|\nabx\theta|^2 \psi
\dd \bq \dd \bx \dd t.\]
Because the $\mathrm{lim.inf}_{n \to \infty}$ of the sum of the two terms is larger than the sum of the two $\mathrm{lim.inf}_{n \to \infty}$ terms, we have that 
\begin{align*}
\mathrm{lim.inf}_{n \to \infty} \int_{Q \times D} \left(2\nu(\theta^n) T_{k,\varepsilon}'(\theta^n)|\bDD(\bv^n)|^2 - T_{k,\varepsilon}''(\theta^n)\kappa(\theta^n)|\nabx\theta^n|^2 \right) \psi \dd \bq \dd \bx \dd t \\
 \quad \geq \int_{Q \times D} \left(2\nu(\theta) T_{k,\varepsilon}'(\theta)|\bDD(\bv)|^2 - T_{k,\varepsilon}''(\theta)\kappa(\theta)|\nabx\theta|^2 \right) \psi \dd \bq \dd \bx \dd t,
\end{align*}
which completes the passage to the limit $n \to \infty$ in the term $\mathcal{T}_1$.

\smallskip

For the term $\mathcal{T}_2$, we use that $0 \leq T'_{k,\varepsilon}(s) s \leq k$ for all $s \geq 0$. We also note that because $\cU_\eta \in C^1(\overline{D})$, we have, analogously as in \eqref{eq:tg},
\[ \int_D (\nabq\, \cU_\eta \otimes \bq) \varphi^n \psi \dd \bq \to \int_D (\nabq\, \cU_\eta \otimes \bq) \varphi \psi \dd \bq 
\quad \mbox{strongly in $L^s(Q;\mathbb{R}^{d \times d})$}
\]
as $n \rightarrow \infty$ for all $s \in [1,\infty)$,  and in particular for $s=\frac{r}{r-1}$ where $r=\frac{4+6\beta}{5+3\beta}$ with $\beta>\frac{5}{6}$. Then, it directly follows from \eqref{eq:v}$_2$ and \eqref{eq:thp}, again invoking Proposition 2.61 on p.183 in \cite{Fonseca}, that 
\begin{align*}
\lim_{n \to \infty} \int_{Q \times D} \left(T_{k,\varepsilon}'(\theta^n) \bDD(\bv^n):\theta^n (\nabq\, \cU_\eta \otimes \bq) \varphi^n \right) \psi \dd \bq \dd \bx \dd t\\
= \int_{Q \times D} \left(T_{k,\varepsilon}'(\theta) \bDD(\bv):\theta (\nabq\, \cU_\eta \otimes \bq) \varphi\right) \psi \dd \bq \dd \bx \dd t.
\end{align*}

For the term $\mathcal{T}_3$ we use that $T'_{k,\varepsilon}(s) \geq 0$ for all $s \geq 0$ and $\nabq\, \cU_\eta \cdot \nabq\, \cU_e \geq 0$ on $D$ (because $U_\eta$ and $U_e$ have been assumed to be monotonically increasing functions on $[0, b/2)$, whereby $\nabq\, \cU_\eta \cdot \nabq\, \cU_e = |\bq|^2\, U_\eta'(|\bq|^2/2)\, U_e'(|\bq|^2/2)\geq 0$ on $D$). Also, thanks to \eqref{eq:thp} and \eqref{eq:33},
\[T_{k,\varepsilon}'(\theta^n)\, \theta^n  \varphi^n \to T_{k,\varepsilon}'(\theta)\, \theta  \varphi \quad \mbox{as $n \to \infty$ a.e.~on $Q \times D$}.\]
Hence, by Fatou's lemma,  
\begin{align*}
\lim_{n \to \infty} \int_{Q \times D} \left( T_{k,\varepsilon}'(\theta^n)\, \theta^n  \varphi^n \nabq\, \cU_\eta \cdot \nabq\, \cU_e \right) \psi \dd \bq \dd \bx \dd t\\
\geq \int_{Q \times D} \left( T_{k,\varepsilon}'(\theta)\, \theta  \varphi \nabq\, \cU_\eta \cdot \nabq\, \cU_e \right) \psi \dd \bq \dd \bx \dd t.
\end{align*}

We are left to deal with the term $\mathcal{T}_4$. We perform partial integration with respect to the variable $\bq$ in the first term in the parentheses using the implicitly imposed homogeneous Dirichlet boundary condition $\varphi^n|_{Q \times \partial D}=0$ (cf.~\eqref{eq:bc0}). Hence, 
\begin{align*}
\mathcal{T}_4 &= \int_{Q \times D}  T_{k,\varepsilon}'(\theta^n)\left( |\nabq\, \cU_e|^2 -\theta^n \Delta_{\bq}\, \cU_e\right) \varphi^n \, \psi \dd \bq \dd \bx \dd t\\
& = \int_{(Q \times D)\cap \{\theta^n \leq k + \varepsilon\}}  T_{k,\varepsilon}'(\theta^n)\left( |\nabq\, \cU_e|^2 -\theta^n \Delta_{\bq}\, \cU_e\right) \varphi^n \, \psi \dd \bq \dd \bx \dd t\\
& = \int_{(Q \times \wtD)\cap \{\theta^n \leq k + \varepsilon\}}  T_{k,\varepsilon}'(\theta^n)\left( |\nabq\, \cU_e|^2 -\theta^n \Delta_{\bq}\, \cU_e\right) \varphi^n \, \psi \dd \bq \dd \bx \dd t\\
&\quad + \int_{(Q \times (D\setminus \wtD))\cap \{\theta^n \leq k + \varepsilon\}}  T_{k,\varepsilon}'(\theta^n)\left( |\nabq\, \cU_e|^2 -\theta^n \Delta_{\bq}\, \cU_e\right) \varphi^n \, \psi \dd \bq \dd \bx \dd t =: \mathcal{T}_{4,1} + \mathcal{T}_{4,2}, 
\end{align*}
for any $\wtD \subset\joinrel\subset D$. The passage to the limit $n \to \infty$ in $\mathcal{T}_{4,1}$ is easy because  $|\nabq\, \cU_e|^2$ and  $\Delta_{\bq}\, \cU_e$ are uniformly continuous (and therefore bounded) functions on $\overline{\wtD}$ and $\theta^n \leq k+\varepsilon$ over the region of integration $(Q \times \wtD)\cap \{\theta^n \leq k + \varepsilon\}$. Thus, by the dominated convergence theorem,
\begin{align*}
\lim_{n \to \infty} \int_{(Q \times \wtD)\cap \{\theta^n \leq k + \varepsilon\}}  T_{k,\varepsilon}'(\theta^n)\left( |\nabq\, \cU_e|^2 -\theta^n \Delta_{\bq}\, \cU_e\right) \varphi^n \, \psi \dd \bq \dd \bx \dd t \\
= \int_{(Q \times \wtD)\cap \{\theta \leq k + \varepsilon\}}  T_{k,\varepsilon}'(\theta)\left( |\nabq\, \cU_e|^2 -\theta\Delta_{\bq}\, \cU_e\right) \varphi \, \psi \dd \bq \dd \bx \dd t
\end{align*}
for any $\wtD \subset\joinrel\subset D$, and for all $k \geq 1$ and $\varepsilon \in (0,1)$.

\smallskip

Concerning $\mathcal{T}_{4,2}$, for $k \geq 1$ and $\varepsilon \in (0,1)$ fixed, we now choose a specific $\wtD \subset\joinrel\subset D$ such that 
\[ |\nabq\, \cU_e|^2 - (k+\varepsilon) \Delta_{\bq}\, \cU_e \geq 0\quad \mbox{on $D \setminus \wtD$}.\]
As in the argument leading to the bound \eqref{eq:up}, this is always possible to achieve thanks to \eqref{eq:5+}$_2$ taking, for example, $\wtD=(1-\delta)D$ with $\delta \in (0,1)$ sufficiently small; indeed because of the radial symmetry of $\cU_e(\bq)=U_e(\frac{1}{2}|\bq^2)$ the expression on the left-hand side of the last inequality is a function of $|\bq|$ alone and $\Delta_{\bq}\,\cU_e/|\nabq\,\cU_e|^2$ converges to $0$ as $|\bq| \to (\sqrt{b})_{-}$ thanks to \eqref{eq:5+}$_2$. Hence, by Fatou's lemma, we have
\begin{align*}
 \lim_{n \to \infty} \int_{(Q \times (D\setminus \wtD))\cap \{\theta^n \leq k + \varepsilon\}}  T_{k,\varepsilon}'(\theta^n)\left( |\nabq\, \cU_e|^2 -\theta^n \Delta_{\bq}\, \cU_e\right) \varphi^n \, \psi \dd \bq \dd \bx \dd t  \\
 \geq \int_{(Q \times (D\setminus \wtD))\cap \{\theta \leq k + \varepsilon\}}  T_{k,\varepsilon}'(\theta)\left( |\nabq\, \cU_e|^2 -\theta \Delta_{\bq}\, \cU_e\right) \varphi \, \psi \dd \bq \dd \bx \dd t.    
\end{align*}
Having dealt with both $\mathcal{T}_{4,1}$ and $\mathcal{T}_{4,2}$ we can pass to the limit in $\mathcal{T}_4$ to infer that 
\begin{align*}
&\lim_{n \to \infty}   \int_{Q \times D} T_{k,\varepsilon}'(\theta^n)\left( \theta^n \nabq \varphi^n \cdot \nabq\, \cU_e + |\nabq\, \cU_e|^2 \varphi^n\right)  \psi \dd \bq \dd \bx \dd t \\
&\qquad \geq \int_{Q \times D}
T_{k,\varepsilon}'(\theta)\left( \theta \nabq \varphi \cdot \nabq\, \cU_e + |\nabq\, \cU_e|^2 \varphi\right)  \psi \dd \bq \dd \bx \dd t.
\end{align*}

Having passed to the limit $n \to \infty$ in each of the terms $\mathcal{T}_1, \mathcal{T}_2, \mathcal{T}_3, \mathcal{T}_4$ that appear on the right-hand side in the weak formulation of the renormalized temperature equation \eqref{eq:ren} we can now pass to the limit $n \to \infty$ in the weak formulation of \eqref{eq:ren} to obtain a renormalized variational inequality involving the mollified cutoff function $T_{k,\varepsilon}$ and its first two derivatives, $T_{k,\varepsilon}'$ and $T_{k,\varepsilon}''$:
\begin{align*}
&\langle \partial_t T_{k,\varepsilon}(\theta) , \psi \rangle - \int_Q \bv\, T_{k,\varepsilon}(\theta) \cdot \nabx \psi \dd \bx \dd t + \int_Q \kappa(\theta) \nabx T_{k,\varepsilon}(\theta)\cdot \nabx \psi \dd \bx \dd t\\
& \quad \geq  \int_Q 2\nu(\theta) T_{k,\varepsilon}'(\theta)|\bDD(\bv)|^2\psi  \dd \bx \dd t \\
&\qquad +\int_Q \left(T_{k,\varepsilon}'(\theta) \bDD(\bv): \int_D \theta (\nabq\, \cU_\eta \otimes \bq) \varphi^n \dd \bq \right) \psi \dd \bx \dd t \\
&\qquad-\int_Q T_{k,\varepsilon}''(\theta)\kappa(\theta)|\nabx\theta|^2 \psi \dd \bx \dd t \\
&\qquad+ \int_Q \left(T_{k,\varepsilon}'(\theta) \int_D \theta\nabq \varphi \cdot \nabq\, \cU_e + \theta \varphi \nabq\, \cU_\eta \cdot \nabq\, \cU_e + |\nabq\, \cU_e|^2 \varphi \dd \bq\right) \psi \dd \bx \dd t,
\end{align*}
which holds for all test functions $\psi \in C([0,T]; W^{1,z}(\Omega))$ with $z>d$ such that $\psi \geq 0$ on $\overline{Q}$, and for all $k \geq 1$ and all $\varepsilon \in (0,1)$.

Alternatively, testing \eqref{eq:ren} with $\psi \in C^1([0,T]; W^{1,z}(\Omega))$ where $z>d$, and where $\psi \geq 0$ on $\overline{Q}$ and $\psi(T,\cdot)=0$, integrating by parts with respect to $t$ in the time derivative term prior to passing to the limit $n\to \infty$, and then letting $n \to \infty$, we have that
\begin{align*}
&-\int_Q T_{k,\varepsilon}(\theta) \partial_t \psi \dd \bx \dd t - \int_Q \bv\, T_{k,\varepsilon}(\theta) \cdot \nabx \psi \dd \bx \dd t + \int_Q \kappa(\theta) \nabx T_{k,\varepsilon}(\theta)\cdot \nabx \psi \dd \bx \dd t\\
& \quad \geq \int_\Omega T_{k,\varepsilon}(\theta_0)\psi(0) \dd \bx   + \int_Q 2\nu(\theta) T_{k,\varepsilon}'(\theta)|\bDD(\bv)|^2\psi  \dd \bx \dd t \\
&\qquad +\int_Q \left(T_{k,\varepsilon}'(\theta) \bDD(\bv): \int_D \theta (\nabq\, \cU_\eta \otimes \bq) \varphi \dd \bq \right) \psi \dd \bx \dd t \\
&\qquad-\int_Q T_{k,\varepsilon}''(\theta)\kappa(\theta)|\nabx\theta|^2 \psi \dd \bx \dd t \\
&\qquad+ \int_Q \left(T_{k,\varepsilon}'(\theta) \int_D \theta\nabq \varphi \cdot \nabq\, \cU_e + \theta \varphi \nabq\, \cU_\eta \cdot \nabq\, \cU_e + |\nabq\, \cU_e|^2 \varphi \dd \bq\right) \psi \dd \bx \dd t,
\end{align*}
for all test functions $\psi \in C^1_0([0,T); W^{1,z}(\Omega))$ with $z>d$ such that $\psi \geq 0$ on $\overline{Q}$, and for all $k \geq 1$ and all $\varepsilon \in (0,1)$.

It remains to remove the $\varepsilon$-mollification from $T_{k,\varepsilon}$ by passing to the limit $\varepsilon \to 0_+$. As $T_{k,\varepsilon}''(\theta) \leq 0$ on $Q$ thanks to the positivity of $\theta$, trivially
\[ \nabx T_{k,\varepsilon}(\theta) = T'_{k,\varepsilon}(\theta) \nabx \theta,\]
and
\[T_{k,\varepsilon}(\theta) \to T_k(\theta)\quad \mbox{and} \quad  T'_{k,\varepsilon}(\theta) \to T'_k(\theta)\]
strongly in $L^p(Q)$ for $p \in [1,\infty)$ and weakly-$\ast$ in $L^\infty(Q)$,
we can drop the penultimate term from the right-hand side and then pass to the limit $\varepsilon \to 0_{+}$ to deduce that 
\begin{align}\label{eq:ref}
\begin{aligned}
&-\int_Q T_{k}(\theta) \partial_t \psi \dd \bx \dd t - \int_Q \bv\, T_{k}(\theta) \cdot \nabx \psi \dd \bx \dd t + \int_Q \kappa(\theta) \nabx T_{k}(\theta)\cdot \nabx \psi \dd \bx \dd t\\
& \quad \geq \int_\Omega T_{k}(\theta_0)\psi(0) \dd \bx   + \int_Q 2\nu(\theta) T_{k}'(\theta)|\bDD(\bv)|^2\psi  \dd \bx \dd t\\
&\qquad +\int_Q \left(T_{k}'(\theta) \bDD(\bv): \int_D \theta (\nabq\, \cU_\eta \otimes \bq) \varphi \dd \bq \right) \psi \dd \bx \dd t \\
&\qquad+ \int_Q \left(T_{k}'(\theta) \int_D \theta\nabq \varphi \cdot \nabq\, \cU_e + \theta \varphi \nabq\, \cU_\eta \cdot \nabq\, \cU_e + |\nabq\, \cU_e|^2 \varphi \dd \bq\right) \psi \dd \bx \dd t,
\end{aligned}
\end{align}
for all test functions $\psi \in C^1_0([0,T); W^{1,z}(\Omega))$ with $z>d$ such that $\psi \geq 0$ on $\overline{Q}$, and for all $k \geq 1$. This is the renormalized variational inequality satisfied by $\theta$ we set out to derive.

\begin{remark}\label{eq:rem4} 
In the absence of additional information concerning the functions $\theta$ and $\varphi$ it does not seem possible to pass to the limit $k \rightarrow \infty$ to remove the cutoff function $T_k$ in \eqref{eq:ref}. In particular, we do not have sufficiently strong a priori bounds at our disposal to ensure that 
\[ \mathcal{I}:=\int_{Q \times D} \left[\theta\nabq \varphi \cdot \nabq\, \cU_e + \theta \varphi \nabq\, \cU_\eta \cdot \nabq\, \cU_e + |\nabq\, \cU_e|^2 \varphi \right]  \dd \bq \dd \bx \dd t< \infty. \]
Looking at this integral, we see that 
\begin{align}\label{eq:bdi}
\begin{aligned} 
\mathcal{I}&= \int_{Q \times D} \frac{1}{\sqrt{\theta \varphi}}\left(\theta\nabq \varphi  + \theta \varphi \nabq\, \cU_\eta  + \nabq\, \cU_e \varphi \right) \cdot \sqrt{\theta \varphi} \,\nabq\, \cU_e  \dd \bq \dd \bx \dd t \\
 &\leq \left(\int_{Q \times D}\frac{1}{\theta \varphi} |\theta\nabq \varphi  + \theta \varphi \nabq\, \cU_\eta  + \nabq\, \cU_e \varphi |^2 \dd \bq \dd \bx \dd t 
 \right)^{\frac{1}{2}}\left(\int_{Q \times D} \theta \varphi |\nabq\,\cU_e|^2  \dd \bq \dd \bx \dd t \right)^{\frac{1}{2}}\\
 & \leq C \left(\int_{Q \times D} \theta \varphi |\nabq\,\cU_e|^2  \dd \bq \dd \bx \dd t \right)^{\frac{1}{2}},
\end{aligned}
\end{align}
where in the transition from the second line to the third line we have used \eqref{eq:18} to bound the first factor. Unfortunately, because no pointwise upper bound on the absolute temperature $\theta$ is available, we have no means of ensuring that the integral appearing in the last line is finite. However, if, in addition to the uniform positive lower bound on $\theta$ implied by the assumed positivity of the initial temperature $\theta_0$, one \underline{assumes} that $\theta$ is (as is physically reasonable) also bounded from above, i.e., that 
\begin{align}\label{eq:con+}
 \|\theta\|_{L^\infty(Q)} < \infty,
\end{align}
then, by revisiting Step 15 in Section \ref{sec:2} under this added assumption, the argument presented there implies (cf., in particular, \eqref{eq:up} with $Q_{\theta_{\mathrm{max}}}$ replaced by $Q$) that the integral in the last line of \eqref{eq:bdi} is bounded by a constant. 
Thus, by using test functions $\psi\in C^1_0([0,T);W^{1,\infty}(\Omega))$ such that $\psi \geq 0$ on $\overline{Q}$, one can pass to the limit $k \rightarrow \infty$ (using the dominated convergence theorem, for example,) to completely remove the cutoff function and deduce the \underline{conditional} variational inequality
\begin{align}\label{eq:con1}
\begin{aligned}
&-\int_Q \theta \partial_t \psi \dd \bx \dd t - \int_Q \bv\, \theta \cdot \nabx \psi \dd \bx \dd t + \int_Q \kappa(\theta) \nabx \theta \cdot \nabx \psi \dd \bx \dd t\\
& \quad \geq \int_\Omega \theta_0\psi(0) \dd \bx + \int_Q 2\nu(\theta) |\bDD(\bv)|^2\psi  \dd \bx \dd t 
+\int_Q \left( \bDD(\bv): \int_D \theta (\nabq\, \cU_\eta \otimes \bq) \varphi \dd \bq \right) \psi \dd \bx \dd t \\
&\qquad+ \int_Q \left(\int_D \theta\nabq \varphi \cdot \nabq\, \cU_e + \theta \varphi \nabq\, \cU_\eta \cdot \nabq\, \cU_e + |\nabq\, \cU_e|^2 \varphi \dd \bq\right) \psi \dd \bx \dd t,
\end{aligned}
\end{align}
for all $\psi\in C^1_0([0,T);W^{1,\infty}(\Omega))$ such that $\psi \geq 0$ of $\overline{Q}$, the condition being that the assumption stated in \eqref{eq:con+} holds.
\end{remark}

\begin{remark}
A natural question to pose is why we chose to work with a renormalization of the temperature equation instead of the evolution equation \eqref{eq:entr} for the entropy $\eta$. The reason is that the weak formulation of the equation \eqref{eq:entr} would involve a term of the form 
\[ -\int_Q \bj_\eta \cdot \nabx z \dd \bx \dd t, \]
where $z \in C^1_0([0,T);W^{1,\infty}(\Omega))$ is a test function and $\bj_\eta$ is the  entropy flux, defined  in \eqref{eq:entf}. Unfortunately, the second summand of $\bj_\eta$, that is 
\[ \theta \nabx \int_D [\cU_\eta \varphi +  \varphi \log \varphi] \dd  \bq,\]
cannot be generally guaranteed to belong to $L^1(Q)$, even if the function $\theta$ were replaced in this term by its truncation $T_k(\theta) \in L^\infty(Q)$. Indeed, while the energy inequality 
\eqref{eq:energy+} implies that the function $\varphi \log \varphi$ is an element of the function space $L^\infty(0,T;L^1(\Omega \times D))$, we have no means of ensuring that
\[ \nabx \int_D \varphi \log \varphi \dd \bq\]
belongs to $L^1(Q)$. In contrast, such a difficulty does not occur in the weak formulation of the renormalized temperature equation \eqref{eq:ren}, where all the terms involved in the weak formulation are well-defined. 
\end{remark}


\subsection{Limit passage in the energy inequality}
Our objective here is to pass to the limit $n \to \infty$ in the sequence of inequalities
\begin{align}
\begin{aligned}\label{eq:en-n}
&\int_\Omega \bigg[\frac{1}{2}|\bv^n(t)|^2 +  \mathcal{H}(\theta^n(t)) \bigg] \dd \bx  + \int_{\Omega \times D} [\,\cU_e \varphi^n(t) + \cU_{\eta} \varphi^n(t) + \mathcal{F}(\varphi^n(t))]\dd \bq \dd \bx \\
&\qquad+ \int_0^t \! \int_\Omega \bigg[\frac{2\nu(\theta^n)|\bDD(\bv^n)|^2}{\theta^n} + \frac{\kappa(\theta^n) |\nabx \theta^n|^2}{(\theta^n)^2}\bigg]\dd \bx \dd s\\ 
&\qquad + \int_0^t \int_\Omega \bigg[\theta^n \int_D \frac{|\nabx \varphi^n|^2}{\varphi^n}\dd \bq + \int_D \frac{4}{\theta^n \varphi^n}\left| \theta^n \nabq\varphi^n + \theta^n \varphi^n \nabq\,\cU_\eta + \varphi^n \nabq \,\cU_e\right|^2\! \dd \bq\bigg] \dd \bx \dd s\\
& \leq  \int_\Omega \bigg[\frac{1}{2}|\bv_0^n|^2 +  \mathcal{H}(\theta_0^n) \bigg] \dd \bx  + \int_{\Omega \times D} [\,\cU_e \varphi_0^n + \cU_{\eta} \varphi_0^n + \mathcal{F}(\varphi_0^n)]\dd \bq \dd \bx \\
&\qquad + \int_0^t \int_\Omega \bff(s) \cdot \bv^n(s)\dd \bx \dd s,
\end{aligned}
\end{align}
where $\mathcal{H}(s):=s - 1 - \log s$ and $\mathcal{F}(s):= s (\log s - 1) + 1$ for $s > 0$, and $\mathcal{F}(0):=1$. We shall do so using a standard argument. We multiply \eqref{eq:en-n} by a test function $\phi \in C^\infty_0((0,T); \mathbb{R}_{\geq 0})$ and integrate over $t \in (0,T)$. Hence, 
\begingroup
\allowdisplaybreaks
\begin{align}
\label{eq:en-n1}
&\int_0^T \left\{\int_\Omega \bigg[\frac{1}{2}|\bv^n(t)|^2 +  \mathcal{H}(\theta^n(t)) \bigg] \dd \bx  + \int_{\Omega \times D} [\,\cU_e \varphi^n(t) + \cU_{\eta} \varphi^n(t) + \mathcal{F}(\varphi^n(t))]\dd \bq \dd \bx \right. \nonumber\\
&\qquad+ \int_0^t \! \int_\Omega \bigg[\frac{2\nu(\theta^n)|\bDD(\bv^n)|^2}{\theta^n} + \frac{\kappa(\theta^n) |\nabx \theta^n|^2}{(\theta^n)^2}\bigg]\dd \bx \dd s
\nonumber \\ 
&\qquad \left. + \int_0^t \int_\Omega \bigg[\theta^n \int_D \frac{|\nabx \varphi^n|^2}{\varphi^n}\dd \bq + \int_D \frac{4}{\theta^n \varphi^n}\left| \theta^n \nabq\varphi^n + \theta^n \varphi^n \nabq\,\cU_\eta + \varphi^n \nabq \,\cU_e\right|^2\! \dd \bq\bigg] \dd \bx \dd s\right\}\phi(t) \dd t \nonumber\\
& \leq  \int_0^T \bigg \{\int_\Omega \bigg[\frac{1}{2}|\bv_0^n|^2 +  \mathcal{H}(\theta_0^n) \bigg] \dd \bx  + \int_{\Omega \times D} [\,\cU_e \varphi_0^n + \cU_{\eta} \varphi_0^n + \mathcal{F}(\varphi_0^n)]\dd \bq \dd \bx \nonumber\\
& \qquad + \int_0^t \int_\Omega \bff(s) \cdot \bv^n(s)\dd \bx \dd s \bigg\} \phi(t) \dd t \qquad \forall\, \phi \in C^\infty_0((0,T);\mathbb{R}_{\geq 0}).
\end{align}
\endgroup
We now pass to the lim.inf as $n \to \infty$ term by term in \eqref{eq:en-n1}.

For the first term in the first pair of square brackets in the first line of \eqref{eq:en-n1}, we use \eqref{eq:va2} to pass to the limit. For the second term in the first pair of square brackets in the first line, we make use of \eqref{eq:thp} in conjunction with Fatou's lemma. For each of the three terms in the second pair of square brackets in the first line, we recall \eqref{eq:33} and again employ Fatou's lemma.

For the terms in the second line of \eqref{eq:en-n1} note that, by \eqref{eq:thp}, $\nu(\theta^n) \to \nu(\theta)$ and $\kappa(\theta^n) \to \kappa(\theta)$ a.e.~on $Q \times D$. By \eqref{eq:34}, $\nu(\theta^n) \bDD(\bv^n) \rightharpoonup \nu(\theta) \bDD(\bv)$ weakly in $L^{\frac{4+6\beta}{5+3\beta}}(Q;\mathbb{R}^{d \times d})$ with $\beta>\frac{5}{6}$. By the dominated convergence theorem and the uniform positive lower bound on the sequence $\theta^n$ over the set $Q$ coming from the minimum principle, we have $1/\sqrt{\nu(\theta^n) \,\theta^n}\to 1/\sqrt{\nu(\theta) \,\theta}$ strongly in $L^1(Q)$, and because this sequence is bounded in $L^\infty(Q)$, we infer its strong convergence in $L^s(Q)$ for all $s \in [1,\infty)$, and therefore, in particular, in 
$L^{\left(\frac{4+6\beta}{5+3\beta}\right)'}(Q)$. Therefore, $\sqrt{\nu(\theta^n)} \bDD(\bv^n)/\sqrt{\theta^n} \rightharpoonup\sqrt{\nu(\theta)} \bDD(\bv)/\sqrt{\theta}$ weakly in $L^1(Q;\mathbb{R}^{d\times d})$. However, the energy inequality \eqref{eq:en-n1} guarantees the boundedness of this sequence in $L^2(Q;\mathbb{R}^{d \times d})$, and therefore the weak convergence of a subsequence, not indicated, in $L^2(Q;\mathbb{R}^{d \times d})$. By the uniqueness of the weak limit, it then follows that 
$\sqrt{\nu(\theta^n)} \bDD(\bv^n)/\sqrt{\theta^n} \rightharpoonup \sqrt{\nu(\theta)} \bDD(\bv)/\sqrt{\theta}$ weakly in $L^2(Q;\mathbb{R}^{d\times d})$. Using this in conjunction with weak lower-semicontinuity in the $L^2(Q;\mathbb{R}^{d \times d})$ norm and Fatou's lemma we can pass to the lim.inf as $n \to \infty$ in the first term in the square brackets in the second line.\footnote{The argument proceeds as follows: 
\begin{align*}\mbox{lim.inf}_{n \to \infty} \int_0^T \bigg[\int_0^t \! \int_\Omega \frac{2\nu(\theta^n)|\bDD(\bv^n)|^2}{\theta^n}\dd \bx \dd s \bigg]\phi(t) \dd t \geq \int_0^T \mbox{lim.inf}_{n \to \infty} \bigg[\int_0^t \! \int_\Omega \frac{2\nu(\theta^n)|\bDD(\bv^n)|^2}{\theta^n}\dd \bx \dd s \bigg]\phi(t) \dd t\\ 
\geq \int_0^T  \bigg[\int_0^t \! \int_\Omega \frac{2\nu(\theta)|\bDD(\bv)|^2}{\theta}\dd \bx \dd s \bigg]\phi(t) \dd t\quad \forall\, \phi \in C^\infty_0((0,T);\mathbb{R}_{\geq 0}), 
\end{align*}
where the first inequality follows by Fatou's lemma and the second inequality follows by weak lower-semicontinunity of the $L^2(Q_t)$ norm for all $t \in (0,T)$, where $Q_t = (0,t) \times \Omega$.} 

Concerning the second term in the square brackets in the second line of \eqref{eq:en-n1}, by \eqref{eq:logfi} we have $\nabx \log \theta^n \rightharpoonup \nabx \log \theta$ weakly in $L^2(Q;\mathbb{R}^d)$ as $n \to \infty$. We claim that we also have 
strong convergence of $\sqrt{\kappa(\theta^n)} \to \sqrt{\kappa(\theta)}$ in $L^2(Q)$, and therefore $\sqrt{\kappa(\theta^n)}\nabx \log \theta^n \rightharpoonup \sqrt{\kappa(\theta)}\nabx \log \theta$ weakly in $L^1(Q;\mathbb{R}^d)$ as $n \to \infty$; thus, by the upper bound on the sequence coming from the energy inequality and uniqueness of the weak limit, we also have $\sqrt{\kappa(\theta^n)}\nabx \log \theta^n \rightharpoonup \sqrt{\kappa(\theta)}\nabx \log \theta$ weakly in $L^2(Q;\mathbb{R}^d)$ as $n \to \infty$. Hence, weak lower-semicontinuity of the  $L^2(Q;\mathbb{R}^d)$ norm implies the desired lower bound on this term when passing to the lim.inf as $n \to \infty$ in \eqref{eq:en-n1}. To prove that $\sqrt{\kappa(\theta^n)} \to \sqrt{\kappa(\theta)}$ in $L^2(Q)$ as $n \to \infty$, note that $|\sqrt{\kappa(\theta^n)} - \sqrt{\kappa(\theta)}|^2 \leq |\kappa(\theta^n) - \kappa(\theta)|$; it therefore suffices to prove that $\kappa(\theta^n) \to \kappa(\theta)$ in $L^1(Q)$ as $n \to \infty$. To this end,
recall that the strong convergence of $\theta^n$ to $\theta$ in $L^1(Q)$ stated in \eqref{eq:str} implies almost everywhere convergence of a subsequence (not indicated) (cf. \eqref{eq:thp}). Thanks to the assumed continuity of $\kappa$, it then follows that $\kappa(\theta^n) \to \kappa(\theta)$ a.e.~on $Q$ as $n \to \infty$. In addition, 
thanks to the upper bound in \eqref{eq:9}, with $\beta>5/6$, and  \eqref{eq:22}, which guarantees boundedness of the sequence $(\theta^n)_{n \geq 1}$ in $L^{\frac{2}{3} + \beta}(Q)$, it follows with 
$s=1 + \frac{2}{3\beta}$ that 
\[\int_{Q}[\kappa(\theta^n)]^s \dd \bx \dd t\leq C\int_{Q}(1 + (\theta^n)^\beta)^s \dd \bx \dd t \leq C\big(1 + \int_Q (\theta^n)^{\frac{2}{3}+ \beta} \dd \bx \dd t\big) \leq C.\]
Because $s>1$, this bound guarantees the uniform integrability of the sequence $(\kappa(\theta^n))_{n \geq 1}$ in $L^1(Q)$, which then, together with the a.e.~convergence of $\kappa(\theta^n) \to \kappa(\theta)$ on $Q$
as $n \to\infty$, implies the strong convergence of $\kappa(\theta^n) \to \kappa(\theta)$ in $L^1(Q)$ by Vitali's theorem. The desired strong convergence of $\sqrt{\kappa(\theta^n)} \to \sqrt{\kappa(\theta)}$ in $L^2(Q)$ as $n \to \infty$ then follows directly. 

\smallskip

For the first term in the square brackets in the third line of \eqref{eq:en-n1}, we use \eqref{eq:tf} in conjunction with weak lower-semicontinuity of the $L^2(Q \times D;\mathbb{R}^d)$ norm and Fatou's lemma, while for the second term in the third line we use \eqref{eq:tfj} and the definition \eqref{eq:bnj} of $\boldsymbol{j}^n_\varphi$ together with weak lower-semicontinuity of the $L^2(Q \times D;\mathbb{R}^d)$ norm and Fatou's lemma.

On the right-hand side, we use the strong convergence of the initial data and the weak convergence of $\bv^n$ to $\bv$ in $L^2(0,T;L^2(\Omega;\mathbb{R}^d))$. This yields
\begin{align*}
&\int_0^T \left\{\int_\Omega \bigg[\frac{1}{2}|\bv(t)|^2 +  \mathcal{H}(\theta(t)) \bigg] \dd \bx  + \int_{\Omega \times D} [\,\cU_e \varphi(t) + \cU_{\eta} \varphi(t) + \mathcal{F}(\varphi(t))]\dd \bq \dd \bx \right.\\
&\qquad+ \int_0^t \! \int_\Omega \bigg[\frac{2\nu(\theta)|\bDD(\bv)|^2}{\theta} + \frac{\kappa(\theta) |\nabx \theta|^2}{(\theta)^2}\bigg]\dd \bx \dd s\\ 
&\qquad \left. + \int_0^t \int_\Omega \bigg[\theta \int_D \frac{|\nabx \varphi|^2}{\varphi}\dd \bq + \int_D \frac{4}{\theta \varphi}\left| \theta \nabq\varphi + \theta \varphi \nabq\,\cU_\eta + \varphi \nabq \,\cU_e\right|^2\! \dd \bq\bigg] \dd \bx \dd s\right\}\phi(t) \dd t\\
& \leq  \int_0^T \bigg \{\int_\Omega \bigg[\frac{1}{2}|\bv_0|^2 +  \mathcal{H}(\theta_0) \bigg] \dd \bx  + \int_{\Omega \times D} [\,\cU_e \varphi_0 + \cU_{\eta} \varphi_0 + \mathcal{F}(\varphi_0)]\dd \bq \dd \bx \\
& \qquad + \int_0^t \int_\Omega \bff(s) \cdot \bv(s)\dd \bx \dd s \bigg\} \phi(t) \dd t \qquad \forall\, \phi \in C^\infty_0((0,T);\mathbb{R}_{\geq 0}).
\end{align*}
This then implies the desired energy inequality
\begin{align*}
&\int_\Omega \bigg[\frac{1}{2}|\bv(t)|^2 +  \mathcal{H}(\theta(t)) \bigg] \dd \bx  + \int_{\Omega \times D} [\,\cU_e \varphi(t) + \cU_{\eta} \varphi(t) + \mathcal{F}(\varphi(t))]\dd \bq \dd \bx \\
&\qquad+ \int_0^t \! \int_\Omega \bigg[\frac{2\nu(\theta)|\bDD(\bv)|^2}{\theta} + \frac{\kappa(\theta) |\nabx \theta|^2}{\theta^2}\bigg]\dd \bx \dd s\\ 
&\qquad + \int_0^t \int_\Omega \bigg[\theta \int_D \frac{|\nabx \varphi|^2}{\varphi}\dd \bq + \int_D \frac{4}{\theta \varphi}\left| \theta \nabq\varphi + \theta \varphi \nabq\,\cU_\eta + \varphi \nabq \,\cU_e\right|^2\! \dd \bq\bigg] \dd \bx \dd s\\
& \leq  \int_\Omega \bigg[\frac{1}{2}|\bv_0|^2 +  \mathcal{H}(\theta_0) \bigg] \dd \bx  + \int_{\Omega \times D} [\,\cU_e \varphi_0 + \cU_{\eta} \varphi_0 + \mathcal{F}(\varphi_0)]\dd \bq \dd \bx + \int_0^t \int_\Omega \bff(s) \cdot \bv(s)\dd \bx \dd s.
\end{align*}
for almost every $t \in (0,T)$. That completes the proof of Theorem \ref{thm:1}. \hfill $\Box$

\section{Conclusions}\label{sec:4}
Our aim in this paper was to prove the weak sequential stability of solutions to a thermodynamically consistent kinetic model of nonisothermal flow of a dilute polymeric fluid, which involves a system of nonlinear partial differential equations coupling the unsteady incompressible temperature-dependent Navier--Stokes equations to a temperature-dependent generalization of the classical Fokker--Planck equation satisfied by a probability density function $\varphi$ modelling the random configuration of polymer molecules in the viscous Newtonian solvent, and an evolution equation for the absolute temperature $\theta$. We showed that sequences of smooth solutions to the initial-boundary-value problem for the model converge to a global-in-time large-data weak solution that satisfies an energy inequality, where the limiting absolute temperature satisfies a renormalized variational inequality. In the absence of a positive pointwise upper bound on the function $\theta$, it generally does not seem possible to remove the renormalization. We have also shown that if $\theta$ is assumed to be bounded above by a positive constant over the space-time domain $Q$, then the renormalization can be eliminated from the variational inequality satisfied by $\theta$.

\bibliographystyle{abbrv}
\bibliography{references}

\vfill\eject

\end{document}